\documentclass[a4paper,11pt,french,leqno]{smfart}

\usepackage{amsmath}
\usepackage{amssymb}
\usepackage{amsthm}
\usepackage{amsopn}
\usepackage{amscd}
\usepackage{mathdots}
\usepackage{rotating}

\usepackage{stmaryrd}

\usepackage{fullpage}
\usepackage[utf8]{inputenc}
 \usepackage[T1]{fontenc}
\usepackage{graphics, xcolor}
\usepackage{textcomp} 
\usepackage[all]{xy}
\usepackage{lscape}
\usepackage{url}
\usepackage{verbatim}
\usepackage{mathrsfs}
  \usepackage[francais]{babel}

\def\commentaire#1{\ifx\cachercommentaires\undefined \textcolor{red}{#1}\else \fi } 

\newcommand{\scrA}{\mathscr{A}}
\newcommand{\scrB}{\mathscr{B}}

\newcommand{\scrQ}{\mathscr{Q}}

\newcommand{\scrX}{\mathscr{X}}

\newcommand{\frS}{\mathfrak{S}}

\newcommand{\frX}{\mathfrak{X}}

\newcommand{\fra}{\mathfrak{a}}

\newcommand{\frd}{\mathfrak{d}}

\newcommand{\frg}{\mathfrak{g}}
\newcommand{\frh}{\mathfrak{h}}

\newcommand{\frk}{\mathfrak{k}}
\newcommand{\frl}{\mathfrak{l}}

\newcommand{\fro}{\mathfrak{o}}
\newcommand{\frp}{\mathfrak{p}}
\newcommand{\frqqq}{\mathfrak{q}}

\newcommand{\frs}{\mathfrak{s}}
\newcommand{\frt}{\mathfrak{t}}
\newcommand{\fru}{\mathfrak{u}}

\newcommand{\bbC}{\mathbb{C}}

\newcommand{\bbR}{\mathbb{R}}

\newcommand{\bbZ}{\mathbb{Z}}

\newcommand{\caC}{\mathcal{C}}

\newcommand{\caE}{\mathcal{E}}

\newcommand{\caG}{\mathcal{G}}

\newcommand{\caL}{\mathcal{L}}
\newcommand{\caM}{\mathcal{M}}

\newcommand{\caP}{\mathcal{P}}

\newcommand{\caR}{\mathcal{R}}

\newcommand{\caX}{\mathcal{X}}

\def\U{\mathbf{U}}

\newcommand{\GL}{\mathbf{GL}}

\newcommand{\Ad}{\mathrm{Ad}}

\newcommand{\Hom}{\mathrm{Hom}}

\newcommand{\Ind}{\mathrm{Ind}}

\newcommand{\tr}{\mathrm{Tr}\, }
\newcommand{\Id}{\mathrm{Id}}

\renewcommand{\Ind}{\mathrm{Ind}}

\newcommand{\SL}{\mathbf{SL}}
\newcommand{\SU}{\mathbf{SU}}
\newcommand{\SO}{\mathbf{SO}}
\newcommand{\Sp}{\mathbf{Sp}}
\newcommand{\Or}{\mathbf{O}}

\newcommand{\bil}[2]{\langle  #1,#2 \rangle }

\newcommand{\Speh}{\mathbf{Speh}}
\newcommand{\sgn}{\mathbf{sgn}}
\newcommand{\Triv}{\mathbf{Triv}}

\newcommand{\Std}{\mathbf{Std}}

\newcommand{\rvline}{\hspace*{-\arraycolsep}\vline\hspace*{-\arraycolsep}}

\newcommand{\GSp}{\mathbf{GSp}}

\theoremstyle{plain}
\newtheorem{thm}{Théorème}[section]
\newtheorem{lemme}[thm]{Lemme}
\newtheorem{cor}[thm]{Corollaire}
\newtheorem{prop}[thm]{Proposition}

\theoremstyle{definition}
 \newtheorem{defi}[thm]{Définition}
 
\newtheorem{rmqs}[thm]{Remarques}
\newtheorem{rmq}[thm]{Remarque}
\newtheorem{exemple}[thm]{Exemple}

\def \dem {\noindent \underline{\sl Démonstration}. }

\begin{document}

\numberwithin{equation}{section}

\title{Séries discrètes des espaces symétriques et paquets d'Arthur}
  \author{Colette Moeglin}
 \address{CNRS, Institut Mathématique de Jussieu } 
 \email{colette.moeglin@imj-prg.fr}

  \author{David Renard  }
 \address{Centre de Mathématiques
 Laurent Schwartz,  Ecole Polytechnique} 
\email{david.renard@polytechnique.edu}

\date{\today}

\begin{abstract} On vérifie certaines conjectures de Sakellaridis et Venkatesh donnant  une description du spectre discret d'une variété 
sphérique $\scrX=G/H$ dans le formalisme d'Arthur-Langlands lorsque   $G$ est un groupe classique réel
 et $\scrX$ est un espace symétrique. On détermine  ensuite explicitement les membres des paquets d'Arthur  concernés contribuant 
au spectre discret, ce qui fait apparaître des résultats de multiplicité un.

 \medskip 
 
\noindent {\bf  Abstract.} We check  Sakellaridis-Venkatesh  conjectures giving a description of  the discrete spectrum of a spherical variety $\scrX=G/H$
in the Langlands-Arthur formalism when $G$ is a classical real group and $\scrX$  is a symmetric space.
Then, we compute explicitly the representations in the relevant   Arthur paquets which appear in the discrete spectrum, and we establish 
some multiplicity one results.
 \end{abstract} 

\maketitle

\section{Introduction}
Le but de cet article est de vérifier certaines  conjectures de Sakellaridis et Venkatesh énoncées dans  \cite{SV}. Ces conjectures
proposent une description  du spectre discret  d'une variété sphérique  $\scrX$ définie sur un corps local en termes de paramètres
d'Arthur-Langlands. En particulier, ils introduisent  la notion de  $L$-groupe ${}^LG_\scrX$   associé à   $\scrX$, dont la 
 construction est poursuivie  dans \cite{KS}.

Le cadre de notre étude est beaucoup plus restreint que celui de {\sl op. cit.} :
les variétés sphériques que nous considérons sont des espaces symétriques $\scrX=G/H$, où $G$ est le
groupe des points réels d'un groupe algébrique connexe réductif $\mathbf G$  défini sur $\bbR$,
  où $\mathbf H=\mathbf G^\sigma$ est le groupe des points fixes d'une involution  
 $\sigma$  de $\mathbf G$ définie sur $\bbR$, et $H$ est un sous-groupe de $\mathbf H(\bbR)$ contenant la composante neutre de celui-ci
  au sens des groupes de Lie, $\mathbf H(\bbR)_e$. En général, nous prendrons $H=\mathbf H(\bbR)$, mais à ce sujet, remarquons
  que si les conjectures sont valides pour un $\scrX=G/H$ avec $H$  soumis à $ \mathbf H(\bbR)_e\subset  H\subset  \mathbf H(\bbR)$, elles seront 
  valides pour  $\scrX=G/H'$ avec $H\subset H'\subset \mathbf H(\bbR) $, et dans certains cas, nous démontrons les conjectures 
  pour un groupe $H$ strictement
  plus petit que $H(\bbR)$. Remarquons aussi que les conjectures doivent être reformulées pour des groupes $H$  trop petits. 
   Ce phénomène est à mettre en relation avec 
  une extension possible de l'énoncé des conjectures à une situation tordue par un caractère unitaire  $\chi$ de $H/H_e$ sur laquelle nous reviendrons ci-dessous.
   De plus, nous supposons que $G$ est un groupe général linéaire, un groupe unitaire, un groupe symplectique
 ou bien un groupe spécial orthogonal (nous dirons simplement que $G$ est un groupe classique).
  
 On note $L^2_d(\scrX)$ la somme des sous-représentations (au sens des $(\mathfrak{g},K)$-modules) irréductibles de $L^2(\scrX)$ 
 et une telle sous-représentation est  appelée  dans la suite {\sl série discrète de $\scrX$}. 
 On suppose donc que ce spectre discret est non trivial, et l'on veut caractériser les 
 séries discrètes de $\scrX$ du point de vue des $L$-groupes et des paramètres d'Arthur.
 On adopte les notations habituelles pour les $L$-groupes de Langlands : $W_\bbR$ est le groupe de Weil de $\bbR$,
 $G{}^\vee$ est le groupe complexe dual  de $\mathbf G(\bbC)$ et $  {}^LG$ est un certain  produit semi-direct $G{}^\vee\rtimes W_\bbR$.
 Il s'agit tout d'abord d'attacher à l'espace symétrique $\scrX$ un $L$-groupe ${}^LG_\scrX$ et un morphisme
   \begin{equation} \label{varphiscrX} \varphi_\scrX:\; {}^LG_\scrX\times \SL(2,\bbC) \longrightarrow {}^LG . \end{equation}

      La conjecture affirme  que si $\pi$  est une série discrète de $\scrX$, alors il existe un paramètre de Langlands discret
 $\phi_\pi : W_\bbR \longrightarrow {}^LG_\scrX$,  
 tel que $\pi$ appartienne au paquet  d'Arthur  de paramètre
 $\psi_\pi=\varphi_\scrX \circ \tilde \phi_\pi :\;  W_\bbR\times \SL(2,\bbC)\longrightarrow {}^LG$,  
 où  $\phi_\pi $ a été étendu en  un morphisme 
 \[\tilde \phi_\pi : W_\bbR\times \SL(2,\bbC) \rightarrow {}^LG_\scrX\times \SL(2,\bbC) \]
 par l'identité sur le  facteur $\SL(2,\bbC)$.
  Résumons ceci par le diagramme commutatif : 
 \begin{equation}\label{factpar} \xymatrix{  W_\bbR\times \SL(2,\bbC) \ar[rr]^{\tilde \phi_\pi} \ar[rrd]_{\psi_\pi}&    
 & {}^LG_\scrX\times \SL(2,\bbC)\ar[d]^{\varphi_\scrX} \\
 && {}^LG }   \end{equation}

On vérifie donc que cet énoncé est vrai dans le   cadre indiqué (espaces symétriques des groupes classiques réels).   
Une bonne partie de la construction du $L$-groupe   ${}^LG_\scrX$  et du morphisme $\varphi_\scrX$ est déjà effectuée 
   dans \cite{SV} et \cite{KS}. En ce qui concerne   ${}^LG_\scrX$, Knop et Schalke déterminent sa composante neutre,
    et  dans tous les cas, sauf un sur lequel on reviendra ci-dessous, nous trouvons  que ${}^LG_\scrX$ est bien  un  $L$-groupe ayant cette composante neutre. 
Il est alors déterminé par le fait supplémentaire que c'est le   $L$-groupe  d'un groupe quasi-déployé ayant de la série discrète    
(ou ce qui revient au même, sur les réels, d'un groupe compact). 
  Dans le cas restant (le cas 3 avec $n$ impair dans la liste des cas donnés étudiés dans l'article, voir  Table 1), on peut 
    soit garder la composante neutre de Knop et Schalke (un groupe symplectique), mais il faut alors prendre pour  ${}^LG_\scrX$  
    un  produit semi-direct avec $W_\mathbb{R}$
    qui n'est pas un L-groupe (mais pas loin, c'est un \og $E$-group \fg \  au sens de \cite {ABV}), soit changer de  
    $L$-groupe et prendre à la place celui d'un groupe orthogonal pair. Concernant le morphisme  $\varphi_\scrX$, 
    ce qui est déjà prévu par Knop et Schalke est sa restriction $\eta_\scrX$ à $\SL(2,\bbC)$.
    Nous donnons dans l'article  la forme explicite de $\varphi_\scrX$ dans chaque cas.

Nous avons mentionné ci-dessus  la possibilité d'étendre l'énoncé des conjectures de Sakellaridis et Venkatesh à un contexte
 \og tordu \fg \, par un caractère unitaire $\chi$ du groupe $H$ trivial sur $H_e$, où l'on remplace  $L^2_d(\scrX)$
 par $L^2_d(\scrX)_\chi$, l'espace des sections de carré intégrable du fibré  en droite défini par $\chi$ sur $\scrX$. 
 Il s'agit alors d'introduire un $L$-groupe associé ${}^LG_{\scrX,\chi}$, le reste des conjectures s'énonçant de la même façon,
 résumée dans le diagramme ci-dessus.
 Nous ne démontrons rien dans le cas où $\chi$ est non trivial, nous espérons y revenir ultérieurement, mais nous faisons ça et là quelques remarques sur 
 ce que nous pensons être vrai. Cela éclaire parfois les résultats obtenus dans le cas où $\chi$ est trivial, en les rendant plus naturels
 lorsqu'on considère ce contexte généralisé.

Donnons maintenant un aperçu du contenu de l'article. Dans la section \ref{SDES},
nous rappelons les résultats de Oshima-Matsuki, Flensted-Jensen, Schlichtkrull et Vogan concernant les séries discrètes d'un espace symétrique réel 
 $\scrX=G/H$.  En particulier le  théorème \ref{thmOM}
donne un critère sur le rang  pour que $\scrX$ admette un spectre discret non trivial. A l'aide des tables \cite{Wiki} et \cite{KS}, Table 3, nous établissons 
alors la liste des espaces symétriques $\scrX=G/H$ avec $G$ classique
ayant un spectre discret. Cette liste, qui comporte treize cas, est la suivante:

\medskip 

\noindent  $1.$   $ \GL(n,\bbR)/\GL(p,\bbR)\times \GL(n-p,\bbR)$,  $2p\leq n$, 

\noindent  $2.$   $\U(p,q)/\Or(p,q)$,   

\noindent  $3.$ $\U(p,q)/\U(r,s)\times \U(r',s')$, $r+r'=p$, $s+s'=q$, $ r\leq r'$, $s\leq s'$,

\noindent  $4.$   $\U(n,n)/\GL(n,\bbC)$,  

\noindent  $5.$   $\U(2p,2q)/\Sp(p,q)$, 

\noindent  $6.$   $\U(n,n)/\Sp(2n,\bbR)$, 

\noindent  $7.$  $\SO(p,q)/\SO(r,s)\times \SO(r',s')$, $p+q=2n+1$, $r+r'=p$, $s+s'=q$,  $ r\leq r'$, $ s\leq s'$,

\noindent  $8.$ $\SO(p,q)/\SO(r,s)\times \SO(r',s')$, $p+q=2n$,   $r+r'=p$, $s+s'=q$,   $r\leq r'$, $s\leq s'$, 

\noindent  $9.$  $\SO(2p,2q)/\U(p,q)$, $(p-1)(q-1)\neq 0\mod 2$.

\noindent  $10.$  $\SO(n,n)/\GL(n,\bbR)$, 

\noindent  $11.$   $\Sp(2n,\bbR)/\Sp(2p,\bbR)\times \Sp(2(n-p),\bbR)$,  $2p\leq n$, 

\noindent  $12.$  $\Sp(4n,\bbR)/\Sp(2n,\bbC)$, 

\noindent  $13.$  $\Sp(2n,\bbR)/ \GL(n,\bbR)$.

 \medskip 
 

Notre point de départ est  une caractérisation précise des séries discrètes  $\pi$ de $\scrX$. 
Dans la littérature disponible, les résultats sont énoncés avec comme hypothèses $G$ semi-simple et connexe,  $H$ connexe, et sont rappelés
dans le théorème \ref{OMSV1}.
D'après \cite{Vog88},    les séries discrètes de $\scrX$ se réalisent 
 comme modules de Vogan-Zuckerman $\pi=A_\frqqq(\pi_L)$, 
  cohomologiquement induits à partir d'un d'un triplet $(\frqqq,L,\pi_L)$,  où $L$ est le centralisateur dans $G$
d'un sous-espace de Cartan compact  $\frt_0$ de l'espace symétrique $\scrX$ (c'est un $c$-Levi de $G$, dans la terminologie de Shelstad),
  $\frqqq$ est une sous-algèbre parabolique  $\theta$-stable de $\frg=\mathrm{Lie}(\mathbf G)$, et $\pi_L$  est un caractère unitaire de $L$  
  dans le fair range relativement à $\frqqq$.
  D'autre part, $K$ étant un sous-groupe compact maximal de 
 $G$ donné par une involution de Cartan $\theta$ commutant avec  $\sigma$, 
  le $K$-type minimal de $\pi$ contient  génériquement des invariants non triviaux sous le groupe $K\cap H$. En général
  on obtient les $\pi$ restants par les foncteurs de translation de Zuckerman, c'est-à-dire par continuation cohérente.
On traduit la condition sur le $K$-type minimal en la condition uniforme $(iiib)$ de la proposition \ref{iii} portant sur $\pi_L$. 
  On étend ensuite ces énoncés aux espaces symétriques de notre liste, pour lesquels $G$ n'est pas forcément semi-simple et connexe
  et $H$ n'est pas forcément connexe. Pour cela, on emploie certaines techniques usuelles et générales, mais 
  pour que l'article reste d'une longueur raisonnable, on utilise aussi des arguments au cas par cas.
   Par exemple, les groupes $G$ de notre liste peuvent être non connexes, mais d'une non connexité raisonnable, ce qui rend les arguments
  plus faciles.

 Cette réalisation des séries discrètes de $\scrX$ en termes d'induction cohomologique  permet,  en utilisant  nos 
   travaux antérieurs ({\sl cf.}  \cite{MR3}, \cite{MR7}) sur les paquets d'Arthur des groupes classiques réels, 
 de fournir pour  chaque  série discrète $\pi$ de $\scrX$ un paramètre d'Arthur explicite $\psi_\pi$ tel que 
  $\pi$ soit dans le paquet d'Arthur correspondant $\Pi(G,\psi_\pi)$.  
  On rappelle dans la section \ref{secparartbp} la forme explicite des paramètres d'Arthur  pour les groupes classiques réels (\cite{MR3}),  puis 
  dans la section \ref{parexpl}, on donne pour chaque espace de la liste, le  paramètre  $\psi_\pi$ contenant la  série discrète $\pi$.
  Pour pouvoir faire ceci, on a besoin de déterminer le groupe $L$, centralisateur dans $G$ du sous-espace de Cartan $\frt_0$.
  Ce calcul un peu laborieux est relégué à la fin de l'article (section \ref{secL}),   les résultats étant reportés dans la Table \ref{table2}.
  Notons que les calculs explicites  en rang un  sont cruciaux, car on s'y ramène par un argument général.

  Par définition, le groupe  $\mathbf L$,  centralisateur dans $\mathbf G$ du sous-espace de Cartan $\frt_0$, 
est conjugué à un sous-groupe de Levi appelé $\mathbf L(\scrX)$ dans \cite{SV}. Le groupe dual 
de ce groupe est un sous-groupe de Levi $\mathbf L_\scrX^\vee$ de $G^\vee$, qui est déterminé dans \cite{KS}, voir Table \ref{table1}.
Un argument général nous dit que la restriction de $\psi_\pi$ à $\SL(2,\bbC)$
 est le  morphisme de Jabcobson-Morosov associé à l'orbite unipotente principale  du  sous-groupe de Levi $\mathbf L_\scrX^\vee$  de $G^\vee$.
 Ainsi on peut  associer  à l'espace symétrique $\scrX$ un  morphisme  $\eta_\scrX : \, \SL(2,\mathbb{C})\rightarrow {}^LG$ 
 (ou plutôt une classe de conjugaison de tel morphismes, ou encore de manière équivalente une orbite unipotente)
 de sorte que nos  paramètres  $\psi_\pi$  aient tous ce morphisme pour restriction à $\SL(2,\mathbb{C})$  (à conjugaison près).  
 Cette construction coïncide avec  celle prévue par Sakellaridis et Venkatesh et explicitée dans 
 \cite{KS},  c'est-à-dire que la restriction  à $\SL(2,\mathbb{C})$ du morphisme conjectural  $\varphi_\scrX$ est bien $\eta_\scrX$.
Ainsi, pour ce qui est de la restriction à $\SL(2,\mathbb{C})$  des morphismes $\psi_\pi$ et $\varphi_\scrX$, la conjecture est établie.
Les détails sont donnés dans la section \ref{SV}.

Une conséquence directe de ceci  est  que la restriction des paramètres $\psi_\pi$  à $W_\mathbb{R}$ se factorise par le commutant 
dans ${}^LG$ de l'image de 
 $\SL(2,\mathbb{C})$, que nous notons $\caG$,  en un morphisme tempéré.  
  Reprenons le diagramme (\ref{factpar}) et simplifions-le en tenant compte de   la partie démontrée sur le facteur $\SL(2,\bbC)$. Si l'on note 
 $\bar \psi_\pi$ et $\bar \varphi_\scrX$ les restrictions respectives de  $\psi_\pi$ et $ \varphi_\scrX$ à  
 $W_\bbR$ et ${}^LG_\scrX$, on obtient  le diagramme commutatif (qu'il s'agit d'établir) :
  \begin{equation}\label{factpar2} \xymatrix{  W_\bbR \ar[rr]^{ \phi_\pi} \ar[rrd]_{\bar \psi_\pi}&    
 & {}^LG_\scrX \ar[d]^{\bar \varphi_\scrX} &&\\
 && \caG \ar@{^{(}->} [r] &  {}^LG}  \end{equation}

 L'étape suivante est alors de calculer explicitement $\caG$. Ceci est l'objet de la section 
 \ref{commutantSL2}.
  Dans les cas 2, 5, 6, 11, 12, 13   et  9, 10  avec $n$ pair,  qui sont les   plus simples,  ce commutant $\caG$ est 
  isomorphe au $L$-groupe d'un autre groupe réductif ayant (nécessairement) de la série discrète. 
Il  nous reste alors seulement à vérifier que c'est bien
   le $L$-groupe de \cite{KS}. 
 Dans les  cas 7, 8  et  9, 10  avec $n$ impair,  le  commutant $\caG$ est un produit direct d'un $L$-groupe d'un groupe
     réductif quasi-déployé ayant de la série discrète
     par un groupe $\caC$   fini  ou  isomorphe à  $\SO(2,\mathbb{C})$. On  vérifie alors  que le $L$-groupe est bien celui de \cite{KS}, et il reste 
 à  montrer que la projection des paramètres d'Arthur sur $ \caC $ est triviale.
 C'est là un point délicat,  qui  peut être faux pour  $H$ trop petit, nous y avons fait allusion au début de cette introduction.  
  ll faut pour cela  utiliser  la propriété clé que (génériquement) le  $K$-type minimal de $\pi$
 a  des invariants sous $K\cap H$ ou la reformulation $(iiib)$ de cette propriété.
  Dans une  version plus générale des conjectures, où l'on introduit un caractère  $\chi$ de $H$ trivial 
 sur $H_e$, on aurait une propriété du même type avec 
 un  morphisme  de $W_\mathbb{R}$ dans $\caC$  éventuellement non trivial, mais ne dépendant que de $\scrX$ et de $\chi$. 
Dans les trois cas restants, 1, 3 et 4, les racines sphériques de $\scrX$ ne sont pas toutes des racines de $\mathbf G$. 
Ceci a pour effet que le commutant $\caG$  est beaucoup plus gros que le groupe 
${}^LG_\scrX$ cherché.

La section \ref{Lplong} donne des constructions de $L$-morphismes utilisés pour construire certains des morphismes $\bar \varphi_\scrX$. 
La fin de la vérification des conjectures fait l'objet de la section \ref{DEM}.

Dans la section \ref{resComp} nous considérons les paramètres d'Arthur $\psi$ déterminés dans la section \ref{parexpl}, 
mais ici nous ne partons  pas d'une  série discrète $\pi$, ni même de la forme réelle  $G$.  Au contraire
on cherche à savoir quelles sont les représentations $\pi$ du paquet $\Pi(\psi)$ qui sont des séries discrètes pour $\scrX$,
pour un certain  $\scrX$. 
L'idée est évidemment que  tout paquet associé à un  morphisme d'Arthur  qui se factorise par ${}^LG_{\scrX}\times \SL(2,\mathbb{C})$ devrait 
contenir au moins une série discrète pour au moins une forme réelle  $\scrX$ pertinente. 
On vérifie essentiellement cette propriété, on renvoie à la section  \ref{resComp}  pour une discussion précise au cas par cas.
Rappelons que
la théorie d'Arthur attache à toute représentation $\pi$ de $G$ dans un paquet $\Pi(G,\psi)$ de paramètre $\psi$ un certain invariant 
$\epsilon(\pi)$, qui est dans les cas qui nous occupent ici un caractère du groupe $A(\psi)$ des composantes connexes
du centralisateur de $\psi$ dans $G^\vee$ (les $A(\psi)$ sont ici des 2-groupes finis). 
Les résultats de \cite{MR3}, \cite{MR7} déterminent (un choix de donnée de Whittaker
pour une forme intérieure pure quasi-déployée de $G$ étant fixée) ces caractères  $\epsilon(\pi)$.
 Il s'avère que  l'on peut  lire sur le caractère $\epsilon(\pi)$ attaché à $\pi$ si cette représentation est  une 
 série discrète de   $\scrX$. En particulier,  si deux représentations dans $\Pi(G,\psi)$ donnent le même caractère
 de $A(\psi)$ et que l'une est une série discrète pour $\scrX$, alors l'autre l'est aussi.
Les résultats sont plutôt jolis et essentiellement indépendant des choix, même si nous n'avons pas trouvé une façon uniforme de le formuler.
Ces calculs nous permettent  de constater la chose suivante : certains paramètres $\psi$
de $G$ sont obtenus comme ci-dessus pour deux espaces symétriques $\scrX=G/H$ et $\scrX'=G/H'$, mais une représentation 
$\pi$ du paquet ne contribue  au plus qu'à un seul des deux espaces symétriques.

Nous remercions très chaleureusement Y. Sakellaridis qui nous a corrigé une erreur dans une première version 
et nous a expliqué l'introduction du produit semi-direct de $\Sp(2k,\mathbb{C})\rtimes W_\mathbb{R}$ dans le cas 3.

\begin{sidewaystable}

\vspace{15cm}

\centering

\begin{tabular}{ | c | c | c  | c | c | c |}  
\hline
&              &&                           &    &  \\
&  $G/H$  & & $\check \frg_X$  &  $\check  \frl _X $ &  ${}^LG_\scrX$   \\
\hline 
 &              &&                           &    &  \\
 1& $  \GL(n,\bbR)/\GL(p,\bbR)\times \GL(n-p,\bbR) $ &   $2p \leq n  $  &   $\frs\frp(2p)$  & $ (\frg\frl_1)^{2p}\times \frg\frl_{n-2p}$ & $\Sp(2p,\bbC)\times W_\bbR $ \\ 
\hline
&              &&                           &   &   \\
2& $  \U(p,q)/\Or(p,q) $ &     & $\frg\frl(p+q)$ & $ (\frg\frl_1)^{p+q} $ & $\GL(p+q,\bbC)\rtimes W_\bbR$   \\
\hline
&              &&                           &     & \\
3& $   \U(p,q)/\U(r,s)\times \U(r',s') $ &   $r\leq r', \, s\leq s'$  &  $\frs\frp(2(r+s))$ & $ (\frg\frl_1)^{2(r+s)}\times \frg\frl_{r'+s'-(r+s)}$ &
$\begin{cases} n  \text { pair }  \Sp(2(r+s),\bbC)\times W_\bbR \\  n  \text { impair  ?? }     \end{cases}$\\
\hline
&              &&                           &     & \\
4& $   \U(n,n)/\GL(n,\bbC)  $ &     & $\frs\frp(2n)$  & $ (\frg\frl_1)^{2n}$  &   $\Sp(2n,\bbC)\times W_\bbR $ \\
\hline
&              &&                           &    &  \\
5 &$  \U(2p,2q)/\Sp(p,q)  $ &     & $\frg\frl(p+q)$ & $ (\frg\frl_2)^{p+q}$&$\GL(p+q,\bbC)\rtimes W_\bbR$ \\
\hline
&              &&                           &    &  \\
6& $   \U(n,n)/\Sp(2n,\bbR)  $ &     & $\frg\frl(n)$ & $( \frg\frl_2)^n$ &$\GL(n,\bbC)\rtimes W_\bbR$\\
\hline
&              &&                           &    &  \\
 7& $\SO(p,q)/\SO(r,s)\times \SO(r',s') $ & $p+q=2n+1  $   & $\frs\frp(2(r+s))$  & $\frs\frp(2(n-(r+s))\times (\frg\frl_1)^{r+s}$& $\Sp(2(r+s),\bbC)\times W_\bbR$\\
 & & $r\leq r',\, s\leq s'  $   &  &&    \\
 \hline
 &              &&                           &     & \\
 8& $\SO(p,q)/\SO(r,s)\times \SO(r',s') $ & $p+q=2n  $   & $\frs\fro(2(r+s)+1)$   & $\frs\fro(2(n-(r+s))\times (\frg\frl_1)^{r+s}$  &$\SO(2(r+s)+1,\bbC)\times W_\bbR$ \\
 & & $r\leq r',\, s\leq s'  $   &  & &\\
 \hline
 &              &&                           &      &\\
 9& $\SO(2p,2q)/\U(p,q)$  &$ p+q=n$& $\frs\frp(2(p+q))$ & $ (\frg\frl_2)^{\lfloor n/2 \rfloor }$ ($+\frs\fro(2) $ si $n$ est impair) &
 $\begin{cases} n  \text { pair }  \Sp(n,\bbC)\times W_\bbR \\  n  \text { impair }  \Sp(n-1,\bbC)\times W_\bbR   \end{cases} $\\
  & &${\scriptstyle (p-1)(q-1)\neq 0\; [2] } $ &  &&   \\
 \hline
 &              &   &                           &    &  \\
 10& $ \SO(n,n)/\GL(n,\bbR)$  &     & $\frs\frp(2n)$  & $ (\frg\frl_2)^{\lfloor n/2 \rfloor }$ ($+\frs\fro(2) $ si $n$ est impair)
  & $\begin{cases} n  \text { pair }  \Sp(n,\bbC)\times W_\bbR \\  n  \text { impair }  \Sp(n-1,\bbC)\times W_\bbR   \end{cases} $ \\
\hline
&              &&                           &      &\\
 11& $ \Sp(2n,\bbR)/\Sp(2p,\bbR)\times \Sp(2(n-p),\bbR) $ & $2p\leq n  $   & $\frs\frp(2p)$  & $ (\frg\frl_2)^p \times \frs\fro(2(n-2p)+1)$ &
 $\Sp(2p,\bbC)\times W_\bbR$ \\
\hline
&              &&                           &    &  \\
 12& $ \Sp(4n,\bbR)/\Sp(2n,\bbC) $ &    & $\frs\frp(2n)$    & $ (\frg\frl_2)^n $  & $\Sp(2n,\bbC)\times W_\bbR$\\
\hline
&              &&                           &     & \\
 13& $ \Sp(2n,\bbR)/\GL(n,\bbR) $ &    & $\frs\fro(2n+1)$ & $ (\frg\frl_1)^{n}$  &   $\SO(2n+1)\times W_\bbR $\\ 
 \hline
\end{tabular}
\bigskip
\caption{Espaces symétriques classiques avec séries discrètes}\label{table1}
\end{sidewaystable}

\begin{sidewaystable}

\vspace{15cm}
\centering

\begin{tabular}{ | c | c  | c | c | }  
\hline
&                         &                           &      \\
&            $G/H$   & L                        & $L\cap H$    \\
\hline 
 &                         &                           &      \\
 1& $  \GL(n,\bbR)/\GL(p,\bbR)\times \GL(n-p,\bbR) $    & $ (\bbC^\times)^p\times \GL(n-2p,\bbR) $    &  $(\bbR^\times)^p\times \GL(n-2p,\bbR)$ \\ 
\hline
&              &                           &      \\
2& $  \U(p,q)/\Or(p,q) $      & $\U(1)^n$ & $ \{ \pm 1\}^{n} $     \\
\hline
&              &                           &      \\
3& $   \U(p,q)/\U(r,s)\times \U(r',s') $  &  $(\U(1)\times \U(1))^{r+s}\times \U(r'-r,s'-s)$ & $ \U(1)^{r+s}\times \U(r'-r,s'-s)$ \\
\hline
&              &                           &      \\
4& $   \U(n,n)/\GL(n,\bbC)  $      & $(\U(1)\times \U(1))^n$ & $ \U(1)^{n}$   \\
\hline
&              &                           &      \\
5 &$  \U(2p,2q)/\Sp(p,q)  $      & $\U(2,0)^p\times \U(0,2)^q$ & $\SU(2,0)^p\times \SU(0,2)^q$ \\
\hline
&              &                           &      \\
6& $   \U(n,n)/\Sp(2n,\bbR)  $      & $\U(1,1)^n$ &  $\SU(1,1)^n$  \\
\hline
&              &                           &     \\
 7& $\SO(p,q)/\SO(r,s)\times \SO(r',s') $    & $\U(1,0)^r\times \U(0,1)^s\times \SO(p-2r,q-2s)$  & $\SO(p-2r,q-2s)$\\
 & &   &   \\
 \hline
 &              &                           &      \\
 8& $\SO(p,q)/\SO(r,s)\times \SO(r',s') $  & $\U(1,0)^r\times \U(0,1)^s\times \SO(p-2r,q-2s)$  & $\SO(p-2r,q-2s)$ \\
 &    &  &     \\
 \hline
 &              &                           &      \\
 
 9a& $\SO(2p,2q)/\U(p,q), \, p,q \text{ pair}$     & $ \U(2,0)^{p/2}\times \U(0,2)^{q/2}$   &   $\SU(2,0)^{p/2}\times \SU(0,2)^{q/2}  $ 
 \\
 \hline
 &              &                           &      \\
 
 9b& $\SO(2p,2q)/\U(p,q), \, p \text{ pair},\,  q \text{ impair} $    & $ \U(2,0)^{p/2}\times \U(0,2)^{q-1/2}  \times \SO(0,2)$ & 
 $\SU(2,0)^{p/2}\times \SU(0,2)^{q-1/2}  \times \SO(0,2)  $ 
 \\
 \hline
 &              &                           &      \\

 10a& $ \SO(n,n)/\GL(n,\bbR), \, n \text{ pair }$       & $\U(1,1)^{n/2} $  & $\SU(1,1)^{n/2} $   \\
\hline
&              &                           &       \\
 10b& $ \SO(n,n)/\GL(n,\bbR), \, n \text{ impair }$       & $\U(1,1)^{n-1/2} \times \SO(1,1)$  & $\SU(1,1)^{n-1/2} \times \SO(1,1) $   \\
\hline
&              &                           &       \\

 11& $ \Sp(2n,\bbR)/\Sp(2p,\bbR)\times \Sp(2(n-p),\bbR) $   & $\U(1,1)^p\times \Sp(2(n-2p),\bbR)$  & $\SU(1,1)^p\times \Sp(2(n-2p),\bbR)$  
  \\
\hline
&              &                           &      \\
 12& $ \Sp(4n,\bbR)/\Sp(2n,\bbC) $ &     $\U(2)^n$    & $ \SU(2)^n$   \\
\hline
&              &                           &      \\
 13& $ \Sp(2n,\bbR)/\GL(n,\bbR) $     & $\U(1)^n$ & $ \{\pm 1\}^{n}$     \\ 
 \hline
\end{tabular}
\bigskip
\caption{$L$ et $L\cap H$}\label{table2}
\end{sidewaystable}

\section{Séries discrètes des espaces symétriques}\label{SDES}

\subsection{Notations pour les espaces symétriques}  \label{EsSy}

Si $M$ est un groupe de Lie,  on note $M_e$ la composante connexe de l'identité dans  $M$. 
Si   $\mathbf M$ un groupe algébrique réductif, on note $\mathbf M^0$ sa composante neutre au sens des groupes algébriques.

Soit $\mathbf G$ un groupe algébrique réductif connexe défini sur $\bbR$ et soit $G=\mathbf G(\bbR)$ le groupe de ses points réels.
Soit $\sigma$ une involution de $\mathbf G$, définie sur $\bbR$. Posons $\mathbf H=\mathbf G^\sigma$ et 
soit $H$ un sous-groupe de $G$ avec 
$ \mathbf H(\bbR)_e\subset H\subset \mathbf H(\bbR)$.

Notons $\scrX$ l'espace symétrique $G/H$.
Les sous-représentations irréductibles de $G$ dans  $L^2(\scrX)$ forment un sous-espace noté $L^2_d(\scrX)$ et  sont appelées séries discrètes de $\scrX$.
Elles ont été décrites grâce aux travaux de nombreux mathématiciens (\cite{FlJ}, \cite{OM},  \cite{Sch}, \cite{Vog88});  en particulier, Vogan
donne une description en termes d'induction cohomologique que nous allons  rappeler ci-dessous, après avoir introduit les notations nécessaires.

\begin{rmq} \label{bigproblem} Dans ces références, 
on suppose souvent $G$ connexe et/ou  semi-simple et parfois aussi $H$ connexe. 
 Il nous faudra  donc  étendre ces résultats  en les adaptant  à nos groupes classiques, dont certains ne sont pas   semi-simples connexes.  
\end{rmq}

Rappelons  d'abord un  critère dû à Oshima et Matsuki pour que  $\scrX$ admette des séries discrètes.
Fixons  une involution de Cartan $\theta$ de $G$ qui commute à $\sigma$, et soit   $K=G^\theta$  le sous-groupe compact maximal de
$G$ associé à $\theta$. On note $\frg_0$, $\frh_0$, $\frk_0$ les algèbres de Lie de $G$, $H$, $K$ respectivement, et plus généralement, l'algèbre de Lie d'un groupe 
de Lie   est notée par la lettre gothique correspondante avec un indice $0$.  Lorsqu'on enlève l'indice $0$, ceci désigne alors l'algèbre de Lie complexifiée.
De plus, les involutions $\theta$ et $\sigma$ donnent lieu à des décompositions :
\[ \frg_0=\frk_0\oplus \frp_0 , \quad \frg_0=\frh_0\oplus \frs_0, \quad 
\frg_0=(\frk_0\cap \frp_0)\oplus  (\frk_0\cap \frs_0)\oplus(\frp_0\cap \frh_0)\oplus(\frp_0\cap \frs_0).\]

Rappelons qu'un sous-espace de Cartan de $\frg_0$ pour $\scrX$ est une 
sous-algèbre abélienne,  composée d'éléments semisimples, contenue   dans   $\frs_0$ et  maximale pour l'inclusion
avec ces propriétés.  Le rang de l'espace symétrique $\scrX$ est égal à la dimension d'un sous-espace de Cartan.
 Notons $\scrX_K=K/K\cap H$ :  c'est un espace symétrique compact.

\begin{thm} [Oshima-Matsuki] \label{thmOM}
L'espace symétrique  $\scrX$ admet des séries discrètes si et seulement si le rang de 
 $\scrX$ est égal au rang de $\scrX_K$,  c'est-à-dire qu'un sous-espace de Cartan   de $\frk_0$ pour $\scrX_K$ est  aussi un sous-espace de Cartan  de $\frg_0$
pour $\scrX$. 
\end{thm}

Un sous-espace de Cartan de $\frg_0$ pour $\scrX$ comme dans le théorème sera appelé sous-espace de Cartan compact.

\begin{rmq}  Comme expliqué dans la remarque \ref{bigproblem}, le théorème ci-dessus est démontré dans la littérature pour un espace symétrique $G/H$ avec  
$G$ connexe et semi-simple. Le cas 
 plus général envisagé ici se réduit au  cas  connexe et semi-simple de manière routinière.
\end{rmq}

\subsection{Listes des espaces symétriques classiques avec séries discrètes}\label{liste}

Nous nous plaçons maintenant dans le  cas où le groupe $\mathbf G$ est un groupe classique, c'est-à-dire que $G$ est 
un groupe général linéaire $\GL(n,\bbR)$, ou bien un groupe unitaire $\U(p,q)$, ou bien un groupe symplectique
$\Sp(2n,\bbR)$, ou bien un groupe spécial orthogonal $\SO(p,q)$.
Pour déterminer   la liste des espaces symétriques  $\scrX=G/H$ avec $G$ comme ci-dessus, et vérifiant de plus la condition 
du théorème \ref{thmOM} d'existence de séries discrètes,  on utilise 
 la classification de Wikipedia \cite{Wiki} pour avoir tous les  $\scrX$ possibles, et la condition de rang se vérifie en utilisant
  la table 3 de \cite{KS} : le rang de $\scrX$ est égal au rang de son algèbre de Lie duale 
$\check g_\scrX$. On peut alors  comparer au rang de l'espace symétrique compact 
$\scrX_K=K/K\cap H$. On obtient  la liste  de 13  cas donnée dans l'introduction.

\subsection{Notations pour l'induction cohomologique}\label{ICo}
Dans ce paragraphe, on garde les mêmes notations que précédemment concernant le groupe $G$, mais on oublie l'espace symétrique.

Soit $\frt_0$ une sous-algèbre  abélienne de $\frk_0$.
Posons $L=\mathrm{Cent}_G(\frt_0)$. C'est un facteur de Levi $\theta$-stable de $G$ (un $c$-Levi dans la terminologie de Shelstad). 
Soit $\frqqq=\frl\oplus \fru$ une sous-algèbre parabolique $\theta$-stable de $\frg$.
On note $\caR_{\frqqq,L,G}^i$ le foncteur d'induction cohomologique de Vogan-Zuckerman ({\sl cf.} \cite{Vgreen}, \S 6.3.1)
en degré $i$, de la catégorie des $(\frl,L\cap K)$-modules vers la catégorie des $(\frg,K)$-modules.
Nous renvoyons à \cite{KnVo}, Def. 0.49 et 0.52 pour les définitions du good range et du fair range. 
Dans cet article, nous allons toujours considérer des inductions cohomologiques à partir de caractères unitaires dans le fair range.
Dans ce contexte, le degré qui nous intéresse particulièrement, et même exclusivement,  est $S=\dim(\fru\cap \frk)$, et  nous 
posons pour tout  caractère unitaire $\pi_L$ de $L$, 
\[A_\frqqq(\pi_L)=\caR^S_{\frqqq, L,G}(\pi_L).\]

\subsection{Séries discrètes de $\scrX$ comme $A_\frqqq(\pi_L) $ }\label{SDX}
Reprenons les notations du paragraphe \ref{EsSy}. Supposons que $\scrX$ vérifie la condition du théorème \ref{thmOM} et admet 
donc des séries discrètes. Fixons un sous-espace de Cartan compact de $\frt_0$  pour $\scrX$
(tous les choix possibles pour  $\frt_0$ sont conjugués sous $K\cap H$).
Comme dans le paragraphe \ref{ICo}, posons $L=\mathrm{Cent}_G(\frt_0)$.
On a dans ce cas $\frl=(\frl\cap \frh)\oplus \frt$, c'est la décomposition en sous-espaces propres pour $\sigma$.
On dit que $L$ et $\frl$ sont associés à l'espace symétrique $\scrX$. 
Notons $\scrQ(\frl)$ l'ensemble des sous-algèbres paraboliques 
$\theta$-stables  $\frqqq$ de $\frg$ de facteur de Levi $\frl$.

Nous allons maintenant donner une description des séries discrètes  pour les espaces symétriques de notre liste 
comme modules $A_\frqqq(\pi_L)$.  Dans la littérature  (\cite{FlJ}, \cite{OM},\cite{Sch}, \cite{Vog88})
cette description est disponible sous les hypothèses $G$ et $H$ connexes et $G$ semi-simple. Nous allons donc partir de 
là et adapter les résultats aux espaces de notre liste.  On fait donc (provisoirement) ces hypothèses dans ce qui suit.

Soit $\frqqq=\frl\oplus \fru \in \scrQ(\frl)$.   On introduit  d'abord un ensemble de paramètres relatif à $\frqqq$ qui va paramétrer la plus grande partie 
des séries discrètes.
 
\begin{defi} \label{defSD'}
Soit  $\caP'(\frqqq)$ l'ensemble des $(\frl,L\cap K)$-modules irréductibles $\pi_L$ vérifiant les conditions $(i)'$,  $(ii)'$ et $(iii)'$ suivantes.

\medskip

$(i)'$  $\pi_L$ est un caractère unitaire de $L$, de différentielle nulle sur $\frh_0\cap \frl_0$.

$(ii)'$ $\pi_L$ est dans le  good range pour $\frqqq$.

 Le bottom-layer ({\sl cf. \cite{KnVo}}, § V. 6) de $A_\frqqq(\pi_L)$ est alors 
est constitué du  $K$-type $V_\mu$ de plus haut poids $\mu=\pi_L\otimes \bigwedge^{\mathrm{top}}(\fru\cap \frp)$. 
Ce $K$-type $V_\mu$ est le $K$-type 
minimal de $A_\frqqq(\pi_L)$.

 $(iii)'$ Le $K$-type  $V_\mu$ de plus haut poids $\mu=\pi_L\otimes \bigwedge^{\mathrm{top}}(\fru\cap \frp)$ 
 a des invariants non triviaux sous $H\cap K$.
\end{defi}

 Ici, ce que nous entendons par \og  plus haut poids \fg \, d'une représentation irréductible $V$ de $K$ 
 est la représentation de $L\cap K$ dans $H^0(\fru\cap \frk,V)$  qui caractérise $V$ (Kostant).

Pour obtenir toutes les séries discrètes, il faut autoriser $\pi_L$ à décrire tout le fair range. Le problème dans ce cas est que le bottom-layer
de $A_\frqqq(\pi_L)$ peut être vide, et la condition $(iii)'$ n'a plus de sens. Pour compléter la définition,  Vogan 
remarque que les séries discrètes de $\scrX$ forment des familles cohérentes (voir \cite{Vog88} pour une définition précise de ces familles). 
On introduit donc un ensemble de paramètres plus grand.

\begin{defi} \label{defSD}
Soit  $\caP(\frqqq)$ l'ensemble des $(\frl,L\cap K)$-modules irréductibles $\pi_L$ vérifiant les conditions $(i)$,  $(ii)$ et $(iii)$ suivantes.

\medskip

$(i)$  $\pi_L$ est un caractère unitaire de $L$, de différentielle nulle sur $\frh_0\cap \frl_0$. 

$(ii)$ $\pi_L$ est dans le  fair range pour $\frqqq$.

  $(iii)$ Le module $A_\frqqq(\pi_L)$ est dans une famille cohérente de modules $A_\frqqq(\pi'_L)$ dont génériquement les paramètres $\pi'_L$
  vérifient les conditions $(ii)'$ et $(iii)'$ de la définition précédente.
\end{defi}

Dans le fair range, il se peut que  $A_\frqqq(\pi_L)$ soit nul.  S'il ne l'est pas, Vogan a démontré (\cite{Vog88}) qu'il est irréductible.

\begin{prop}\label{iii}
La condition  $(iii)'$ de la définition \ref{defSD'} est équivalente  (sous les conditions (i)' et (ii)') à la suivante : 

$(iiib)$  le caractère $\pi_L\otimes \bigwedge^{\mathrm{top}}(\fru\cap \frp)$ de $L\cap K$  est 
trivial sur $L\cap H\cap K$.

\noindent La condition  $(iii)$ de la définition \ref{defSD} est équivalente à (iiib)  (sous les conditions (i) et (ii)).

\end{prop}

\dem  La première assertion est démontrée en utilisant la dualité de Flensted-Jensen et le théorème de Cartan-Helgason, 
ceci se trouve dans  \cite{FlJ}, \cite{Sch},\cite{OM}.  Cette condition étant stable pour les familles  cohérentes, on obtient la seconde assertion.\qed

\begin{thm}[Oshima-Matsuki-Schlichtkrull-Vogan]  \label{OMSV1}Soit $\scrX=G/H$ un espace symétrique avec  $G$ et $H$ connexes et $G$  semi-simple.
On a une décomposition 
\[L^2(\scrX)= \bigoplus_{   \frqqq\in \scrQ(\frl), \pi_L\in \caP(\frqqq)} A_\frqqq(\pi_L). \]
D'autre part les $A(\pi_L)$, sont irréductibles ou nuls.
\end{thm}
La première assertion est due à Oshima-Matsuki, avec une  description différente des $A(\pi_L)$. Le lien avec la description en termes de modules 
cohomologiquement induits est due à Schlichtkrull  et Vogan, et l'irréductibilité est due à  à Vogan.
Pour faire le lien avec la littérature, la remarque suivante peut s'avérer utile.

\begin{rmq} \label{Gconnexe}  L'action de $\frl\cap \frk\cap \frh$ sur $\bigwedge^{\mathrm{top}}(\fru\cap \frp)$ est triviale.
On étend le sous-espace de Cartan $\frt_0$ en une sous-algèbre de Cartan $\frt_0\oplus \frt'_0$ de $\frk_0$
avec $\frt'_0$ dans $\frh_0\cap \frk_0$. Le sous-espace $\frt'_0$ est donc une sous-algèbre de Cartan de 
$\frd_0:=\frh_0\cap \frk_0\cap \frl_0$. Comme  $\bigwedge^{\mathrm{top}}(\fru\cap \frp)$ est une représentation de dimension $1$, l'action 
de  $[\frd,\frd]$ sur $\bigwedge^{\mathrm{top}}(\fru\cap \frp)$ est triviale. Il reste donc à voir que l'action de $\frt'$ est triviale.
L'argument est donné dans \cite{Sch}, en bas de la page 140.
\end{rmq}

Voyons comment le théorème \ref{OMSV1} peut s'étendre à des espaces symétriques plus généraux. 

\medskip

\noindent {\bf  Première étape}. Elle consiste à obtenir une description analogue pour $\scrX=G/H$ avec les mêmes hypothèses sur $G$, mais
sans hypothèse sur $H$ autre que $(G^\sigma)_e=H_e \subset H\subset G^\sigma$.

Remarquons que l'on a une action par translation à droite de $H$ sur $L^2(G/H_e)$ commutant avec l'action par translation à gauche de $G$. 
En effet, si $f\in L^2(G/H_e)$,  on a quels que soient $g\in G$, $h'\in H$ et $h\in H_e$, 
\[  (r(h').f)(gh)=   f(ghh')=f(gh'({h'}^{-1}hh')=f(gh')= (r(h').f)(g).\] 
 Cette action est triviale sur $H_e$. On a donc une action de $H/H_e$ sur $L^2(G/H_e)$ commutant avec l'action par translation à gauche de $G$. 
 Fixons une série discrète  $\pi=A_\frqqq(\pi_L)$  (suffisamment régulière,  disons $\pi_L\in \caP'(\frqqq)$
 de   $L^2(G/H_e)$ et considérons l'espace 
 \[ \scrA= \Hom_{\frg,K} (\pi, L^2(G/H_e)).\]
 
 On sait que $\pi$  a multiplicité un dans  $L^2(G/H_e)$, donc $\scrA$ est de dimension $1$. 
 L'action de $H/H_e$ sur $\scrA$ héritée de l'action sur  $L^2(G/H_e)$ est donc un caractère de ce groupe, et l'on en conclut 
 qu'il existe un caractère $\chi$ de $H/H_e$ tel que pour tout 
 $ \phi\in \scrA$, 
 \[[h\cdot \phi : \, v\mapsto r(h)\cdot \phi(v)] = [(\chi(h) \phi)  : \, v\mapsto  \chi(h) \phi(v) ] \] 
 et donc $\phi(v)$ est dans l'espace noté précédemment $L^2(G/H)_\chi$.
 Par réciprocité de Frobenius, on a 
  \[ \scrA= \Hom_{\frg,K} (\pi, L^2(G/H_e)) \simeq \scrB=\Hom_{\frh,K\cap H_e} (\pi,\chi).  \]
et l'action de $H/H_e$ sur $\scrA$ se transporte en une action sur $\scrB$ (par ce même caractère $\chi$).
 On a un morphisme de restriction non nul  au $K$-type minimal $V_\mu$ de $\pi$  : 
 \[  \scrB=\Hom_{\frh,K\cap H_e} (\pi,\chi)   \longrightarrow   \Hom_{K\cap H_e} (V_\mu,\chi). \]
 On en conclut que $V_\mu$ possède des vecteurs se transformant par  $\chi$ sous l'action de $H$.
 On a donc obtenu  la première partie du résultat suivant : 
 
 \begin{prop} \label{GHnc}
 Soit $\scrX=G/H$ un espace symétrique avec $G$   semi-simple connexe, admettant des séries discrètes tordues par un caractère $\chi$ de $H/H_e$.
 Alors
  \[L^2(\scrX)_\chi= \bigoplus_{   \frqqq\in \scrQ(\frl), \pi_L\in \caP(\frqqq)_\chi   } A_\frqqq(\pi_L). \]
où l'espace des paramètres  $\caP(\frqqq)_\chi $ est défini comme dans les définitions \ref{defSD'} et \ref{defSD}
en remplaçant les conditions $(iii)'$ et $(iii)$ respectivement  par 

$(iii)' 1$.   Le $K$-type  $V_\mu$ de plus haut poids $\mu=\pi_L\otimes \bigwedge^{\mathrm{top}}(\fru\cap \frp)$ 
 a des vecteurs non nuls se transformant selon le caractère $\chi$ pour l'action de  $H\cap K$.

$(iiib) 1$.  Le caractère $\pi_L\otimes \bigwedge^{\mathrm{top}}(\fru\cap \frp)$ de $L\cap K$  coïncide avec le caractère $\chi$ sur 
 $L\cap H\cap K$.

 \end{prop}
 
\dem  Il nous faut maintenant reformuler la condition $(iiib)$ et il faut revenir à la démonstration de l'équivalence entre $(iii)$' et $(iiib)$
dans le cas connexe. On utilise la dualité de  Flensted-Jensen entre $(G, H_e,K)$ et un triplet similaire $(G^d, K^d,H^d)$
de groupes ayant pour complexification respective $\mathbf{G}$, $\mathbf K$ et $\mathbf H^0$ ($\mathbf H^0$ est la composante neutre du 
groupe réductif $\mathbf H =\mathbf G^\sigma)$.  De plus  $H^d$ est un sous-groupe compact maximal de $G^d$. 
On voit alors  $V_\mu$  comme une représentation  algébrique de   $\mathbf K$ ayant des invariants sous $\mathbf H^0\cap \mathbf K$
et par restriction, comme une représentation de dimension finie de 
$K^d$ ayant des invariants sous $H^d\cap K^d$.  Or  $K^d \cap H^d$ est un sous-groupe compact maximal de 
$K^d$. On remarque que $\mathbf L\cap \mathbf K$ est un sous-groupe de Levi de $\mathbf K$. Le théorème de 
Cartan-Helgason nous dit alors que le plus haut poids $\mu$ est invariant sous $\mathbf L\cap \mathbf K \cap \mathbf H^0$.

Ceci constitue  la démonstration de l'équivalence entre $(iii)'$ et $(iiib)$ et il s'agit de l'étendre  avec $\mathbf H$ à la place de $\mathbf H^0$.
 Pour cela  on fait agir $\mathbf H\cap  \mathbf K$ sur $\Hom_{\mathbf H^0\cap \mathbf K}(V_\mu,1)$ ce qui par 
restriction donne une action de  $\mathbf (\mathbf H\cap \mathbf K)(\bbC)/(\mathbf H^0\cap \mathbf K)(\bbC) $ sur l'espace 
$\Hom_{( \mathbf L\cap \mathbf K^0 \cap \mathbf H)(\mathbb{C})}(\mu,1)$. La conclusion est alors évidente. 
  \qed

 \medskip

\noindent{\bf Deuxième étape}. L'étape suivante de la généralisation consiste à autoriser $G$ à être réductif, mais toujours connexe. On écrit alors 
 $G=Z(G)_eG_1$ où $G_1$ est  semi-simple connexe et $Z(G)$ est le centre de $G$. Les groupes $G_1$ et $Z(G)_e$ sont stables par l'involution 
 $\sigma$. On fait tout d'abord  l'hypothèse suivante, vérifiée dans certains  cas de notre liste où $G$ est non réductif. 
 
 \noindent {\bf Hypothèse 1} : $Z(G)_e\subset H_e$.

  Posons $H_1=G_1\cap H$. On a alors un difféomorphisme, induit par l'inclusion de $G_1$ dans $G$ :
   \[  G_1/ H_1   \longrightarrow G/H . \]
 L'injectivité est immédiate puisque l'on a quotienté par $G_1\cap H$ et la surjectivité est conséquence de l'hypothèse.
 Pour tout caractère $\chi$ de $H_1/H_{1,e}$, on a donc une description de 
 $L_d^2(\scrX)_\chi$ comme représentation de $G_1$ donné par la proposition \ref{GHnc}.
 D'autre part, grâce à l'hypothèse 1, $Z(G)_e$ agit trivialement sur     $L_d^2(\scrX)_\chi$. Enfin
 \[ H_1/ H_{1,e}= G_1\cap H/  (G_1\cap H)_e  \twoheadrightarrow    G_1\cap H/  G_1\cap H_e \simeq H/H_e.  \]
 On peut donc identifier les caractères de $H/H_e$ avec des caractères de $H_1/ H_{1,e}$.
 Les modules entrant dans la décomposition de  $L_d^2(\scrX)_\chi$ comme représentation de $G$
 sont donc les $A_\frqqq(\pi_L) $ de la décomposition comme représentation de $G_1$,
 étendu   par l'action  triviale du centre $Z(G)_e$.

 \begin{exemple} Pour étudier le  premier exemple de notre liste, 
 on considère donc déjà l'espace symétrique  
 \[G/H=\GL^+(n,\bbR)/(\GL(n,p)\times \GL(n-p,\bbR))^+,\]
 où le $+$  indique que l'on prend  le sous-groupe des matrices de déterminant strictement positif.
 Le groupe $Z(G)_0$ est alors le sous-groupe des matrices scalaires $\lambda I_n$ avec $\lambda>0$, 
 $G_1$ est le groupe $\SL(n,\bbR)$, $H_1=\mathbf{S} (\GL(n,p)\times \GL(n-p,\bbR)) $  et l'hypothèse 1 est clairement vérifiée.
 On a de plus ici $H_1/ H_{1,e} \simeq H/H_e\simeq \bbZ/2\bbZ$.  Les cas 3 et 4 de la liste sont aussi des exemples où l'hypothèse 
 1 est vérifiée, et de plus dans ces cas, les groupes sont tous connexes.
 
 Remarquons que la  condition sur l'action du  centre est déjà dans la condition $(i)$ de  la définition de l'espace des paramètres, sous l'hypotèse 
 ci-dessus   $Z(G)_e\subset H_e$
 \end{exemple}
 
 Traitons maintenant les cas restants de notre liste  où $G$ n'est pas semi-simple. Commençons par la variante du cas   2 avec $H$ connexe,
  $G/H=\U(p,q)/\SO(p,q)$. On a alors $G_1=\SU(p,q)$,   $Z(G)=Z(G)_e=\U(1)$. Ce qui va servir ici est la compacité de ce centre.
   Posons $\scrX_1=G_1/\SO(p,q)$. Considérons l'application naturelle, surjective
 \[  \U(1)\times \scrX_1= \U(1)\times \SU(p,q)/\SO(p,q) \longrightarrow  \U(p,q)/\SO(p,q)=\scrX. \]
 C'est un fibré principal de fibre $\boldsymbol{\mu}_n$, le groupe des racines $n$-ième de l'unité.
Le théorème \ref{OMSV1} donne une description du spectre discret de $\scrX_1$.
On en déduit facilement une description du spectre discret de $\U(1)\times \scrX_1$ comme représentation de $\U(1)\times \SU(p,q)$, puis du spectre discret  
 de $\scrX$ comme représentation de $\U(1)\times_{\boldsymbol{\mu}_n}\SU(p,q)=\U(p,q)$. 
 On passe de $\U(p,q)/\SO(p,q)$ à  $\U(p,q)/\Or(p,q)$ comme dans la première étape ci-dessus.
 
 La conclusion pour ce cas est donc que la description de $L_d^2(\scrX)$ du théorème \ref{OMSV1} est valide sans aucun changement de la définition de 
 l'espace des paramètres. La réduction ci-dessus s'applique aussi de la même façon aux cas 5 et 6.
 
 \begin{rmq} Le cas général pour un groupe $G$ réductif connexe peut sans doute être obtenu par une combinaison des cas considérés ci-dessus, c'est-à-dire 
 $Z(G)_e\subset H$ ou $Z(G)_e$ compact. En effet la condition sur le rang de $\scrX$ du théorème \ref{thmOM} implique que la 
partie de  $Z(G)_e$ qui n'est pas dans $H$ est compacte.
 \end{rmq}
 
 \medskip
 
 \noindent{\bf Troisième étape}.  Il s'agit maintenant d'étendre les résultats aux groupes $G$ non connexes.
 Remarquons que pour les cas de notre liste, la non-connexité est limitée puisque $\vert G/G_e\vert \leq \bbZ/2\bbZ$ (cas 1, 7, 8, 9,10).
   Il est clair que $L^2_d(G/H_e)$ est l'induite de   $L^2_d(G_e/H_e)$ de $G_e$ à $G$.
 On a donc une description de $L^2_d(G/H_e)$ comme somme directe de modules $\Ind_{G_e}^G(A_\frqqq(\pi_{L_e}))$. 
D'autre part, on a aussi $\vert L/L_e \vert \leq 2$, et si $\pi_{L_e}$ est un caractère unitaire de $L_e$, alors 
$\Ind_{L_e}^L(\pi_{L_e})$ est somme de caractères $\pi_L$ de $L$. 
Par \cite{KnVo}, on a  $\Ind_{G_e}^G(A_\frqqq(\pi_{L_e}))= \bigoplus  A_\frqqq(\pi_{L})$ où la somme est sur les caractères 
$\pi_L$ apparaissant dans $\Ind_{L_e}^L(\pi_{L_e})$. La description du $K$-type minimal
des   $A_\frqqq(\pi_{L})$ de plus haut poids $\pi_L\otimes \bigwedge^{\mathrm{top}}(\fru\cap \frp)$ est similaire.
On en déduit que la description du spectre discret $L^2_d(G/H_e)$ similaire à celle de  $L^2_d(G/H_e)$, 
avec la même définition de l'espace des paramètres.

\medskip

 \noindent{\bf Quatrième  étape}. On passe maintenant de $G/H_e$ avec $G$ dans notre liste et non connexe
 à $G/H$ avec $H$ simplement soumis à $G^\sigma_0\subset H\subset G^\sigma$. La réduction se fait comme dans la première étape.

\section{Paramètres d'Arthur des séries discrètes de $\scrX$}\label{secparartSD}

Dans cette section, on donne la forme générale des paramètres d'Arthur (de bonne parité) des groupes classiques, puis 
pour chaque série discrète $\pi$  d'un espace symétrique  $\scrX$ de la liste, nous donnons un paramètre d'Arthur explicite 
$\psi_\pi$ tel que le paquet d'Arthur attaché à ce paramètre contienne $\pi$.  
On note  $R[a]$ la représentation irréductible algébrique de $\SL(2,\bbC)$ de dimension $a$.

\subsection{Paramètres d'Arthur de bonne parité}\label{secparartbp} Ce qui suit est tiré de \cite{MR3} auquel nous renvoyons le lecteur pour plus de détails.

\subsubsection{$\GL(n,\bbR)$  }\label{GLnR}

Dans ce paragraphe, $G=\GL(n,\bbR)$, dont le  $L$-groupe est 
$\GL(n,\bbC)\times W_\bbR$. Un paramètre d'Arthur
$\psi: \, W_\bbR \times \SL(2,\bbC)\rightarrow \GL(n,\bbC)$ peut se voir comme une représentation de
dimension $n$ de $W_\bbR \times \SL(2,\bbC)$ et cette représentation est semi-simple. On peut donc la décomposer
en somme directe de représentations irréductibles $\psi_i$, chaque $\psi_i$ étant un produit tensoriel extérieur
d'une représentation irréductible de $W_\bbR$, disons $\phi_i$, et d'une représentation $R[a_i]$ de  $\SL(2,\bbC)$.
Les représentations irréductibles de $W_\bbR$ sont de dimension $1$ où $2$, et rappelons 
que la restriction d'un paramètre d'Arthur à $W_\bbR$ est un paramètre de Langlands tempéré, ce doit donc être le cas 
de ces irréductibles.  Celles de dimension 2 sont notées
$\delta(s_1,s_2)$. Elles sont paramétrées par  $s_1$, $s_2\in \bbC$, avec $s_1-s_2\in \bbZ_{>0}$ et 
$s_1+s_2\in i\bbR$. Par la correspondance de Langlands, on associe à $\delta(s_1,s_2)$
une série discrète de $\GL(2,\bbR)$ que l'on note aussi  $\delta(s_1,s_2)$ et dont le caractère infinitésimal, 
identifié de la manière usuelle à un ensemble de deux nombres complexes, est $\{s_1,s_2\}$.
Les représentations irréductibles de dimension 1 de $W_\bbR$ sont notées $\eta(\epsilon,s)$, avec $\epsilon\in \{0,1\}$ et $s\in i\bbR$. 
Elles sont associées par la correspondance de Langlands aux caractères unitaires $x\mapsto \sgn(x)^\epsilon\vert x \vert^s $ de $\GL(1,\bbR)$.
Ainsi, on décompose un paramètre d'Arthur $\psi$ en une somme directe
\[ \psi=\bigoplus_{i=1}^m \left(  \delta(s_{i,1},s_{i,2} )\boxtimes R[a_i]  \right)\oplus \bigoplus_{j=1}^u  \left( \eta(\epsilon_j,s_j)\boxtimes R[a'_j]\right).  \]
avec \begin{equation} \label{partn}2\sum_{i=1}^m  a_i + \sum_{j=1}^u  a'_j=n.\end{equation} 

Le paquet d'Arthur associé à $\psi$ est un singleton, consistant en une représentation unitaire irréductible notée $\pi^{\GL}(\psi)$ de 
$\GL(n,\bbR)$. Cette représentation peut se décrire de la manière suivante.
La donnée de la série discrète $\delta(s_{i,1},s_{i,2} )$ et de l'entier $a_i$
détermine une représentation $\Speh(\delta(s_{i,1},s_{i,2} ),a_i)$ de $\GL(2a_i,\bbR)$
et le caractère $\eta(\epsilon_j,s_j)$ de $\bbR^\times $ donne par composition
avec le déterminant un caractère, encore noté $ \eta(\epsilon_j,s_j)$, de $\GL(a'_j,\bbR)$.
Soit $P=P_{2a_1,\ldots,2a_m,a'_1,\ldots,a'_u}$ le sous-groupe parabolique standard de $\GL(n,\bbR)$ associé à la partition
(\ref{partn}) de $n$. Son facteur de Levi standard est $M=\left( \prod_{i=1}^m \GL(2a_i,\bbR)\right) \times \left(\prod_{j=1}^u\GL(a'_j,\bbR)\right)$.
La représentation $\pi^{\GL}(\psi)$ est alors l'induite parabolique  (irréductible) de $P$ à $G$ de la représentation 
\[ \left( \boxtimes_{i=1}^m \Speh(\delta(s_{i,1},s_{i,2} ),a_i)   \right) \boxtimes \left(\boxtimes_{j=1}^u \eta(\epsilon_j,s_j)  \right).\]
Le paramètre de Langlands de  $\pi^{\GL}(\psi)$  est $\phi_\psi$, où $\phi_\psi$ est le paramètre de Langlands
attaché au paramètre d'Arthur $\psi$ ({\sl cf.} \cite{Art82}, p. 10).

\subsubsection{Groupes unitaires}\label{Gu}
Dans ce paragraphe, $G$ est un groupe unitaire de rang $n$. Son $L$-groupe est donc
$\GL(n,\bbC)\rtimes W_\bbR$. Considérons un paramètre d'Arthur pour $G$ :
\[\psi: \, W_\bbR\times \SL(2,\bbC)\longrightarrow \GL(n,\bbC)\rtimes W_\bbR \] 
que l'on suppose de bonne parité ({\sl cf.} \cite {MR7}, Déf. 2.1). Cela signifie que la  restriction de $\psi$ à
$W_\bbC\times \SL(2,\bbC)$ (encore notée $\psi$) est de la forme 
\begin{equation}\label{parunit} \psi= \bigoplus_{i=1}^R (\chi_{m_i}\boxtimes R[a_i])\end{equation}
où $\sum_{i=1}^R a_i=n$, $m_i\in \bbZ$, $a_i\in \bbZ_{>0}$,  $\chi_{m_i}(z)=\left(  \frac{z}{\bar z}  \right)^{\frac{m_i}{2}}$, avec la condition de parité 
$m_i+a_i  \cong n \mod 2$. 
On suppose, ce qui est loisible, que  les $m_i$ sont rangés dans l'ordre décroissant : $m_1\geq m_2\geq \cdots\geq m_R$.
Notons qu'il n'existe, à conjugaison par $G^\vee=\GL(n,\bbC)$ près, qu'une seule manière d'étendre (\ref{parunit}) à $W_\bbR$.

\subsubsection{Groupes symplectiques et orthogonaux}
  Si   $G=\Sp(2n,\bbR)$,  ${}^LG=\SO(2n+1,\bbC)\times W_\bbR$. Si $G=\SO(p,q)$,  $p+q=2n+1$,  ${}^LG=\Sp(2n,\bbC)\times W_\bbR$.
  Notons   \[\Std_G: \,{}^LG \longrightarrow \GL(N,\bbC) \] 
  la représentation standard de ${}^LG$,  avec respectivement $N=2n+1$ et $N=2n$.
 Si  $G=\SO(p,q)$ avec $p+q=2n$, on a alors ${}^LG=\SO(2n,\bbC)\rtimes W_\bbR$, l'action de $W_\bbR$ sur $\SO(2n,\bbC)$ étant triviale 
si  $n-p \cong 0 \mod 2$ et si $n-p \cong 1 \mod 2$, cette  action se factorise par $\mathrm{Gal}(\bbC/\bbR)$, 
l'élément non trivial agissant par un automorphisme extérieur, réalisé par l'action adjointe d'un élément de $\Or(2n,\bbC)\setminus \SO(2n,\bbC)$, 
  en préservant un épinglage. On a donc dans les deux cas un morphisme 
  \begin{equation}\label{LSOO} \zeta_{SO}: {}^LG=\SO(2n,\bbC)\rtimes W_\bbR\longrightarrow \Or(2n,\bbC), \end{equation}
 et l'on définit la représentation standard de $ {}^LG$ en composant ce morphisme $\zeta_{SO}$ avec une réalisation  de 
 $\Or(2n,\bbC)$ dans $\GL(2n,\bbC)$ : $\Std_G: \,{}^LG \longrightarrow \GL(2n,\bbC)$. On pose dans ce cas $N=2n$.

Soit $G$ l'un de ces groupes classiques.
 Composons un  paramètre d'Arthur pour $G$,  \[ \psi : \, W_\bbR\times \SL(2,\bbC)\longrightarrow  W_\bbR={}^LG  \]
 avec la représentation standard de ${}^LG$,  $\Std_G: \,{}^LG \longrightarrow \GL(N,\bbC)$.
 Notons  encore $\psi$ cette composition. On dit que $\psi$ est de bonne parité si  $\psi$ se décompose sous la forme : 
      \begin{equation}\label{parclass} 
       \psi= \bigoplus_{i=1}^R\left(  \delta\left(\frac{m_i}{2}, -\frac{m_i}{2} \right) \boxtimes R[a_i]\right)
\oplus \bigoplus_{j=1}^S \left( \sgn_{W_\bbR}^{b_j} \boxtimes R[c_j]\right),
\end{equation} 
 où $m_i\in \bbZ_{>0}$, $a_i\in \bbZ_{>0}$,  $2\sum_{i=1}^Ra_i+\sum_{j=1}^S c_j=N$, 
avec selon le type de groupe, les conditions de parité suivantes : 
\[\bullet\;  G=\Sp(2n, \bbR) :   m_i+a_i  \cong 1 \mod 2 \, (i=1,\ldots, R),\,  c_j \cong 1 \mod 2,\,  (j=1,\ldots, S),\qquad \qquad\qquad{} \]
 $\sum_{j=1}^sb_j  \cong f \mod 2$, où  $f $ est le nombre de  $ a_i $  impairs.
\[\bullet\;  G=\SO(p,q), \, p+q=2n+1 :  m_i+a_i  \cong 0 \mod 2 \, (i=1,\ldots, R),\,  c_j \cong 0 \mod 2,\,  (j=1,\ldots, S).\]
\[\bullet\;   G=\SO(p,q), \, p+q=2n :  m_i+a_i   \cong 1 \mod 2 \, (i=1,\ldots, R),\,  c_j \cong 1 \mod 2,\,  (j=1,\ldots, S).\]
 $\sum_{j=1}^sb_j   \cong  f+n-p  \mod 2$, où $f$ est le nombre de $a_i$ impairs.

\subsection{Séries discrètes de $\scrX$ et leur paramètre}\label{parexpl}

Reprenons les notations des sections \ref{EsSy} et \ref{SDX} relatives à un espace symétrique
$\scrX$ admettant des séries discrètes.
 On a donc un sous-espace de Cartan  compact $\frt_0$ et un $c$-Levi $L=N_G(\frt_0)$ dont le complexifié est 
$\mathbf L=N_{\mathbf G}(\frt_0)$.
Les séries discrètes de $\scrX$ sont des modules $A_\frqqq(\pi_L)$ où $\frqqq$ est une sous-algèbre parabolique $\theta$-stable de 
$\frg$ et $\pi_L$ une représentation de $L$ vérifiant les hypothèses de la définition \ref{defSD}. Lorsque $G$ est un groupe 
classique (c'est-à-dire dans cet article, un groupe général linéaire, un groupe unitaire ou un groupe spécial orthogonal), 
les résultats de \cite{MR3} et \cite{MR7}  
montrent que l'on a le résultat suivant, déjà bien connu dans le cas des groupes linéaires généraux.

\begin{prop} \label{resSL2} Soit $\caL$ le sous-groupe de Levi dual de $\mathbf L$ dans $G^\vee$. Soit $\pi$ une  série discrète de $\scrX$.
Alors $\pi$ est contenue   dans un  paquet d'Arthur de paramètre
\[\psi_\pi: \, W_\bbR \times \SL(2,\bbC)\longrightarrow {}^LG \]
tel que la restriction de $\psi_\pi$ à $\SL(2,\bbC)$ soit à conjugaison près  le morphisme de Jacobson-Morozov associé à l'orbite unipotente régulière
de $\caL$.
\end{prop}

Dans chaque cas de la liste,  les sous-groupe de Levi $\mathbf L$ et $\caL$ sont donnés dans la section  \ref{SV}.
La forme réelle explicite $L$ de $\mathbf L$ est plus délicate à déterminer. Les calculs, un peu fastidieux, sont données dans la section \ref{secL}
et les résultats sont résumés dans la table \ref{table2}. 
La restriction de $\psi_\pi$ à $\SL(2,\bbC)$ est donc déterminée (à conjugaison près). On procède alors cas par cas.

 \subsection{Le cas 1: $\scrX=\GL(n,\mathbb{R})/\GL(p,\mathbb{R})\times \GL(n-p,\mathbb{R})$,  $2p\leq n$ } Nous allons montrer
le résultat suivant.

\begin{prop} Soit $\pi$ une série discrète de $\scrX$.  Alors il existe des  séries discrètes $(\delta_i)_{i=1,\ldots,p}$  
de $\GL(2,\mathbb{R})$ de caractère central trivial et toutes distinctes,     telles que la représentation $\pi$ soit obtenue par l'induction parabolique
 à partir de la représentation $\left(\bigotimes_{i=1}^p \delta_i\right)\otimes \Triv_{\GL(n-2p,\mathbb{R})}$
 d'un sous-groupe de Levi standard $M= \GL(2,\bbR)^p\times\GL(n-2p,\mathbb{R})$ de $\GL(n,\bbR)$.
\end{prop}

Avant de démontrer cette proposition, on en tire le corollaire suivant, conséquence du fait  
que  pour $\GL(n,\mathbb{R})$, le paramètre d'Arthur est facilement obtenu à partir du  paramètre de Langlands
(voir section \ref{GLnR}). 

\begin{cor}
Le paramètre d'Arthur d'une série discrète  $\pi$ de $\scrX$  réalisée par induction parabolique comme dans la proposition   est
\[\psi_\pi =
\bigoplus_{i=1}^p (\delta_i \boxtimes R[1])\oplus (\Triv_{W_\bbR} \boxtimes R[n-2p]),
\]
où les $\delta_i$, $i=1,\ldots, p$  sont ici  les représentations de dimension deux de $W_\mathbb{R}$ qui paramètrent les séries discrètes de $\GL(2,\bbR)$
notées de la même façon. De plus,  ces séries discrètes sont de caractère 
central trivial et les $\delta_i$ sont donc à valeurs dans $\SL(2,\mathbb{C})\times W_\mathbb{R}$
c'est-à-dire (avec les notations de \ref{GLnR}), $\delta_i=\delta\left( \frac{m_i}{2}, -\frac{m_i}{2} \right)$, avec  les $m_i$ entiers impairs.
\end{cor}

\noindent \underline{\sl Démonstration de la proposition}. On sait  par la proposition \ref{resSL2} ci-dessus, et le fait que 
$\caL\simeq \GL(1,\bbC)^{2p}\times \GL(n-2p,\bbC)$ 
que le paramètre  d'Arthur $\psi_\pi$ ci-dessus  a pour restriction
 à $\SL(2,\mathbb{C})$ le morphisme de Jacobson-Morozov correspondant à l'élément unipotent régulier de $\caL$. 
 Ainsi $\psi_\pi$ est nécessairement de la forme:
\[ \psi_\pi=  (\phi_d  \boxtimes R[1]) \oplus (\eta(\epsilon,s)  \boxtimes R[n-2p]),\] 
où $\epsilon\in \{0,1 \}$  et où $\phi_d $ est un morphisme tempéré de $W_\mathbb{R}$ à valeurs dans
 $\GL(2p,\mathbb{C})\times W_\mathbb{R}$. On commence par montrer que ce morphisme tempéré est en fait une somme de séries discrètes. 

Ici, exceptionnellement on utilise l'induction cohomologique comme Vogan l'a normalisée  dans \cite{Vog86}. L'avantage majeur
 est que cette induction cohomologique commute à l'induction parabolique (\cite{Vog86}, Theorem 17.6). 
 On sait a priori que $\pi=\widetilde A_\frqqq(\tilde \pi_L)$, les $\, \tilde{}\, $ étant mis ici pour distinguer cette induction cohomologique
 normalisée de celle utilisée précédemment.
   Ici $L\simeq (\bbC^\times)^p\times \GL(n-2p,\bbR)$ (voir section \ref{secL}). 
Via cet isomorphisme, le caractère unitaire $\tilde \pi_L$ de $L$ est donné par des caractères unitaires $\chi_i$, $i=1,\ldots,  p$ de $\bbC^\times$
et d'un caractère unitaire   $\eta$ de $\GL(n-2p,\mathbb{R})$. La condition $(i)$ de la définition \ref{defSD} implique que $\eta$ est 
le caractère trivial ou le signe du déterminant 
et que chaque caractère $\chi_i$ est de  la forme $\chi_i(z)=(z/\overline{z})^{m_i/2}$ où les  $m_i$ sont entiers. 

La condition de fair range sur l'induite   cohomologique $\pi=\widetilde A_\frqqq(\tilde \pi_L)$
 s'exprime (pour un bon choix de $\mathfrak{q}$,  mais pour  $G=\GL(n,\bbR)$ tous les choix sont conjugués sous $K$, et ce choix n'a donc aucune conséquence)
 par  $m_1>\cdots >m_p >0$. 
En notant $\delta_i=\delta\left(\frac{m_i}{2},-\frac{m_i}{2}\right)$ (notations de \ref{GLnR}), la commutation de l'induction
 cohomologique et de l'induction parabolique dit que
 $\pi$ est obtenu par induction parabolique
 à partir de la représentation $\left(\bigotimes_{i=1}^p \delta_i\right)\otimes (\sgn(\det))^\epsilon$
 du  sous-groupe de Levi standard $M= \GL(2,\bbR)^p\times \GL(n-2p,\mathbb{R})$ de $G$. Ici $\epsilon\in \{0,1\}$.
On veut  maintenant montrer pour terminer la démonstration  que les $m_i$ sont tous impairs et que $\epsilon=0$.
Le caractère infinitésimal est dans un système de coordonnées bien choisi
\[ \left( \frac{m_1}{2}, \cdots, \frac{m_p}{2}, \frac{n-2p-1}{2}, \cdots, -\frac{n-2p-1}{2}, -\frac{m_p}{2} \cdots,-\frac{m_1}{2} \right).
\]
Le $K$-type minimal de $\pi$  est donné par Vogan en \cite{Vog86},  6.5 (a) et (6.12) 
 et l'on voit  que ce $K$-type minimal   est la représentation de $K=\Or(n,\mathbb{R})$ de plus haut poids
$$( 2m_1+1, \cdots, 2m_p+1, \underbrace{0,\cdots, 0}_{[n/2]-p})_z ,
$$
où $z$ est relié à $\epsilon$ de la façon suivante: si $n$ est impair $z$ se calcule avec le caractère central de $\pi$ et si $n$ est pair
 $z=0$ si et seulement si $\epsilon=0$. On renvoie à \cite{Vog86} pour les notations ci-dessus concernant les représentations
 des groupes orthogonaux.

On exploite maintenant le fait que $\pi$ est une série discrète de $\scrX$  et donc que son $K$-type minimal a des invariants sous $K\cap H$ 
(condition $(iii)$ de la définition \ref{defSD})
c'est-à-dire sous $\Or(p)\times \Or(n-p)$. On utilise d'abord la condition $(i)$ de \cite{Sch} page 138, qui donne que les $m_i$ sont tous de même parité.
 Le vecteur de plus haut poids doit être invariant sous $\Or(1)^p\times \Or(n-2p)$ ce qui force d'abord  les $m_i+1$ à être pairs. 
 Ensuite si $n$ est pair et $n>2p$, par  définition de $z$, l'invariance sous $\Or(n-2p)$ force  $z=0$. Si $n$ est impair, le caractère central est précisément $z$ 
 et doit être trivial,  d'où encore $z=0$. Ainsi on a bien montré que les $m_i$ sont des entiers impairs  et 
 les $\delta_i$ sont  des séries discrètes de $\GL(2,\mathbb{R})$ de caractère central trivial. 
 On a aussi montré que $\epsilon=0$ si $n$ est pair et si $n$ est impair cela est forcé par le calcul du caractère central.
\qed

\medskip 

 \begin{rmq} Il est aussi possible de donner la fin de la démonstration  sans recourir à l'induction 
 cohomologique subtilement normalisée de \cite{Vog86}. On écrit comme dans les autres cas $\pi=A_\frqqq(\pi_L)$
 avec l'induction cohomologique usuelle, et on traduit la condition $(iiib)$ de  la proposition  \ref{iii}  en une condition 
 sur $\pi_L$ en calculant explicitement le caractère $\bigwedge^{\mathrm{top}}(\fru\cap \frp)$ de $L\cap K\cap H$. Ceci fait apparaître
 un facteur $\sgn(\det)^p$ sur le facteur $\GL(n-2p,\bbR)$ qui disparait lorqu'on prend en compte la non commutativité de l'induction 
 parabolique et de l'induction cohomologique pour calculer le paramètre de Langlands de $\pi$.
 \end{rmq}
 
 \medskip
 
   Ecrivons de manière explicite le paramètre $\psi_\pi$ obtenu :
 \begin{equation}\label{psigl} 
 \psi_\pi=\bigoplus \left(  \delta\left(\frac{m_i}{2},-\frac{m_i}{2}\right)\boxtimes R[1]  \right)\oplus \left(  \Triv_{W_\bbR}\boxtimes R[n-2p]  \right).
 \end{equation}
 avec les $m_i$ entiers impairs distincts.

 \begin{rmq} 
 Il est vraisemblable que les paramètres similaires mais 
avec des $m_i$ pairs ou bien une partie unipotente  $ \sgn_{W_\bbR}\boxtimes R[n-2p] $ 
pourraient apparaître dans le cas des conjectures généralisées avec un caractère $\chi$ de $H/H_e$ non trivial. 
 \end{rmq}

\subsection{Le cas 2: $\scrX=\U(p,q)/\Or(p,q) $,  $n=p+q$}   Dans ce cas, le groupe $L$ est un tore compact,
 et les séries discrètes de $\scrX$ sont des séries discrètes de $G=\U(p,q)$.
La restriction des paramètres $\psi_\pi$ à $\SL(2,\bbC)$ est donc triviale et avec les notations de (\ref{parunit}), 
  \begin{equation}\label{psicas2} \psi_\pi= \bigoplus_{i=1}^n (\chi_{m_i}\boxtimes R[1])\end{equation}
où $m_i\in \bbZ$,   $\chi_{m_i}(z)=\left(  \frac{z}{\bar z}  \right)^{\frac{m_i}{2}}$, avec la condition de parité 
$m_i \cong n-1 \mod 2$.

\subsection{Le cas 3: $\scrX= \U(p,q)/\U(r,s)\times \U(r',s')$, $ r\leq r'$, $s\leq s'$}

On a    ici  $$L=\U(1)^{2(r+s)} \times \U(r-r',s-s').$$ 
 Si  $\pi=A_\frqqq(\pi_L)$ est une série discrète de $\scrX$,  alors son paramètre $\psi_\pi$ est de la forme
\[  \psi_\pi=\bigoplus_{i=1}^{2(r+s)}  (\chi_{m_i} \boxtimes R[1])   \oplus \left(\chi_{m_0} \boxtimes R[n-2(r+s)] \right) \]
 où les $\chi_i$ sont comme dans (\ref{parunit}), pour certains entiers $m_i$ non nuls.
 
Voyons maintenant ce que donne la condition $ (i)$ de la définition \ref{defSD}.
 Comme  $L\cap K$ est  ici un produit de groupes unitaires, donc connexes, 
 Le caractère $\pi_L$ doit être  trivial sur $L \cap  H$.
 
Le groupe  $L\cap H$ contient le facteur  $\U(p-2r,q-2s)$. La restriction de $\pi_L$ à ce facteur doit donc être triviale, ce qui implique d'après \cite{MR7} que $m_0=0$. 
La forme des caractères unitaires de $L$ triviaux sur 
$L\cap  H$  (voir section  \ref{secL}, le calcul en rang un du cas 3)
 montre alors  que  les $m_i$ se regroupent par paires $m_i$, $-m_i$ et ils sont entiers, de la parité de $n-1$. 
  Pour conclure, le paramètre $\psi_\pi$ est de la forme
  \begin{equation}\label{psicas3}
   \psi_\pi=\bigoplus_{i=1}^{r+s}  ((\chi_{m_i}\oplus \chi_{-m_i}) \boxtimes R[1])   \oplus \left(\chi_0 \boxtimes R[n-2(r+s)] \right) \end{equation}
où les $m_i$ sont entiers, de la parité de $n-1$.

\subsection{Le cas 4: $\scrX= \U(n,n)/\GL(n,\bbC)$}  On a    ici  $L=(\U(1) \times \U(1))^{n}$.  C'est donc un tore compact.
 Si  $\pi=A_\frqqq(\pi_L)$ est une série discrète de $\scrX$,  alors son paramètre $\psi_\pi$ est de la forme
\[  \psi_\pi=\bigoplus_{i=1}^{2n}  (\chi_{m_i} \boxtimes R[1])    \]
 où les $\chi_i$ sont comme dans (\ref{parunit}), pour certains entiers $m_i$ non nuls.
 Le groupe $L\cap H$ est isomorphe à $\U(1)^n$ ou chaque facteur $\U(1)$ se plonge  dans le facteur $\U(1) \times \U(1)$ comme expliqué dans la section 
    \ref{secL}. 
 Comme dans le cas précédent , la condition $(i)$ nous dit alors que les $m_i$ se regroupent par paires $m_i$, $-m_i$ et ils sont entiers, 
de la parité de $n-1$.    On a donc finalement 
  \begin{equation}\label{psicas4}
 \psi_\pi=\bigoplus_{i=1}^{n}  ((\chi_{m_i}\oplus \chi_{-m_i}) \boxtimes R[1])  \end{equation}
avec les $m_i$ entiers impairs.  

\subsection{Le cas 5:    $\U(2p,2q)/\Sp(p,q)$}   On a ici $L=\U(2,0)^p\times \U(0,2)^q$.  
Si  $\pi=A_\frqqq(\pi_L)$ est une série discrète de $\scrX$,  alors son paramètre $\psi_\pi$ est de la forme
 \begin{equation}\label{psicas5}
 \psi_\pi=\bigoplus_{i=1}^{2n}  (\chi_{m_i} \boxtimes R[2])  \end{equation}
et les $m_i$ sont ici pairs car $n=2p+2q$ est pair.

 \subsection{Le cas 6:   $\U(n,n)/\Sp(2n,\bbR)$ } On a ici $L=\U(1,1)^n$.  
Si  $\pi=A_\frqqq(\pi_L)$ est une série discrète de $\scrX$,  alors son paramètre $\psi_\pi$ est de la forme
 \begin{equation}\label{psicas6}
 \psi_\pi=\bigoplus_{i=1}^{2n}  (\chi_{m_i} \boxtimes R[2])     \end{equation}
et les $m_i$ sont ici pairs car $2n$ est pair.

 \subsection{Le cas 7:  $\SO(p,q)/\SO(r,s)\times \SO(r',s')$, $p+q=2n+1$, $r+r'=p$, $s+s'=q$,  $ r\leq r'$, $ s\leq s'$}
 On a ici $L= \U(1,0)^r\times \U(0,1)^s\times \SO(p-2r,q-2s)$.
  Si  $\pi=A_\frqqq(\pi_L)$ est une série discrète de $\scrX$,  alors son paramètre $\psi_\pi$ est de la forme
 \[
  \psi_\pi=\bigoplus_{i=1}^{r+s}  \left(\delta\left( \frac{m_i}{2},  - \frac{m_i}{2}\right) \boxtimes R[1] \right)  \oplus \left(\sgn_{W_\bbR}^\epsilon \boxtimes R[2(n-r-s)]  \right) .   \]
où les $m_i$ sont impairs et $\epsilon\in \{\pm 1\}$.   Montrons que de plus $\epsilon=0$.

Notons  $L_0=\SO(p-2r,q-2s)$ le  facteur   spécial orthogonal de $L$.
Notons $\eta_{SO}$ le caractère non trivial d'un groupe orthogonal $\SO(k,\ell)$ non connexe, c'est-à-dire avec 
$k\ell\neq 0$ (si $k\ell=0$, $\SO(k,\ell)$ est connexe et $\eta_{SO}$ est alors par convention le caractère trivial).
Utilisons maintenant la condition $(iii)$  (ou  de manière équivalente $(iiib) $)   de la définition \ref{defSD}. On peut conclure de deux manières.
Il découle de \cite{MR3} que 
l'élément $\pi=A_\frqqq(\pi_L)$ est de restriction $\eta_{SO}^\epsilon$ au facteur $L_0$. 
La condition $(iiib)$ nous dit que l'on veut que $\pi_L\otimes \bigwedge^{\mathrm{top}}(\mathfrak{u}\cap 
\mathfrak{p})$ soit trivial sur $L\cap H\cap K$. Or un calcul explicite de  $\bigwedge^{\mathrm{top}}(\mathfrak{u}\cap 
\mathfrak{p})$ (qui nécessite d'introduire beaucoup de notations, en particulier des systèmes de racines explicites, et qui de fait est le calcul
fait dans l'appendice de  \cite{MR3} pour prouver l'assertion précédente)
montre que $L\cap K\cap H$ agit trivialement sur cet espace.
 Ainsi la restriction de $\pi_L$ à $L\cap H\cap K$ doit être triviale, et ceci impose $\epsilon=0$.
L'autre manière consiste à utiliser \cite{SpVo}, Thm. 4.23,  qui donne le lien entre paramètre de Langlands et le $K$-type 
minimal de plus haut poids $\pi_L\otimes \bigwedge^{\mathrm{top}}(\mathfrak{u}\cap 
\mathfrak{p})$ du module $\pi$. La forme de ce paramètre de Langlands imposé par la condition $(iii)$ 
montre alors directement que l'on doit avoir $\epsilon=0$ dans le paramètre d'Arthur $\psi_\pi$ ci-dessus. 
Finalement, on a obtenu :
\begin{equation}\label{psicas7} 
  \psi_\pi=\bigoplus_{i=1}^{r+s}  \left(\delta\left( \frac{m_i}{2},  - \frac{m_i}{2}\right) \boxtimes R[1] \right)  \oplus \left(\Triv_{W_\bbR} \boxtimes R[2(n-r-s)]  \right) .   \end{equation}

 \subsection{Le cas 8:  $\SO(p,q)/\SO(r,s)\times \SO(r',s')$, $p+q=2n$,   $r+r'=p$, $s+s'=q$, $s+s'=q$,  $ r\leq r'$, $ s\leq s'$}
 On a ici $L= \U(1,0)^r\times \U(0,1)^s\times \SO(p-2r,q-2s)$.
  Si  $\pi=A_\frqqq(\pi_L)$ est une série discrète de $\scrX$,  alors son paramètre $\psi$ est de la forme
 \begin{equation}\label{psicas8} 
  \psi_\pi=\bigoplus_{i=1}^{r+s}  \left(\delta\left( \frac{m_i}{2},  - \frac{m_i}{2}\right) \boxtimes R[1] \right)  \oplus \left( \eta_1 \boxtimes R[2(n-r-s)-1]  \right)
  \oplus \left(\eta_2 \boxtimes R[1] \right)   \end{equation}
où les $m_i$ sont pairs.
Le même argument utilisant la condition $(iii)$ ou $(iiib)$ de la définition \ref{defSD} donne ici $\eta_1,\eta_2\in \{\Triv_{W_\bbR}, \sgn_{W_\bbR}\}$ avec $\eta_1\eta_2=\sgn_{W_\bbR}$ si $p$ et $q$ sont de la parité de $n$, 
$\eta_1\eta_2=\Triv_{W_\bbR}$ si $p$ et $q$ sont de la parité de $n-1$.

 \subsection{Le cas 9:   $\SO(2p,2q)/\U(p,q)$} On a ici $L= \U(2,0)^{p/2}\times \U(0,2)^{q/2}$ si $p$ et $q$ sont pairs 
 et   $ \U(2,0)^{p/2}\times \U(0,2)^{q-1/2}  \times \SO(0,2)$ si $p$ est pair et  $q$  impair.  
 Si  $\pi=A_\frqqq(\pi_L)$ est une série discrète de $\scrX$,  alors son paramètre $\psi_\pi$ est de la forme
\begin{equation}\label{psicas9a} 
  \psi_\pi=\bigoplus_{i=1}^{\frac{n}{2}}  \left(\delta\left( \frac{m_i}{2},  - \frac{m_i}{2}\right) \boxtimes R[2] \right)  \end{equation}
si $n$ est pair et 
\begin{equation}\label{psicas9b}   \psi_\pi=\bigoplus_{i=1}^{\frac{n-1}{2}}  \left(\delta\left( \frac{m_i}{2},  - \frac{m_i}{2}\right) \boxtimes R[2] \right)  \oplus \left(\sgn_{W_\bbR}\boxtimes R[1]\right) \oplus
 \left( \Triv_{W_\bbR} \boxtimes R[1]\right)   \end{equation}
si $n$ est impair.
Dans les deux cas, les $m_i$ sont impairs.



 \subsection{Le cas 10:   $\SO(n,n)/\GL(n,\bbR)$}  On a ici  $L=\U(1,1)^{n/2} $   si $n$ est pair et 
  $L=\U(1,1)^{n-1/2}\times \SO(1,1) $ si $n$ est impair.  Si  $\pi=A_\frqqq(\pi_L)$ est une série discrète de $\scrX$,  alors son paramètre $\psi_\pi$ est de la forme

\begin{equation}\label{psicas10a} 
  \psi_\pi=\bigoplus_{i=1}^{\frac{n}{2}}  \left(\delta\left( \frac{m_i}{2},  - \frac{m_i}{2}\right) \boxtimes R[2] \right)  \end{equation}
si $n$ est pair et 
\[   \psi_\pi=\bigoplus_{i=1}^{\frac{n-1}{2}}  \left(\delta\left( \frac{m_i}{2},  - \frac{m_i}{2}\right) \boxtimes R[2] \right) 
 \oplus \left(\sgn_{W_\bbR}^\epsilon\boxtimes R[1]\right) \oplus
 \left( \sgn_{W_\bbR}^\epsilon \boxtimes R[1]\right)   \]
si $n$ est impair, avec   $\epsilon\in \{0,1\}$. Dans les deux cas les $m_i$ sont impairs.
 Le même argument utilisant la condition $(iiib)$ de la définition \ref{defSD} que dans le cas 7  donne ici 
$\epsilon=0$. D'où finalement : 
\begin{equation}\label{psicas10b}   \psi_\pi=\bigoplus_{i=1}^{\frac{n-1}{2}}  \left(\delta\left( \frac{m_i}{2},  - \frac{m_i}{2}\right) \boxtimes R[2] \right) 
 \oplus \left(\Triv_{W_\bbR} \boxtimes R[1]\right) \oplus
 \left( \Triv_{W_\bbR} \boxtimes R[1]\right)   \end{equation}

 \subsection{Le cas 11:    $\Sp(2n,\bbR)/\Sp(2p,\bbR)\times \Sp(2(n-p),\bbR)$,  $2p\leq n$} Soit $\pi=A_\frqqq(\pi_L)$  une série discrète de $\scrX$. Comme
  ici 
 $ L= \U(1,1)^p\times  \Sp(2(n-2p),\bbR)$,  le paramètre $\psi_\pi$ est de la forme :  
 \begin{equation}\label{psicas11}  \psi_\pi=\bigoplus_{i=1}^{p}  \left(\delta\left( \frac{m_i}{2},  - \frac{m_i}{2}\right) \boxtimes R[2] \right) \oplus
 \left( \Triv_{W_\bbR}\boxtimes R[2n-4p+1]\right) ,    \end{equation}
où  les $m_i$ sont impairs.

 \subsection{Le cas 12:   $\Sp(4n,\bbR)/\Sp(2n,\bbC)$} Soit $\pi=A_\frqqq(\pi_L)$  une série discrète de $\scrX$.  Comme on a ici $L= \U(2)^n$,   
le paramètre $\psi_\pi$ est de la forme : 
  \begin{equation}\label{psicas12}  \psi_\pi=\bigoplus_{i=1}^{n}  \left(\delta\left( \frac{m_i}{2},  - \frac{m_i}{2}\right) \boxtimes R[2] \right) \oplus
 \left( \Triv_{W_\bbR}\boxtimes R[1]\right)     \end{equation}
où  les $m_i$ sont  impairs.

\subsection{Le cas 13: $\scrX=\Sp(2n,\bbR)/\GL(n,\bbR) $}   Dans ce cas, le groupe $L$ est un tore compact,
 et les séries discrètes de $\scrX$ sont des séries discrètes de $G=\Sp(2n,\bbR)$.  Si $\pi=A_\frqqq(\pi_L)$ est une série discrète de $\scrX$,
    alors son paramètre $\psi_\pi$ est de la forme
  \begin{equation}\label{psicas13} \psi_\pi= \bigoplus_{i=1}^n \left(\delta\left( \frac{m_i}{2}, -\frac{m_i}{2}\right)   \boxtimes R[1]\right)   \oplus 
  ( \sgn_{W_\bbR}^n \boxtimes R[1])
  \end{equation}
où les  $m_i$ sont pairs.

\section{Conjectures de Sakellaridis et Venkatesh : première réduction} \label{SV}\label{debut}\label{liste2}

Soient $G$, $H$, $\scrX$, etc. comme dans la section  \ref{EsSy}. 
Comme expliqué dans l'introduction  Sakellaridis et Venkatesh introduisent dans \cite{SV} 
 un formalisme du $L$-groupe associé à l'espace symétrique $\scrX$.
Ces constructions sont précisées par Knop et Schalke dans \cite{KS}.  Donnons des précisions sur leurs constructions.
Ils associent d'abord    à $\scrX$ un sous-groupe parabolique noté $\mathbf P(\scrX)$ 
de $\mathbf G$ et un facteur de Levi  $\mathbf L(\scrX)$ de celui-ci. Le sous-groupe $\mathbf P(\scrX)$ est un sous-groupe parabolique
$\sigma$-déployé minimal de $\mathbf G$ et $\mathbf L(\scrX)=\mathbf P(\scrX)\cap \sigma(\mathbf P(\scrX))$.
Soit $\mathbf A$ un tore maximal de $\mathbf L(\scrX)$ et posons  
$ \mathbf A_\scrX= \mathbf L(\scrX)/\mathbf L(\scrX)\cap \mathbf H= \mathbf A /\mathbf A \cap \mathbf H$.
C'est un tore d'algèbre de Lie $\fra_\scrX$. 
Il résulte aussi des constructions que les algèbres de Lie $\frt$ et $\fra_{\scrX}$ sont isomorphes.
 Ainsi le rang de $\scrX$ est égal à la dimension de $\mathbf A_\scrX$.

\medskip

Notons $G^\vee$ le groupe dual de $\mathbf G$. On fixe une réalisation ${}^LG= G^\vee\rtimes W_\bbR$ du $L$-groupe de 
$\mathbf G$, l'action de $W_\bbR$ sur $G^\vee$ se factorise par $\mathrm{Gal}(\bbC/\bbR)$ et préserve un épinglage 
$(B^\vee, T^\vee,\{\caX_\alpha\}_\alpha)$ de $G^\vee$.
On fixe aussi un épinglage de $\mathbf G$, ce qui nous permet d'identifier racines simples de $\mathbf G$ et  coracines de $G^\vee$.
Le sous-groupe de Levi  $\mathbf L(\scrX)$ de $\mathbf G$ du paragraphe précédent détermine donc via la dualité et le choix des épinglages,
 un sous-groupe de Levi  $\mathbf L_\scrX ^\vee$ de $G^\vee$. On note 
  \begin{equation}\label{phiSL2}\eta_\scrX: \, \SL(2,\bbC)\longrightarrow \mathbf L_\scrX ^\vee \subset G^\vee\end{equation}
   un morphisme de Jacobson-Morosov associé à l'orbite unipotente principale de $\mathbf L_\scrX ^\vee$.

\medskip

Supposons que $\scrX$ vérifie le critère d'existence de séries discrètes du théorème \ref{thmOM}
\begin{prop} \label{egalL}
 Les groupes $\mathbf L$ et $\mathbf L(\scrX)$ (identifiés à leur groupes de points complexes)
sont conjugués sous $\mathbf G$. Dualement, les groupes $\caL$ et $\mathbf L_\scrX^\vee$ sont donc conjugués dans $G^\vee$.
\end{prop}
Ceci résulte de la remarque 2.1.1 de \cite{SV}.

\medskip

 Une des conditions portant sur ${}^L G_{\scrX}$ est que le rang de   $G_{\scrX}^\vee $ 
 est égal à celui de l'espace symétrique $\scrX$. 
 D'autre part, si $G_{\scrX}^\vee $ est déterminé, alors ${}^L G_{\scrX}$
l'est aussi par le fait que nous voulons l'existence de paramètres de Langlands {\sl discrets} $\phi_d: \, W_\bbR\rightarrow  {}^L G_{\scrX}$.
En effet  ${}^L G_{\scrX}$   est  alors aussi le  $L$-groupe d'un groupe quasi-déployé  admettant  des séries discrètes.  Or
  pour un groupe réductif complexe connexe, il n'y a qu'une seule classe de forme réelle quasi-déployée ayant des séries discrètes, celle qui est 
  une forme intérieure de la forme compacte.   La structure de produit semi-direct de ${}^L G_{\scrX}$ est alors fixée par cette condition.
  Or $G_{\scrX}^\vee $ a été déterminé par Knop et Schalke. Ceci fixe donc a priori ${}^L G_{\scrX}$.
Comme conséquence de  la  proposition \ref{resSL2} et de la proposition ci-dessus,  on obtient le corollaire suivant : 

\begin{cor}\label{SL2} Soit $\pi$ une  série discrète de $\scrX$. Alors  $\pi$ est contenue dans un  paquet d'Arthur de paramètre
$\psi_\pi$  tels que la restriction de $\psi$ à $\SL(2,\bbC)$ soit à conjugaison près  le morphisme de Jacobson-Morozov
$\eta_\scrX$  associé à l'orbite unipotente régulière
de $\mathbf L_\scrX^\vee$.
\end{cor}

Ce corollaire constitue une étape importante dans la vérification de la conjecture, puisqu'il dit que 
la restriction de $\psi_\pi$  à $\SL(2,\bbC)$ est bien obtenue via le morphisme (\ref{phiSL2}).

\medskip

Dans la table  \ref{table1}, on donne dans chacun des 13 cas de  la liste donnée dans l'introduction  :

---  en troisième colonne, les conditions sur les variables définissant $\scrX$, 

--- en quatrième  colonne,  l'algèbre de Lie  $\check \frg_X$ du groupe $G^\vee_\scrX$,

---  en cinquième colonne,  l'algèbre de Lie  $\check \frl_X$ du groupe $\mathbf L^\vee_\scrX$.

--- en sixième colonne le groupe ${}^LG_\scrX$ (sauf dans le cas 3, $n$ impair, qui fera l'objet d'une discussion spéciale).

Les  colonnes 4 et 5 sont issues de la table 3 de \cite{KS}, le $L$-groupe ${}^LG_\scrX$ est déterminé par les considérations ci-dessus.
\begin{rmq}\label{detL}
L'algèbre de Lie  $\check \frl_X$ détermine le groupe des points complexes $\mathbf L$ du $c$-Levi $L$ 
de $G$  introduit dans la section \ref{SDX}.  En effet, le groupe $\mathbf L^\vee_\scrX$
est dual du groupe $\mathbf L (\scrX)$, ce dernier étant conjugué dans $\mathbf G$
à $\mathbf L$ (remarque \ref{egalL}). L'inspection de la liste montre que  $\check \frl_X$ est toujours
un produit de $\frg\frl$ de rang un ou deux et d'une algèbre de Lie classique de même type que $\check \frg$.
\end{rmq}

\section{Plongements de $L$-groupes}\label{Lplong}

Dans cette section, nous construisons des $L$-morphismes qui nous seront utiles dans la suite de cet article, ou éventuellement dans un article ultérieur
portant sur la généralisation des conjectures de Sakellaridis et Venkatesh lorsque l'on considère un caractère  $\chi$ de $H/H_e$  non trivial.

\subsection{Le groupe de Weil}
Le groupe de Weil $W_\bbR$ est engendré par son sous-groupe $W_\bbC=\bbC^\times$ et un élément $j$, avec les relations
\begin{equation}\label{relWeil} jzj^{-1}=\bar z , \, (z\in \bbC^\times), \quad j^2=-1.
\end{equation}

On note $\Triv_{W_\bbR}$ le caractère trivial de $W_\bbR$ et $\sgn_{W_\bbR}$ la caractère quadratique, trivial sur $W_\bbC=\bbC^\times$ et valant $-1$ en l'élément $j$.

\subsection{ Plongements    $ \Sp(2p,\mathbb{C})\times W_\mathbb{R} \longrightarrow {}^L \U_n$ et
 $\SO(2p,\mathbb{C})\rtimes_p W_\mathbb{R} \longrightarrow {}^L\U_n$}\label{PlClU}
On note ${}^L \U_{n}=\GL(n,\bbC)\rtimes W_\bbR$ le $L$-groupe d'un groupe unitaire de rang $n$. Donnons une réalisation explicite de ce $L$-groupe.
Il s'agit de  définir l'action de $W_\bbR$ sur $\GL(n,\bbC)$. Cette action se factorise par la projection   de $W_\bbR$ sur 
$\mathrm{Gal}(\bbC/\bbR)$ et l'on choisit l'action de l'élément non trivial de ce groupe de Galois préservant
l'épinglage usuel, à savoir $g\mapsto w_n({}^tg^{-1})w_n^{-1}$, où 
$w_n=\begin{pmatrix}  \ldots & &&\ldots &1\\
\ldots &&& -1&  \\
&&\iddots&&\\
&1&\ldots &&\\
(-1)^{n-1}&&\ldots&&
\end{pmatrix}$.
Remarquons que ${}^t w_n=w_n^{-1}=w_n $ si $n$ est impair et $ {}^t w_n=w_n^{-1}=-w_n$ si $n$ est pair.

De même, donnons une réalisation explicite  du $L$-groupe 
d'un groupe spécial orthogonal pair ayant une série discrète, c'est-à-dire un $\SO(p,q)$, $p+q=2n$, 
$p$ et $q$ pairs. Notons 
$\SO(2n,\bbC)\rtimes_n W_\bbR$
 ce $L$-groupe.
Cette action se factorise par la projection   de $W_\bbR$ sur 
$\mathrm{Gal}(\bbC/\bbR)$,   l'action de l'élément non trivial de ce groupe de Galois étant
 $g\mapsto Tg T^{-1}$, où $T=T_{2n}$ est la matrice diagonale
ayant des coefficients $1$ et $-1$ qui alternent. Cette conjugaison est un automorphisme intérieur si $n$ est pair, et extérieur si $n$ est impair.
Remarquons que que cette construction donne naturellement un morphisme 
\begin{equation}\label{xiso} \SO(2n,\bbC)\rtimes_n W_\bbR \longrightarrow \Or(2n,\bbC) . \end{equation}
Pour tout entier $m$,  et tout nombre complexe $z$, on a posé 
$\chi_m(z) = \left( \frac{z}{\bar z} \right)^{\frac{m}{2}} = z^m ( z\bar z)^{-\frac{m}{2}}$.

\begin{prop}\label{PSOU}
 Fixons des entiers $n$, $p$, $a$, $a'$  et $b$ tels que $n\geq 2p$,  $b\cong 0\mod 2$, $a \cong n \mod 2$ et $a' \cong n-1 \mod 2$. 

(i) Il existe une inclusion de $L$-groupes 
$\xi_{Sp}:  \,  \Sp(2p,\mathbb{C})\times W_\mathbb{R} \longrightarrow {}^L \U_n $
dont la restriction  à $W_\mathbb{C}=\mathbb{C}^\times $ est donnée à conjugaison près par :
$
z\mapsto (\underbrace{\chi_{a}(z), \cdots, \chi_{a}(z)}_{2p}, \underbrace{\chi_{b}(z), \cdots, \chi_{b}(z)}_{n-2p})$.

(ii) Il existe une inclusion de $L$-groupes 
$ \xi_{SO} :\,  \SO(2p,\mathbb{C})\rtimes_p W_\mathbb{R} \longrightarrow {}^L\U_n$
 dont la restriction  à $W_\mathbb{C}=\mathbb{C}^\times $  est donnée à conjugaison près par :
$
z\mapsto (\underbrace{\chi_{a'}(z), \cdots, \chi_{a'}(z)}_{2p}, \underbrace{\chi_{b}(z), \cdots, \chi_{b}(z)}_{n-2p})$.

\end{prop}

\dem Le $(i)$ est classique, on remarque que la matrice $w_{2p}$ qui définit la réalisation de ${}^L\U_{2p}$ donnée ci-dessus
est antisymétrique.  L'automorphisme $g\mapsto w_{2p}{}^tg^{-1}w_{2p}^{-1}$ de $\GL(2p,\bbC)$ a donc pour groupe des points fixes
un groupe symplectique, que nous choisissons ici comme réalisation du groupe $\Sp(2p,\bbC)$. On prolonge cette inclusion de 
$\Sp(2p,\bbC)$ dans $\GL(2p,\bbC)$ par l'identité sur le facteur $W_\bbR$ pour obtenir une injection 
$\Sp(2p, \mathbb{C}) \hookrightarrow {}^L  \U_{2p}$. 
Ensuite, on  envoie 
${}^L\U_{2p}$ dans ${}^L \left( \U_{2p} \times \U_{n-2p}\right)$ naturellement 
 sur le premier facteur puis ${}^L \left( \U_{2p}\times \U_{n-2p}\right) $ dans ${}^L\U_n$ en utilisant les entiers $a$ et $b$ 
 ({\sl cf.} \cite{waldspurgertransfert}).

Pour $(ii)$, commençons par le cas où  $n$ est impair. On réalise  $\Or(2p,\mathbb{C})$ 
comme un sous-groupe de $\GL(n,\bbC)$ de la manière suivante : considérons la matrice 
$J_{n,p}=
   \begin{pmatrix}
 0    &\rvline & 0    &\rvline &  w_p \\
\hline 
 0    &\rvline & I_{n-2p}    &\rvline & 0 \\
 \hline
{}^t  w_p  &\rvline &  0  &\rvline & 0 
  \end{pmatrix}
$ qui est symétrique dans $\GL(n,\bbC)$ et réalisons $\Or(2p,\bbC)$ comme le sous-groupe des matrices 
$g$ de   $\GL(n,\bbC)$ de la forme   
$g=\begin{pmatrix}
A &\rvline & 0    &\rvline & B \\
\hline 
 0    &\rvline & I_{n-2p}    &\rvline & 0 \\
 \hline
C  &\rvline &  0  &\rvline & D
  \end{pmatrix}$ telles que $g=J_{n,p} {}^tg^{-1}J_{n,p}^{-1}$.
Comme $w_n=
   \begin{pmatrix}
 0   &\rvline & 0    &\rvline &  w_p \\
\hline 
0  &\rvline &(-1)^p w_{n-2p}    &\rvline & 0 \\
 \hline
{}^t  w_p  &\rvline & 0  &\rvline & 0
  \end{pmatrix}$, il est clair que $\Or(2p,\bbC)$ ainsi réalisé est inclu dans le groupe des  points fixes 
de l'automorphisme $g\mapsto w_{n}{}^tg^{-1}w_{n}^{-1}$ de $\GL(n,\bbC)$.
En particulier $\Or(2p,\bbC)$ commute à ${}^L\U_{n-2p}=\GL(n-2p,\bbC)\rtimes W_\bbR$ réalisé  dans le bloc central, et l'on a donc une inclusion 
\[\iota:\,  \Or(2p,\bbC) \times {}^L\U_{n-2p}  \hookrightarrow  \GL(n,\bbC) \rtimes W_\bbR. \] 
L'entier $b-a'$ est pair et définit un isomorphisme 
$ \xi_{b-a'}: \, {}^L\U_{n-2p}\rightarrow {}^L\U_{n-2p}$
 et de même pour $a'$ et ${}^L\U_n$
 qui donne un isomorphisme  
 $\xi_{a'}: \, {}^L\U_{n}\rightarrow {}^L\U_{n} $
(par torsion  respectivement par le caractère de $W_\mathbb{R}$ qui paramètre le caractère $\det^{\frac{b-a'}{2}}$ de $\U_{n-2p}$
et $\det^{\frac{a'}{2}}$ de  $ {}^L\U_{n}$, voir  \cite{waldspurgertransfert}).
On  étend ainsi  l'inclusion $\Or(2p,\bbC)\hookrightarrow \GL(n,\bbC)$ en tordant $\iota$ à la source et au but : 
\[  \Or(2p,\bbC)\times {}^L\U_{n-2p} \stackrel{  \Id\times \xi_{b-a'} }{\longrightarrow }  \Or(2p,\bbC) \times {}^L\U_{n-2p}  \longrightarrow  {}^L\U_n 
 \stackrel{  \xi_{a'} }{\longrightarrow} {}^L\U_n  , \]
On compose ensuite avec (\ref{xiso}) pour obtenir le 
morphisme $\xi_{SO}$ de  $(ii)$.

Si $n$ est pair, la construction précédente doit être légèrement modifiée. A la place de $J_{n,p}$ on prend la matrice 
$J'_{n,p}$ qui est antidiagonale avec tous les coefficients égaux à $1$. On réalise alors 
$\Or(2p,\bbC)$ dans $\GL(n,\bbC)$ avec la même recette,   $J'_{n,p}$ remplaçant $J_{n,p}$.
Posons $T'_{2p}= 
   \begin{pmatrix}
 0   &\rvline & 0    &\rvline &  w_p \\
\hline 
0  &\rvline & I_{n-2p}    &\rvline & 0 \\
 \hline
- {}^t  w_p  &\rvline & 0  &\rvline & 0
  \end{pmatrix}$. 
 Définissons $\xi_{SO}$ sur le facteur $\bbC^\times$ de $W_\bbR$ par  
\[ z\mapsto  ( (\underbrace{\chi_{a'}(z), \cdots, \chi_{a'}(z)}_{p}, \underbrace{\chi_{b}(z), \cdots, \chi_{b}(z)}_{n-2p},\underbrace{\chi_{a'}(z), \cdots, \chi_{a'}(z)}_{p}), z)
\in \GL(n,\bbC)\rtimes W_\bbR, 
\]
et sur l'élément $j$ de $W_\bbR$ par $ \xi_{SO}(j)=(T'_{2p} w_n,j ).$ On vérifie que $ \xi_{SO}(j)^2=(1,-1)$, et pour tout $z\in \bbC^\times$,
$  \xi_{SO}(j)  \xi_{SO}(z)  \xi_{SO}(j)^{-1} =\xi_{SO}(\bar z) $. 
Ainsi on a bien un morphisme de $W_\bbR$ dans ${}^L\U_n$.
Il reste à vérifier que l'inclusion, notons-là $\iota$,   donnée ci-dessus de $\SO(2p,\bbC)$ dans $\GL(n,\bbC)$ et le morphisme 
 $\xi_{SO}$  défini sur $W_\bbR$ sont compatibles aux relations dans $\SO(2p,\bbC)\rtimes_p W_\bbR$ et ${}^L\U_n$.
 D'une part, $\xi_{SO}(\bbC^\times )$ et $\iota(g)$ commutent. D'autre part, on a dans 
 $\SO(2p,\bbC)\rtimes_p W_\bbR$,
\[   (1,j)(g,1)(1,j)^{-1}= (T_{2p} gT_{2p}^{-1} , 1) , \] 
 et ceci s'envoie donc sur $(\iota (T_{2p}gT_{2p}^{-1}), 1)$. 
 De plus,  
\begin{align*} & \xi_{SO}(j) (\iota(g),1) \xi_{SO}(j)^{-1}=  (T'_{2p} w_n,j ) (\iota(g),1) (T'_{2p} w_n,j )^{-1}   
=(T'_{2p} w_n   w_n {}^t\iota(g)^{-1} w_n^{-1} (T'_{2p} w_n)^{-1},1)\\
=&(T'_{2p} {}^t\iota(g)^{-1} (T'_{2p})^{-1},1)=  
(T'_{2p}  (J'_{n,p})^{-1}\iota(g) J'_{n,p} (T'_{2p})^{-1},1).  \end{align*}
 En effet $  {}^t\iota(g)^{-1}= (J'_{n,p})^{-1}\iota(g) J'_{n,p} $, c'est la réalisation choisie de $\SO(2p,\bbC)$.
Il s'agit donc de vérifier que $\iota (T_{2p}gT_{2p}^{-1})=T'_{2p}  (J'_{n,p})^{-1}\iota(g) J'_{n,p} (T'_{2p})^{-1}$. Or 
$T'_{2p}  (J'_{n,p})^{-1}$ est la matrice diagonale dont les coefficients sont 
$  (\underbrace{1,-1,\ldots ,(-1)^{p-1}   }_p , \underbrace{1,\ldots, 1}_{n-2p},  \underbrace{ (-1)^{p-1}  , (-1)^{p-2}, \ldots,-1,-1}_p) $
et $\iota$ entrelace bien l'action par conjugaison de $T_{2p}$ et de cette matrice.\qed 

\bigskip

\subsection{Commutant de $ R[2] \oplus \ldots \oplus R[2] $  dans ${}^L\U_{2n}$ }\label{CR2}

Plongeons  $\GL(n,\bbC)\times \SL(2,\bbC)$ dans $\GL(2n,\bbC)$ en prenant le produit tensoriel des 
représentations naturelles de ces deux groupes, soit 
\begin{equation*} \label{iota} \iota: \GL(n,\bbC)  \times \SL(2,\bbC) \rightarrow \GL(2n,\bbC), \quad \left( g, \begin{pmatrix}
a&b\\c&d\end{pmatrix}\right)
\mapsto    \begin{pmatrix}
 a g   &\rvline &  bg \\
\hline 
cg &\rvline & dg \end{pmatrix}.\end{equation*}

Calculons le commutant  $\caG$ de $\iota(\SL(2,\bbC))$ dans ${}^L\U_{2n}$. Le commutant dans $\GL(2n,\bbC)$ est  exactement $\iota(\GL(2n,\bbC))$.
On a
\[w_{2n} {}^t \iota \left( {\begin{pmatrix} a&b\\c&d\end{pmatrix}} \right)^{-1}w_{2n}^{-1}= \begin{pmatrix}
aI_n& (-1)^{n-1}bI_n\\(-1)^{n-1}cI_n&dI_n\end{pmatrix}\]
et ainsi si $n$ est impair, 
$ \caG=  \{ (\iota(g),w)\in {}^L\U_{2n} , \; g\in \GL(n,\bbC), \; w\in W_\bbR \}$. 
Ce groupe est  le produit semi-direct de $\iota(\GL(n,\bbC))$ et de $W_\bbR$, où $W_\bbC=\bbC^\times \subset W_\bbR$ agit 
trivialement sur $\iota(\GL(n,\bbC))$ et $W_\bbR\setminus W_\bbC$ agit par l'automorphisme $\iota(g)\mapsto w_{2n}{}^t\iota(g)^{-1}w_{2n}^{-1}$.
On a un isomorphisme 
\begin{equation}\label{xiuplus} (n \text{ impair}) \quad 
\xi_+: \, {}^L\U_n\longrightarrow \caG\subset {}^L\U_{2n}, \quad (g,w)\mapsto  \left( \begin{pmatrix}
 g   &\rvline &0   & \\
\hline 
0&\rvline & g \end{pmatrix},w\right).
\end{equation}
Ceci est bien défini car pour tout $g\in \GL(n,\bbC)$, $ w_{2n} {}^t \iota(g)^{-1} w_{2n}^{-1}  =\iota(w_n {}^t {g}^{-1}w_n^{-1}) $. 
En effet, comme $\iota({}^tg^{-1}) ={}^t\iota(g)^{-1}$ pour tout $g\in \GL(n,\bbC)$, il reste à vérifier que $\Ad(w_{2n}^{-1}\iota(w_n ))$
agit trivialement sur l'image de $\iota$, ce qui est le cas, car  $w_{2n}^{-1}\iota(w_n )= \begin{pmatrix}
 0   &\rvline &(-1)^{n}I_n   & \\
\hline 
I_n&\rvline & 0  \end{pmatrix}$.

Si $n$ est pair 
$\caG= \left\{ (\iota(g)\begin{pmatrix} -I_n&0\\0&I_n\end{pmatrix},w)\in {}^L\U_{2n} , \; g\in \GL(n,\bbC), \; w\in W_\bbR \right\}$. 
On a encore un isomorphisme entre ${}^L\U_n$ et $\caG$ : 
\begin{equation}\label{xiumoins1}  (n \text{ pair})
   \quad \xi_-: \, {}^L\U_n\longrightarrow  \caG\subset  {}^L\U_{2n}, \quad (g,w)\mapsto  (\iota(g)c(w),w), \end{equation}
   où $c:\, W_\bbR\rightarrow \GL(2n,\bbC)$ est défini par 
$ c(z)=\left(\frac{z}{\bar z}\right)^{\frac{1}{2}}I_{2n}, \; (z\in \bbC)$, $c(j)= \begin{pmatrix}
  -I_n  &\rvline &  0  & \\
\hline 
  0&\rvline & I_n \end{pmatrix}$, $c(zj)=c(z)c(j)$.

\section{Le commutant du $\SL(2)$ et le plongement  ${\varphi_\scrX}:\,  {}^LG_\scrX\times \SL(2,\bbC)\rightarrow {}^LG$ .} \label{commutantSL2}

Reprenons les résultats obtenus dans la section \ref{SV}.
Soit $\pi$ une série discrète d'un espace $\scrX$ de la liste, et soit $\psi_\pi$ le paramètre d'Arthur que nous lui associons. 
  D'après le corollaire \ref{SL2}, la restriction de $\psi_\pi$ à  $\SL(2,\bbC)$ est (à conjugaison près) le morphisme 
  de Jacobson-Morosov $\eta_\scrX$  associé à l'orbite unipotente régulière du  
 sous-groupe  de Levi $\mathbf L^\vee_\scrX$ de $G^\vee$. 
Ce groupe (ou plutôt son algèbre de Lie  $\check  \frl _\scrX $, ce qui le détermine)  est  donné dans la table 1. 
Nous allons calculer le commutant $\caG$  de cette image de $\SL(2,\bbC)$ dans ${}^LG$. Ceci nous permet de vérifier que le groupe 
de Langlands dual ${}^LG_\scrX$ prévu par  \cite{SV} et \cite{KS} est bien le bon, et à 
 terme, de déterminer   le  morphisme  $\varphi_\scrX$ de (\ref{varphiscrX}) 
  (ou de manière équivalente, puisque la restriction au facteur $\SL(2,\bbC)$ est déjà donnée, le morphisme 
 $\bar \varphi_\scrX$ de (\ref{factpar2})). En effet l'image par ce morphisme de ${}^LG_\scrX$ doit  être contenue  dans $\caG$.  
  
\medskip

Nous commençons par les cas les plus simples, ceux où le commutant $\caG$ est isomorphe au groupe de Langlands dual ${}^LG_\scrX$ prévu par 
\cite{SV} et \cite{KS}.

\begin{prop}  \label{Gcastressimples}  Dans les cas 2, 5,  6,  9 avec $n=p+q$ pair, 10 avec $n$ pair, 11,  12 et 13 de la table \ref{table1},  le commutant  
$\caG$ est isomorphe à  $^LG_\scrX$. Ceci fixe le morphisme $\bar \varphi_\scrX$ de
  (\ref{factpar2}).
\end{prop}

\dem {\bf Cas 2} :  $\U(p,q)/\Or(p,q) $,  $n=p+q$.  On a  $ {}^L G=\GL(n,\bbC)\rtimes W_\bbR$ et $ \mathbf L_\scrX^\vee=\GL(1,\bbC)^n$ qui  est un tore
 donc de $\SL(2)$ principal  trivial. 
Le $L$-groupe ${}^LG_\scrX$ est  égal à ${}^L G$ et 
   le morphisme (\ref{varphiscrX}) égal à l'identité des deux $L$-groupes.

\medskip

 {\bf Cas 5} :  $ \U(2p,2q)/\Sp(p,q)  $.  Posons $n=p+q$. On a  $ {}^LG=\GL(2n,\bbC)\rtimes W_\bbR$ et $  \mathbf L^\vee_\scrX=\GL(2,\bbC)^{n}$. 
Le commutant $\caG$  est isomorphe à 
${}^L\U_n=\GL(n,\bbC)\rtimes W_\bbR$. Les calculs sont donnés dans la section \ref{CR2} où  nous définissons 
un isomorphisme  entre ${}^L\U_n$  et $\caG$. La forme de cet isomorphisme
dépend de la parité de $n$, il est noté $\xi_-$  si  $n$ est  pair, et  $\xi_+$ 
si  $n$ est impair.
Ainsi  ${}^LG_\scrX=\GL(n,\bbC)\rtimes W_\bbR={}^L\U_n$, et 
le  morphisme $\bar \varphi_\scrX$ de (\ref{factpar2}) est donné par 
$\xi_-$ ou $\xi_+$ selon la parité de $n$.

\medskip 

 {\bf Cas 6} :   $ \U(n,n)/\Sp(2n,\bbR)  $.  Ce cas est identique au précédent

  \medskip
 
  {\bf Cas 9} :   $\SO(2p,2q)/\U(p,q)$ avec   $n=p+q$ pair.   L'espace symétrique $\scrX$ admet des séries discrètes
 sauf dans le cas où $p$ et $q$ sont tous les deux impairs, donc si $n=p+q$ est pair, $p$ et $q$ sont pairs.
 Le groupe $G$  est une  forme intérieure de $\SO(n,n)$. On a alors 
 $  {}^LG =\SO(2n,\bbC)\times W_\bbR$ et $\mathbf L^\vee_\scrX=\GL(2,\bbC)^{\frac{n}{2}}$.
 Le  commutant  du $\SL(2)$ principal de  $\mathbf L^\vee_\scrX$ dans $\GL(2n,\bbC)$ dans $\SO(2n,\bbC)$ est  
$\Sp(n,\bbC)$  (\cite{CMcG}, Thm. 6.1.3). Ceci fixe  ${}^LG_\scrX= \Sp(n,\bbC)\times W_\bbR$.
Le commutant $\caG$ est donc  isomorphe à ${}^LG_\scrX= \Sp(n,\bbC)\times W_\bbR$ et 
ceci fixe le morphisme $\bar \varphi_\scrX$ de
  (\ref{factpar2}).

 \medskip

  {\bf Cas 10}   :   $\SO(n,n)/\GL(n,\bbR)$ avec $n=p+q$ pair. Ce cas est similaire au cas précédent. Ici ${}^LG=\SO(2n,\bbC)\times W_\bbR$, et 
 $ \mathbf L^\vee_\scrX=\GL(2,\bbC)^{\frac{n}{2}}$.
  Le commutant $\caG$ du $\SL(2)$-principal de  $ \mathbf L^\vee_\scrX$ 
  est $\caG= \Sp(n,\bbC)\times W_\bbR$.  Ainsi  $ {}^LG_\scrX= \Sp(n  , \bbC)\times W_\bbR$ et et 
ceci fixe le morphisme $\bar \varphi_\scrX$ de
  (\ref{factpar2}).

 \medskip

  {\bf Cas 11}:   $\scrX= \Sp(2n,\bbR)/\Sp(2p,\bbR) \times \Sp(2(n-p),\bbR) $,  $2p\leq n$.
 On a ici 
   \[{}^LG=\SO(2n+1,\bbC)\times W_\bbR,  \quad  \mathbf L^\vee_\scrX=  \GL(2,\bbC)^{p}\times   \SO(2n-4p+1,\bbC).\]
Le commutant dans  $\SO(2n+1,\bbC)$ du $\SL(2)$ principal de  $\mathbf L^\vee_\scrX$ est 
$\caG=\Sp(2p,\bbC)\times \{  I_{2n-4p+1}\}$ ({\sl cf.} \cite{CMcG} Thm. 6.1.3).
Ainsi   $ {}^LG_\scrX= \Sp(2p, \bbC) \times W_\bbR$ et 
ceci fixe le morphisme $\bar \varphi_\scrX$ de
  (\ref{factpar2}).

\medskip 

 {\bf Cas 12}:  $\scrX= \Sp(4n,\bbR)/\Sp(2n,\bbC)  $. 
Ce cas est analogue au précédent, en remplaçant $2n$ par $4n$ et $p$ par $n$.
\medskip 

 {\bf Cas 13}:  $\scrX=\Sp(2n,\bbR)/\GL(n,\bbR) $.  
 On a 
  ${}^LG= \SO(2n+1, \bbC) \times W_\bbR$ et $  \mathbf L^\vee_\scrX=     \GL(1,\bbC)^{n}$ qui est un tore, donc de   $\SL(2)$ principal trivial. 
Le morphisme  $\varphi_\scrX$ de (\ref{varphiscrX}) est  l'identité des $L$-groupes:  ${}^LG_\scrX= \SO(2n+1, \bbC) \times W_\bbR={}^LG$. \qed

\begin{prop}   \label{Gcassimples}   Dans les cas 7 et 8, le commutant  $\caG$  de l'image de 
$\SL(2,\mathbb{C})$ dans $^LG$ est isomorphe à 
 $^LG_\scrX\times \{\pm 1\}$. Le morphisme $\varphi_\scrX$ de (\ref{varphiscrX}) est alors fixé 
 à une torsion près dans $\{\pm 1\}$.
\end{prop}

\dem {\bf Cas 7} :  $\SO(p,q)/\SO(r,s)\times \SO(r',s') $, $p+q=2n+1$. 
On pose $m=r+s$.  On a ici 
${}^L G=\Sp(2n,\bbC)\times W_\bbR$ et  $ \mathbf L^\vee_\scrX=  \GL(1,\bbC)^{m}\times  
\Sp(2(n-m) ,\bbC)$. 
Le commutant du $\SL(2)$ principal de  $\mathbf L^\vee_\scrX$  dans $\Sp(2n,\bbC)$ est  
$\Sp(2m,\bbC)\times \{ \pm I_{2n-2m}\}$ (\cite{CMcG}, Thm. 6.3.1).
On a donc ${}^LG_\scrX=\Sp(2m,\bbC)\times W_\bbR$ et 
$\caG$  est isomorphe à ${}^LG_\scrX\times \{ \pm1 \}$.

Ceci donne de manière évidente  un morphisme :
\begin{equation}\label{pho7}  \bar \varphi_\scrX: \; {}^LG_\scrX= \Sp(2m,\bbC)\times W_\bbR 
 \longrightarrow {}^LG=\Sp(2n,\bbC)\times W_\bbR
  \end{equation}
dont  la restriction à $W_\bbR$ est l'identité de ce facteur. 
  
  On a aussi la possibilité de tordre ce morphisme par $\sgn_{W_\bbR}$   à valeurs dans le facteur $\{\pm 1\}$ de $\caG$, c'est-à-dire que 
  l'on définit un morphisme 
  \begin{equation}\label{pho7'}  \bar \varphi'_\scrX: \; {}^LG_\scrX = \Sp(2m,\bbC)\times W_\bbR 
 \longrightarrow {}^LG=\Sp(2n,\bbC)\times W_\bbR
  \end{equation}
avec même restriction à   $\Sp(2m,\bbC)$ que $\varphi$, la restriction   
  de $\bar \varphi'_\scrX$ à $W_\bbR$ étant
  \[w\mapsto (  \sgn_{W_\bbR},w)\in  \{ \pm I_{2n-2m}\} \times W_\bbR\subset \caG\times W_\bbR.\]

\medskip

{\bf Cas 8} :   $\SO(p,q)/\SO(r,s)\times \SO(r',s') $, $p+q=2n  $. On pose $m=r+s$. 
 Comme $p+q=2n$, $p$ et $q$ sont de même parité et le $L$-groupe
 de $\SO(p,q)$ est $\SO(2n,\bbC)\times W_\bbR$ si cette parité est aussi celle de $n$, et 
$\SO(2n,\bbC)\rtimes W_\bbR$ si cette parité est celle de $n+1$, l'action de $W_\bbR$ étant induite par l'automorphisme extérieur de $\SO(2n,\bbC)$
en fixant un épinglage.   On a ici $ \mathbf L^\vee_\scrX=  \GL(1,\bbC)^{m}\times  \SO(2(n-m),\bbC)$.

Rappelons que l'orbite
unipotente principale de $\SO(2(n-m),\bbC)$ est associée à la partition $(1,2(n-m)-1)$ de $2n$,  
 et son commutant  dans $\SO(2n,\bbC)$ est 
$\mathbf S(\Or(2m+1,\bbC)\times \{ \pm I_{2n-2m-1}   \})$ ({\sl cf. \cite{CMcG}},Thm. 6.3.1). 
On a donc ${}^LG_\scrX=\SO(2m+1,\bbC)\times W_\bbR$.

Si ${}^LG=\SO(2n,\bbC)\times W_\bbR$,  on voit immédiatement que $\caG$ est  isomorphe 
à ${}^LG_\scrX\times \{\pm 1\}$, 
et si  ${}^LG=\SO(2n,\bbC)\rtimes W_\bbR$,  comme l'image de $\SL(2,\bbC)$ est 
 dans le groupe des points fixes de $\SO(2n,\bbC)$ sous l'action de $W_\bbR$, 
  on a la même conclusion.
Comme en (\ref{pho7}) et (\ref{pho7'}), on définit 
\begin{equation}\label{pho8}  \bar \varphi_\scrX,\bar  \varphi'_\scrX: \; {}^LG_X = \Sp(2m,\bbC)\times W_\bbR 
 \longrightarrow {}^LG=\SO(2n,\bbC)\rtimes W_\bbR.
  \end{equation}
où $\bar \varphi_\scrX$ est donnée par les inclusions évidentes et $\bar \varphi'_\scrX$ est obtenu en tordant $\bar \varphi_\scrX$ par le caractère 
  $\sgn_{W_\bbR}$   à valeurs dans le facteur $\{\pm 1\}$ de $\caG$. \qed

\begin{prop}    \label{Gcastsimples}     Dans les cas 9 et 10  avec  $n$ impair, le  commutant  $\caG$  de l'image de 
$\SL(2,\mathbb{C})$ dans $^LG$ est isomorphe à  $^LG_\scrX\times \SO(2,\bbC)$. 
Dans le  cas 9 avec $n$ impair,  ceci fixe le morphisme $\varphi_\scrX$ de
  (\ref{varphiscrX}). Dans le cas   10 avec $n$ impair, le morphisme $\varphi_\scrX$ de
  (\ref{varphiscrX})  est fixé  à une torsion près dans $\{\pm 1\}$. \end{prop}

\dem 
{\bf Cas 9} avec   $n=p+q$ impair.
Dans ce cas  $p$ et $q$ sont  de parité différente, le groupe $G$ est forme intérieure de $\SO(n+1,n-1)$. On a alors  
$ {}^LG =\SO(2n,\bbC)\rtimes W_\bbR$ et 
  $  \mathbf L^\vee_\scrX=\GL(2,\bbC)^{\frac{n-1}{2}} \times \SO(2,\bbC)$.
Le commutant $\caG$ du $\SL(2)$ principal  de  $\mathbf L^\vee_\scrX$  dans $\SO(2n,\bbC)$  est $\Sp(n-1,\bbC)\times \SO(2,\bbC)$
(\cite{CMcG}, Thm. 6.1.3).
On remarque que  ce $\SL(2)$ principal  est dans les points fixes de l'automorphisme extérieur  par lequel agit $W_\bbR$ 
sur  $\SO(2n,\bbC)$. On a alors $\caG=\left(\Sp(n-1,\bbC)\times \SO(2,\bbC)\right)\rtimes W_\bbR$ et  ${}^LG_\scrX= \Sp(n-1, \bbC)\times W_\bbR$. 

\begin{rmq}\label{notor} Le seul $L$-morphisme de $W_\bbR$ à valeurs dans $\SO(2,\bbC)\rtimes W_\bbR$ trivial sur $\bbR^\times_+$ est le morphisme trivial.
Le morphisme $\bar \varphi_\scrX$ est donc fixé, il n'y a pas possibilité  de le tordre par un morphisme à valeurs dans le facteur  
$\SO(2,\bbC)\rtimes W_\bbR$ de $\caG$.
\end{rmq}

{\bf Cas 10} avec $n$ impair :  $\SO(n,n)/\GL(n,\bbR)$. Ici ${}^LG=\SO(2n,\bbC)\times W_\bbR$, et 
 $    \mathbf L^\vee_\scrX=\GL(2,\bbC)^{\frac{n-1}{2}}\times \SO(2,\bbC)$ et  $\caG= \left(\Sp(n-1,\bbC) \times \SO(2,\bbC)\right) \times W_\bbR$.
 Ainsi $ {}^LG_\scrX= \Sp(n-1, \bbC)\times W_\bbR$. 
Les $L$-morphismes de $W_\bbR$ à valeurs dans $\SO(2,\bbC)\times W_\bbR$ triviaux sur $\bbR^\times_+$ sont $\Triv_{W_\bbR}$
 et $\sgn_{W_\bbR} $. On a donc un morphisme $\bar \varphi_\scrX$ défini par les inclusions évidentes, et 
  de manière similaire à (\ref{pho7}) et (\ref{pho7'}), on a 
  possibilité de tordre le morphisme 
$\bar \varphi_\scrX$  par le caractère signe de   $W_\bbR$ à valeurs dans le facteur $\SO(2,\bbC)$ du commutant pour obtenir un morphisme $\bar \varphi'_\scrX$.
\qed

\medskip

Les cas restants sont ceux où les racines sphériques de $\scrX$ ne sont pas toutes des racines. On les traite maintenant.

{\bf Cas 1} :  $\GL(n,\bbR)/\GL(p,\bbR)\times \GL(n-p,\bbR)$, $2p\leq n$.
On a ici :
    \[ {}^LG= \GL(n,\bbC)\times W_\bbR, \qquad \mathbf L^\vee_\scrX=  \GL(1,\bbC)^{2p}\times  
\GL(n-2p,\bbC).\]
 Le  commutant du $\SL(2)$ principal de  $\mathbf L^\vee_\scrX$dans $\GL(n,\bbC)$ est 
   $\GL(2p)\times Z(\GL(n-2p,\bbC))$ et son commutant dans ${}^LG$ est donc 
   $\caG=\left( \GL(2p,\bbC)\times Z(\GL(n-2p,\bbC))\right) \times W_\bbR$. 

\medskip 

 {\bf Cas 3} : $\U(p,q)/\U(r,s)\times \U(r',s')$.
Posons  $p+q=n$ et $m=r+s$. On a ici :
    $ {}^LG= \GL(n,\bbC)\rtimes W_\bbR, \; \mathbf L^\vee_\scrX=  \GL(1,\bbC)^{2m}\times  
\GL(n-2m,\bbC)$.
    Le  commutant du $\SL(2)$ principal de  $\mathbf L^\vee_\scrX$ dans $\GL(n,\bbC)$ est 
   $\GL(2m)\times Z(\GL(n-2m,\bbC))$.
    A conjugaison près, on s'arrange pour que l'image du $\SL(2)$ principal de  $\mathbf L^\vee_\scrX$
 soit stable par l'automorphisme $g\mapsto w_n{}^tg^{-1}w_n$ qui définit le $L$-groupe de $\U(n)$ (voir section \ref{Lplong}).   
    Le commutant  $\caG$  est alors 
   $ \caG=\left( \GL(2m)\times Z(\GL(n-2m,\bbC))  \right)\rtimes W_\bbR$.

\medskip

{\bf Cas 4} :   $\U(n,n)/\GL(n,\bbC)  $. On a ici ${}^LG= \GL(2n,\bbC)\rtimes W_\bbR, \, \mathbf L^\vee_\scrX=\GL(1,\bbC)^{2n}$.
Ainsi   $\mathbf L_\scrX^\vee$ est un tore et son $\SL(2)$ principal est trivial. 
Son commutant est donc  ${}^LG$.

\section{ Fin de la démonstration de la conjecture}\label{DEM}

Dans cette section, nous  finissons  la démonstration de la conjecture de 
Sakellaridis et Venkatesh commencée dans les sections  \ref{debut} et   \ref{commutantSL2}. 
Dans les cas couverts par la proposition \ref{Gcastressimples},  le commutant $\caG$   est 
isomorphe à $^LG_\scrX$, ce qui fixe le morphisme $\varphi_\scrX$  et  la conjecture est essentiellement tautologique, 
puisque le paramètre d'Arthur $\psi_\pi$  se factorise nécessairement par  $\varphi_\scrX$.  
 Il reste à vérifier que  le morphisme $\phi_\pi$, a priori tempéré,  est en fait discret.
On détaille les cas particuliers  pour  
montrer  les torsions éventuelles quand l'isomorphisme entre $^LG_\scrX$ et 
le commutant n'est pas l'identité sur $W_\mathbb{R}$.

{\bf Cas 2 et 13} : dans ces cas    $L$ est un tore  compact. Les représentations
$A_\frqqq(\pi_L)$ sont alors des séries discrètes de $G$. Comme on a identité des $L$-groupes ${}^LG_\scrX={}^LG$, 
les paramètres $\psi_\pi$ (resp.  (\ref{psicas2}) et  (\ref{psicas13}))  sont des paramètres de Langlands discrets et il n'y a rien de plus à démontrer.

{\bf Cas  5 et 6} : les paramètres sont de la forme 
$  \psi_\pi =\bigoplus_{i=1}^n  (\chi_{m_i} \boxtimes R[2])$ 
({\sl cf.}(\ref{psicas5}) et (\ref{psicas6})).
Le paramètre $\phi_d= \bigoplus_{i=1}^n  \chi_{m_i} $ est un paramètre de Langlands discret de $\U_n$ si $n$ est impair et 
$\phi_d=\bigoplus_{i=1}^n  \chi_{m_i-1} $  est   un paramètre de Langlands discret de $\U_n$ si $n$ est pair. 

On a vu dans la section \ref{commutantSL2} que ${}^LG_\scrX$ est isomorphe  à ${}^L\U_n$.
On envoie  alors  ${}^LG_\scrX=  {}^L\U_n$  dans ${}^LG$ par le morphisme $\xi_+$ (si $n$ est impair)
et $\xi_-$ (si $n$ est pair), définis respectivement en (\ref{xiuplus}) et (\ref{xiumoins1}), ce qui définit le morphisme $\varphi_\scrX$ de (\ref{varphiscrX}).
On a alors la propriété de factorisation  $\psi_\pi =\varphi_\scrX \circ \phi_d$.
Remarquons que la torsion  sur $W_\bbR$ utilisée dans la définition du  morphisme $\xi_-$ a bien l'effet voulu 
sur le caractère infinitésimal, en le translatant de $1/2$, ou autrement dit, les  $\chi_{m_i-1}$
deviennent bien  $\chi_{m_i}$ après composition par $\varphi_\chi$.

{\bf Cas  11 et 12} :   les paramètres  $\psi_\pi$ sont respectivement  de la forme   (\ref{psicas5}) et (\ref{psicas6}).
Le paramètre discret $\phi_d=\bigoplus_{i=1}^{p}  \left(\delta\left( \frac{m_i}{2},  - \frac{m_i}{2}\right) \boxtimes R[2] \right)$, avec les $m_i$ impairs,  
est à valeurs dans $\Sp(2p,\bbC)\times W_\bbR= {}^LG_\scrX$ et on a la factorisation  $\psi_\pi =\varphi_\scrX \circ \phi_d$
(dans le cas 12, remplacer $p$ par $n$).

{\bf Cas 9 et 10 avec $n$ pair} :  les paramètres  $\psi_\pi$ sont  de la forme 
$ \psi_\pi=\bigoplus_{i=1}^{\frac{n}{2}}  \left(\delta\left( \frac{m_i}{2},  - \frac{m_i}{2}\right) \boxtimes R[2] \right)  $
avec les $m_i$ impairs.
({\sl cf.}(\ref{psicas9a}) et (\ref{psicas10a})).
Le paramètre 
$ \phi_d=\bigoplus_{i=1}^{ \frac{n}{2} }  \delta\left( \frac{m_i}{2},  - \frac{m_i}{2}\right) $
est de bonne parité pour un groupe orthogonal impair de rang $ \frac{n}{2}  $, dont le $L$-groupe est 
$\Sp(n  ,\bbC)\times W_\bbR={}^LG_\scrX$.
Avec le morphisme $\varphi_\scrX$ défini dans ce cas dans la section \ref{commutantSL2}, on a bien la  propriété de factorisation 
$\psi_\pi=\varphi_\scrX \circ \phi_d$.

On passe maintenant au cas couverts par la proposition \ref{Gcastsimples}.

{\bf Cas 9 et 10 avec $n$ impair} :  les paramètres  $\psi_\pi$ sont respectivement  de la forme 
\[ \psi_\pi=\bigoplus_{i=1}^{\frac{n-1}{2}}  \left(\delta\left( \frac{m_i}{2},  - \frac{m_i}{2}\right) \boxtimes R[2] \right) 
 \oplus \left(\sgn_{W_\bbR} \boxtimes R[1]\right) \oplus
 \left( \Triv_{W_\bbR} \boxtimes R[1]\right) 
  \]
  \[ \psi_\pi=\bigoplus_{i=1}^{\frac{n-1}{2}}  \left(\delta\left( \frac{m_i}{2},  - \frac{m_i}{2}\right) \boxtimes R[2] \right) 
 \oplus \left(\Triv_{W_\bbR} \boxtimes R[1]\right) \oplus
 \left( \Triv_{W_\bbR} \boxtimes R[1]\right) 
  \]
  avec les $m_i$ impairs. ({\sl cf.} (\ref{psicas9b}) et (\ref{psicas10b})).  On pose donc 
$ \phi_d=\bigoplus_{i=1}^{ \frac{n-1}{2} }  \delta\left( \frac{m_i}{2},  - \frac{m_i}{2}\right) $.
C'est un paramètre  de bonne parité pour un groupe orthogonal impair de rang $ \frac{n-1}{2} $, dont le $L$-groupe est 
$\Sp(\frac{n-1}{2} ,\bbC)\times W_\bbR={}^LG_\scrX$.
Avec les morphismes $\varphi_\scrX$ défini dans ce cas dans la section \ref{commutantSL2} (dans le cas 9, le morphisme
$\varphi_\scrX$ ést déjà fixé, {\sl cf.} remarque \ref{notor}, dans le cas 10, on prend le morphisme évident    $\varphi_\scrX$ sans la torsion sur le facteur $W_\bbR$).
On a bien la  propriété de factorisation  $\psi_\pi=\varphi_\scrX \circ \phi_d$.

\medskip

On continue par les  cas couverts par la proposition \ref{Gcassimples}.

 {\bf Cas  7}    : les paramètres  $\psi_\pi$ sont    de la forme    (\ref{psicas7}), c'est-à-dire
  \begin{equation*}
  \psi=\bigoplus_{i=1}^{r+s}  \left(\delta\left( \frac{m_i}{2},  - \frac{m_i}{2}\right) \boxtimes R[2] \right)  \oplus \left(\Triv_{W_\bbR} \boxtimes R[2(n-r-s)]  \right)    \end{equation*}
où  $m_i$ sont des entiers impairs et $\epsilon\in \{0, 1\}$.
 Le paramètre
  $ \phi_d=\bigoplus_{i=1}^{r+s  }  \delta\left( \frac{m_i}{2},  - \frac{m_i}{2}\right) $
 est un paramètre discret pour $\SO(2(r+s)+1)$, dont le $L$-groupe est $\Sp(2(r+s),\bbC)\times W_\bbR={}^LG_\scrX$.
 Le paramètre $\psi_\pi$ se factorise donc en $\psi_\pi=\varphi_\scrX \circ \phi_d$, où pour $\varphi_\scrX$ nous prenons le  morphisme défini
 en  (\ref{pho7}). 
 
\medskip

{\bf Cas 8}:   les paramètres  $\psi_\pi$ sont    de la forme    (\ref{psicas8}), c'est-à-dire   \begin{equation*}\label{psisop}
   \psi=\bigoplus_{i=1}^{r+s}  \left(\delta\left( \frac{m_i}{2},  - \frac{m_i}{2}\right) \boxtimes R[2] \right)  \oplus \left(\eta_1 \boxtimes R[2(n-r-s)-1] \right)
 \oplus \left(\eta_2 \boxtimes R[1] \right) ,   \end{equation*}
où  $m_i$ sont des entiers pairs et $\eta_1,\eta_2 \in \{\sgn_{W_\bbR}, \Triv_{W_\bbR}\}$ avec 
 $\eta_1\eta_2=\Triv_{W_\bbR} $ si $p$ et $q$ sont de la parité de $n$ et  $\eta_1\eta_2=\sgn_{W_\bbR} $ si $p$
et $q$ sont de la parité de $n-1$.

Le paramètre
  $ \phi_d=\bigoplus_{i=1}^{r+s  }  \delta\left( \frac{m_i}{2},  - \frac{m_i}{2}\right) $
 est un paramètre discret pour un groupe symplectique  de rang $r+s$, dont le $L$-groupe est $\SO(2(r+s)+1,\bbC)\times W_\bbR={}^LG_\scrX$.
 Le paramètre $\psi$ se factorise donc en $\psi_\pi=\varphi_\scrX \circ \phi_d$, où pour $\varphi_\scrX$ est le morphisme défini dans la section 
 \ref{commutantSL2} sans torsion par $\sgn_{W_\bbR}$ ({\sl cf.} (\ref{pho8})).

\begin{rmq} Nous conjecturons que les   paramètres (\ref{psicas7}), (\ref{psicas8})  et (\ref{psicas10b}) avec des parties unipotentes non triviales et les morphismes
$\varphi_\scrX$  de la section \ref{commutantSL2} tordus par  $\sgn_{W_\bbR}$ apparaissent dans le cadre des conjectures
 généralisées avec le   caractère $\chi$ de $H$ non trivial. Nous espérons revenir sur ces questions dans un article ultérieur.
\end{rmq}

Nous traitons maintenant  les cas restants.

{ \bf Cas 1}  :  les paramètres  $\psi_\pi$ sont    de la forme    (\ref{psigl}), c'est-dire 
\[\psi_\pi=\bigoplus_{i=1}^p  \left(  \delta\left(\frac{m_i}{2},-\frac{m_i}{2}\right)\boxtimes R[1]  \right)\oplus \left(  \Triv_{W_\bbR}\boxtimes R[n-2p]  \right).\]
avec les $m_i$ entiers impairs distincts.
On note $\bar \varphi_\scrX$ l'inclusion naturelle 
\[\bar \varphi_\scrX:\;  \Sp(2p,\mathbb{C})\times W_\mathbb{R} \longrightarrow 
\left(\GL(2p,\mathbb{C})\times \GL(n-2p,\mathbb{C})\right) \times W_\mathbb{R} \hookrightarrow  \GL(n,\mathbb{C}) \times W_\mathbb{R} \]
qui résulte de la décomposition de $\mathbb{C}^n$ en la somme de $\mathbb{C}^{2p}\oplus \mathbb{C}^{n-2p}$. 
On pose alors $\phi_d=\bigoplus_{i=1}^p  \left(  \delta\left(\frac{m_i}{2},-\frac{m_i}{2}\right)\boxtimes R[1]  \right)$. Comme les $m_i$ sont impairs,
 ceci est bien un paramètre à valeurs dans   ${}^LG_\scrX= \Sp(2p,\mathbb{C})\times W_\mathbb{R}$ et l'on a 
la propriété de factorisation  $\psi_\pi=\varphi_\scrX \circ \phi_d$.

{ \bf Cas 3}  :   les paramètres  $\psi_\pi$ sont    de la forme    (\ref{psicas3}), c'est-dire 
 \[ \psi_\pi=\bigoplus_{i=1}^{r+s}  ((\chi_{m_i}\oplus \chi_{-m_i}) \boxtimes R[1])   \oplus \left(\chi_0 \boxtimes R[n-2(r+s)] \right) \]

La représentation  $\chi_{m_i}\oplus \chi_{-m_i}$ de $\bbC^\times $ s'étend de manière unique en la représentation
  $\delta\left( \frac{m_i}{2},-\frac{m_i}{2} \right)$ de $W_\bbR$ cette représentation étant à valeurs dans $\SL(2,\bbC)$ si $m_i$ est impair
  (c'est-à-dire $n$ pair)
  et à valeurs dans $\Or(2,\bbC)$ si $m_i$ est pair  (c'est-à-dire $n$ impair).
  Posons 
  \[\phi_d = \bigoplus_{i=1}^{r+s} \delta\left( \frac{m_i}{2},-\frac{m_i}{2} \right) :\, W_\bbR\longrightarrow  \begin{cases}
  \Sp(2(r+s),\bbC)\times W_\bbR \text{ si }  n \text{ pair }\\ \Or(2(r+s),\bbC)\times W_\bbR \text{ si }  n \text{ impair }
  \end{cases}\]
  
Si $n$ est pair, on pose  ${}^LG_\scrX=\Sp(2(r+s),\bbC)\times W_\bbR$
comme prévu par \cite{KS}, et le morphisme $\bar \varphi_\scrX$ de (\ref{factpar2}) est  le morphisme $\xi_\Sp$ de 
 la proposition \ref{PSOU}, avec $a=b=0$. On a 
la propriété de factorisation voulue $\psi_\pi=\varphi_\scrX \circ \phi_d$.

Dans le cas où  $n$ est impair, une première façon d'obtenir un énoncé valide est de prendre  
 ${}^LG_\scrX=\SO(2(r+s)\bbC)\rtimes_{r+s}W_\bbR$ (voir  section \ref{Lplong} où l'on a donné
 la réalisation du $L$-groupe $\SO(2(r+s)\bbC)\rtimes_{r+s}W_\bbR$ d'un groupe orthogonal
  pair de rang $r+s$ admettant des séries discrètes).
On a  la propriété de factorisation voulue $\psi_\pi=\varphi_\scrX \circ \phi_d$, mais ${}^LG_\scrX$
n'est pas le groupe prévu dans \cite{SV} et \cite{KS}, car celui-ci a pour composante neutre $\Sp(2(r+s),\bbC)$.

  \medskip
  
  La seconde façon d'obtenir un énoncé valide dans le cas 3, $n$ impair, nous a été suggérée par Sakellaridis. Posons $k=r+s$ pour alléger les écritures.
  On garde la composante neutre $G^\vee_\scrX=\Sp(2k,\bbC)$ de ${}^LG_\scrX$ déterminée par Knop et Schalke, et 
  l'on prend pour  ${}^LG_\scrX$ un produit semi-direct $\Sp(2k,\mathbb{C})\rtimes W_\mathbb{R}$, qui n'est pas un $L$-groupe au sens de Langlands
  \cite{Langl}
  (le scindage $W_\bbR\rightarrow {}^LG_\scrX$ n'est pas \og  distingué \fg).
  Une telle extension de la notion de $L$-groupe a été introduite dans \cite{ABV} sous la terminologie
  \og E-group \fg. Les auteurs montrent que les paramètres de Langlands construits
  avec ces $L$-groupes paramètrent certaines  représentations projectives ({\sl loc. cite},  Thm. 10.4).
  On forme  le produit semi-direct  
  $\Sp(2k,\mathbb{C})\rtimes W_\mathbb{R}$ où $W_\mathbb{C}$ agit trivialement  et où $j\in W_\bbR$
  agit par conjugaison par un élément de $\GSp(2k,\mathbb{C})$ de norme symplectique $-1$.

\begin{lemme} Ce produit semi-direct est naturellement isomorphe à un sous-groupe de $\caG$, le  commutant de l'image de $\SL(2,\mathbb{C})$. 
\end{lemme}

\dem Le commutant $\caG$ a été calculé à la fin de la section \ref{commutantSL2}. 
On part d'une décomposition 
 $\mathbb{C}^n=\mathbb{C}^k\oplus \mathbb{C}^{n-2k}\oplus \mathbb{C}^k$,
de sorte que  l'image de $\SL(2,\mathbb{C})$ est un sous-groupe de $\GL(n-2k,\mathbb{C})$. 
On a ainsi $\GL(2k,\mathbb{C})\rtimes W_\mathbb{R}\subset \caG$. Ici le produit semi-direct est défini par
$(\forall z\in \mathbb{C}^\times,\, \forall g\in \GL(2k,\mathbb{C})),\;   z\cdot g=g$
et $
(\forall g\in \GL(2k,\mathbb{C})),\; j\cdot g=t  w_{2k} {}^tg^{-1} w_{2k}^{-1}t$,
où $w_{2k}$ est la matrice antisymétrique ayant des $1$ et $-1$ alternant sur l'antidiagonale   ({\sl  cf.}  notations de la section \ref{PlClU})    et $t$ est l'élément de
 $\GL(2k,\mathbb{C})$ qui agit par $-1$ sur le premier $\mathbb{C}^k$ et par $1$ sur le dernier $\mathbb{C}^k$. Ainsi
pour tout $g\in \GL(2k,\mathbb{C})$, $ j\cdot g= tgt$.
Evidemment $t$ est un élément de $\GSp(2k,\mathbb{C})$ de rapport de similitude $-1$. Ainsi $W_\mathbb{R}$
 laisse invariant $\Sp(2k,\mathbb{C})$  dans son action par conjugaison et $j$  agit via la conjugaison par $t$. Cela démontre le lemme. \qed

\begin{prop} Soit $\psi_\pi$ le  paramètre d'Arthur d'une série discrète de $\scrX$ (on est dans le cas 3, $n$ impair). La restriction de $\psi$
 à $W_\mathbb{R}$ se factorise, à conjugaison près,  par le produit semi-direct $\Sp(2k,\mathbb{C})\rtimes W_\mathbb{R}$,  et  donc en particulier 
 satisfait la conjecture de Sakellaridis et Venkatesh avec  ${}^LG_\scrX=\Sp(2k,\mathbb{C})\rtimes W_\mathbb{R}$.
\end{prop}
\dem
Soit $\psi_\pi$ un tel   paramètre. On a vu ci-dessus que la restriction de $\psi_\pi$ à $W_\mathbb{R}$   est orthogonale de dimension $2k$,  
plus précisément est la somme de $k$ représentations orthogonales  $\delta\left( \frac{m_i}{2}, -\frac{m_i}{2} \right)$ de dimension deux de $W_\mathbb{R}$ 
toutes distinctes (les $m_i$ sont pairs, ont les suppose rangés dans l'ordre décroissant).
 A l'aide de ces paramètres, on  construit un morphisme $\phi'$ de $W_\mathbb{R}$ dans $\GL(2k,\mathbb{C})\rtimes W_\mathbb{R}$, en posant
pour tout $ z\in \mathbb{C}^\times$
\[  \phi'(z)=\left(\left((z/\overline{z})^{m_1},\cdots, (z/\overline{z})^{m_1},\underbrace{1,\cdots,1}_{n-2k},(z/\overline{z})^{-m_k},\cdots,(z/\overline{z})^{-m_1}\right),z \right),\]
et $ \phi'(j)=(t J_{2k},j)$.
 Pour vérifier que l'on a bien construit un morphisme, il faut remarquer que les $m_i$ étant des entiers pairs $\phi'(-1)=1$ et que $(tw_{2k})^2=1$ ce qui donne 
 $\phi'(j)^2=\phi'(j^2)=\phi'(-1)$.

 Il est clair que $\phi'(z)\in \Sp(2k,\mathbb{C})\times W_\mathbb{C}$ et que $\phi'(j)\in (\GSp(2k,\mathbb{C}),j)$.
  Le sous-groupe de $\GL(2k,\mathbb{C})\rtimes W_\mathbb{R}$ engendré par $\Sp(2k,\mathbb{C})$ et $\phi'(j)$ est indépendant des paramètres
   $m_i$ et est isomorphe au produit semi-direct $\Sp(2k,\mathbb{C})\rtimes W_\mathbb{R}$ par 
   l'application qui est l'identité sur $\Sp(2k,\mathbb{C})$ et envoie $(t,j)\in (\GSp(2k,\mathbb{C}),j)$ sur $(1,j)\in \Sp(2k,\mathbb{C})\rtimes W_\mathbb{R}$.
On note $\phi$ le composé de $\phi'$ avec cet isomorphisme suivi de l'inclusion de
 $\Sp(2k,\mathbb{C})\rtimes W_\mathbb{R}$ dans le commutant dans $^L\U_n$ de l'image de $\SL(2,\mathbb{C})$.
  On obtient donc un morphisme $\psi$ de $W_\mathbb{R}\times \SL(2,\mathbb{C})$ dans $^L\U_n$. Ce morphisme est un paramètre d'Arthur discret.
   Le paquet d'Arthur qu'il détermine est uniquement déterminé par la restriction de $\psi$ à $W_\mathbb{C}\times \SL(2,\mathbb{C})$. 
   Il est clair que cette restriction est bien conjuguée de la restriction à ce même groupe du paramètre d'Arthur  $\psi_\pi$ dont on est parti. \qed

    \begin{rmq} Il est vraisemblable que 
  les plongements plus généraux obtenus dans la proposition \ref{PSOU} et la possibilité d'avoir le facteur 
$  \chi_{m_0} \boxtimes R[n-2(r+s)] $  avec $m_0\neq 0$ dans le paramètre $\psi_\pi$ servent à étendre les conjectures au cas 
où le caractère $\chi$ de $H/H_e$ n'est pas trivial.
  \end{rmq}

  {\bf Cas 4}  : les paramètres  $\psi_\pi$ sont    de la forme    (\ref{psicas4}), c'est-dire 
 $ \psi_\pi=\bigoplus_{i=1}^{n}  ((\chi_{m_i}\oplus \chi_{-m_i}) \boxtimes R[1]$,
où les $m_i$ sont des entiers impairs.
La représentation  $\chi_{m_i}\oplus \chi_{-m_i}$ de $\bbC^\times $ s'étend de manière unique en la représentation
  $\delta\left( \frac{m_i}{2},-\frac{m_i}{2} \right)$ de $W_\bbR$ cette représentation étant à valeurs dans $\SL(2,\bbC)$ si $m_i$ est impair.
Posons     $\phi_d = \bigoplus_{i=1}^{n} \delta\left( \frac{m_i}{2},-\frac{m_i}{2} \right) :\, W_\bbR\longrightarrow 
  \Sp(2n,\bbC)\times W_\bbR$. 
On pose  ${}^LG_\scrX=\Sp(2n,\bbC)\times W_\bbR$
comme prévu par \cite{KS}, et le morphisme $\bar \varphi_\scrX$ de (\ref{factpar2}) est  le morphisme $\xi_\Sp$ de 
 la proposition \ref{PSOU}, avec $a=b=0$ (avec  $2n$ à la place de $n$). On a 
la propriété de factorisation  $\psi_\pi=\varphi_\scrX \circ \phi_d$. \qed

\section{Les paquets  d'Arthur attachés aux  séries discrètes des espaces symétriques}\label{resComp}

Dans ce qui précède, nous avons attaché à chaque série discrète $\pi$  d'un des espaces symétriques $\scrX$
que nous considérons un paramètre d'Arthur $\psi_\pi$ tel que $\pi$ est dans le paquet d'Arthur $\Pi(G,\psi_\pi)$
et nous avons vérifié que ce paramètre se factorise bien de la façon prédite par les conjectures de Sakellaridis et Venkatesh.

Dans cette section, nous prenons une perspective légèrement différente : étant donné un tel paquet 
$\Pi(G,\psi)$, quels sont les éléments de ce paquet qui sont des séries discrètes d'un des  espaces symétriques $\scrX$ ?
Un résultat intéressant est que si l'on fixe un   $G$ et un tel   paramètre $\psi$, 
un élément du paquet ne contribue qu'à un seul espace symétrique parmi ceux possibles.

\subsection*{Les caractères  $\epsilon(\pi)$  de  $A(\psi)$}\label{caracAp}
Soit 
$\psi: W_\bbR\times \SL(2,\bbC)\rightarrow {}^LG$
 un paramètre d'Arthur. Posons $S_\psi=\mathrm{Centr}_{G^\vee}(\psi)$ et 
$A(\psi) =S_\psi/(S_\psi)^0$.   Lorsque les $A(\psi)$ sont abéliens, ce qui est le cas des groupes considérés ici, la théorie d'Arthur
attache à toute représentation $\pi$ du paquet d'Arthur $\Pi(\psi)$ une représentation de dimension finie $\epsilon(\pi)$ de $A(\psi)$.
Pour les groupes unitaires, spéciaux orthogonaux ou symplectiques (et  pour $\GL(n,\bbR)$ où les $A(\psi)$ sont triviaux)
on sait d'après \cite{MR7}, \cite{MR3} que les $\epsilon(\pi)$ sont de dimension $1$, c'est-à-dire, pour tout $\pi \in \Pi(\psi)$, 
$\epsilon(\pi)\in \widehat{A(\psi)}$. Ces   $\epsilon(\pi)$  ont de plus été déterminés de manière explicite, dépendant de certains choix
(celui d'une forme intérieure pure quasi-déployée  de référence, et d'une donnée de Whittaker pour celle-ci).
Soit $\psi$ un paramètre de bonne parité pour un groupe $G$ comme dans la section \ref{secparartbp}, donc 
comme en (\ref{parunit}) ou   (\ref{parclass}). Alors 
$A(\psi)\simeq   (\bbZ/2\bbZ)^R$.
Un caractère de ce groupe sera donc donné par une suite $(\epsilon_1, \ldots , \epsilon_R)$ 
à valeurs dans $\bbZ/2\bbZ$. 
On calcule ci-dessous ces invariants $\epsilon(\pi)$ pour les séries discrètes  $\pi$ des espaces symétriques   et le paramètre 
d'Arthur qui leur est attaché  et l'on en tire  quelques conséquences.

\subsection*{Cas 1 : $\scrX=\GL(n,\bbR)/\GL(p,\bbR)\times \GL(n-p,\bbR)$,   $2p\leq n$}
Les paquets d'Arthur sont des singletons, et les groupes $A(\psi)$ sont triviaux. Un paramètre d'Arthur  $\psi$ détermine
donc une unique représentation  $\pi(\psi)$.   Les paramètres considérés sont de la forme
 \[  \psi =
\bigoplus_{i=1}^p \left(\delta\left(  \frac{m_i}{2},   \frac{m_i}{2} \right) \boxtimes R[1] \right)\oplus (\Triv_{W_\bbR} \boxtimes R[n-2p]) \]
les $m_i$ étant des entiers impairs distincts.
La  représentation unitaire irréductible $\pi(\psi)$ de $\GL(n,\bbR)$  attachée à un tel paramètre ne contribue donc qu'à un seul des espaces symétriques
ci-dessus (lorsque $p$ varie).

\subsection*{Cas 2 :  $\scrX=\U(p,q)/ \Or(p,q)$, $n=p+q$}\label{etapi2} 
Nous avons vu que les séries discrètes de $\scrX$ sont des séries discrètes de $G$. 
Ce qui les distingue  parmi celles-ci est que leur $K$-type minimal a des invariants non triviaux sous $K\cap H=\Or(p)\times \Or(q)$.
Un paquet de Langlands de séries discrètes de $\U(p,q)$ est déterminé par un caractère infinitésimal entier régulier, 
c'est-à-dire dans un système de coordonnées usuelles, un $n$-uplet $\lambda=(\lambda_1,\ldots, \lambda_n)$
où les $\lambda_i$ forment une suite strictement décroissante constitués d'entiers si $n$ est impair, et 
de demi-entiers non entiers si $n$ est pair. Notons $\Pi(\lambda)$ ce paquet de Langlands.
 Remarquons que  $\rho=\left( \frac{n-1}{2}, \frac{n-3}{2}, \ldots , -\frac{n-1}{2} \right)$ de sorte 
que $\lambda-\rho$ est toujours un n-uplet d'entiers.

Les représentations à l'intérieur d'un paquet de Langlands $\Pi(\lambda)$ de séries discrètes de paramètre de Langlands $\psi$
sont paramétrées par le caractère de $A(\psi)$ qui leur est  attaché. Pour la normalisation de cette paramétrisation,  
 on  suppose que la forme intérieure pure de 
$\U(n)$ qui contient la représentation correspondant au  caractère trivial de  $A(\psi)$  est $\U(\frac{n+1}{2},\frac{n-1}{2})$ si n est impair et 
 $\U(\frac{n}{2},\frac{n}{2})$ si $n$ est pair.

\begin{prop}
Si $n=p+q$ est pair, et  si  $\lambda-\rho$ est un n-uplet d'entiers tous de même parité  $\omega$, les séries discrètes 
de $\U(p,q)$ dans le paquet $\Pi(\lambda)$  sont toujours des séries discrètes pour $\U(p,q)/\SO(p,q)$, et elles le sont pour 
$\U(p,q)/\Or(p,q)$ si $p$  a la parité $\omega$.

Si $n$ est impair, alors  quel que soit $\lambda$, il existe exactement une forme intérieure 
$\U(p,q)$ (à permutation près de $p$ et $q$, et $p+q=n$) et une série discrète $\pi$  de caractère infinitésimal $\lambda$ 
pour cette forme intérieure qui est une série discrète pour 
$\U(p,q)/\SO(p,q)$. Cette série discrète $\pi$ est celle qui correspond au caractère  $\epsilon(\pi)$ de $A(\psi)$ donné par 
 $\epsilon(\pi)_i=(-1)^{\lambda_i+\frac{n-1}{2}}$ avec 
$p$ le nombre d'entiers $i$ tel que $\lambda_i+\frac{n+1}{2}-i$ est pair.
\end{prop}
Remarquons aussi que la série discrète $\pi^-$ de $\U(q,p)$ de caractère infinitésimal $\lambda$
apparaissant dans  l'espace symétrique $\U(q,p)/\SO(q,p)$
 (qui s'identifie à $\pi$ si l'on identifie $\U(p,q)$ et $\U(q,p)$) correspond pour  ce choix de paramétrisation, 
 au caractère opposé du $A(\psi)$, c'est-à-dire
  $\epsilon(\pi^-)_i=(-1)^{\lambda_i+\frac{n+1}{2}}$.

\medskip 

\dem 
Pour démontrer cette proposition, on explicite  le $K$-type minimal d'une  série discrète $\pi$ d'un paquet
$\Pi(\lambda)$ en fonction de $\epsilon(\pi)$. Ce calcul est un peu technique. Pour que $\pi$ soit 
une série discrète de $\U(p,q)/\SO(p,q)$, le  plus haut poids de ce $K$-type  
 doit être formé d'entiers ayant tous même parité  et  pour être  une série discrète de $\U(p,q)/\Or(p,q)$ ces entiers doivent  être tous pairs.
  On fixe $\pi$ dans $\Pi(\lambda)$ avec $\epsilon(\pi)=(\epsilon_1,\cdots, \epsilon_n)$. On pose 
$$p:=|\{i; \epsilon_i=(-1)^{i-1}\}|.$$
On note $L$ le tore compact dont le $i$-ème  facteur $\U(1)$ est dans $\U(p)$ si et seulement si 
$\epsilon_i=(-1)^{i-1}$. Avec cela on détermine la sous-algèbre parabolique $\theta$-stable $\mathfrak{q}=\frl\oplus \fru$  de $\frg\frl(n,\mathbb{C})$ 
permettant de réaliser $\pi$ comme $A_\frqqq(\lambda-\rho)$.
Le plus haut poids du $\U(p)\times \U(q)$-type minimal est alors  $\lambda-\rho+2\rho(\mathfrak{u}\cap \mathfrak{p})$.

 Pour tout $i\in [1,n]$, on pose 
$\tilde{\mu}_i:=\lambda_i+(n+1)/2-i=(\lambda-\rho)_i$ et 
\begin{align*}
\mu_i:=\tilde{\mu}_i&+|\{j\in [1,n]; \epsilon_i(-1)^{i-1}\neq \epsilon_j(-1)^{j-1}\hbox{ et }\lambda_j <\lambda_i\}|\\ &-
|\{j\in [1,n]; \epsilon_i(-1)^{i-1}\neq \epsilon_j(-1)^{j-1} \hbox{ et }\lambda_j >\lambda_i \}|.
\end{align*}
Le plus haut poids du $K$-type minimal de $\pi$  est alors 
$\left( \mu_i; \epsilon_i=(-1)^{i-1}\right); \left(\mu_i;\epsilon_i\neq (-1)^{i-1}\right)$.

Distinguons les cas :  si $n$ est pair, pour que les $\mu_i$ aient tous même parité, il faut et il suffit que 
 les $\tilde{\mu}_i$ aient tous même parité. Cela permet de conclure  facilement dans ce  cas.

Si $n$ est impair,  pour que les $\mu_i$ aient tous même parité, il faut et il suffit que 
que les $\tilde{\mu}_i$ pour $i$ tel que $\epsilon_i=(-1)^{i-1}$ aient tous même parité et soient 
de parité opposé aux $\tilde{\mu}_i$ tel que $\epsilon_i\neq (-1)^{i-1}$.
Il  faut donc nécessairement que $(p,n-p)$ soit à l'ordre près
$$
(|\{i;\tilde{\mu}_i\equiv 0[2]\}|,|\{i;\tilde{\mu}_i\equiv 1[2]\}|),
$$
et  si l'on prend $p=|\{i;\tilde{\mu}_i\equiv 0[2]\}|$,  il faut aussi que 
 $\epsilon_i=(-1)^{i-1}$ si et seulement si $(-1)^{\lambda_i+(n+1)/2-i}=1$, 
 c'est-à-dire que pour tout $i$ on doit avoir $\epsilon_i=(-1)^{\lambda_i+(n-1)/2}$. 
 Par contre si l'on prend $p=|\{i;\tilde{\mu}_i\equiv 1[2]\}|$ alors $\epsilon_i=(-1)^{i-1}$ 
 si et seulement si $(-1)^{\lambda_i+(n+1)/2-i}=-1$, c'est-à-dire que pour tout $i$ on doit avoir $\epsilon_i=(-1)^{\lambda_i+(n+1)/2}$. \qed

\subsection*{Cas 3 :  $\scrX=\U(p,q)/ \U(r,s)\times \U(r',s')$}\label{etapi3} Notons  cet espace symétrique de manière plus précise
par $\scrX_{r,s}$, ce qui nous permet de faire varier ces deux entiers. 
On a vu en (\ref{psicas3}) que les paramètres $\psi$ sont de la forme
$  \psi=\bigoplus_{i=1}^{r+s}  ((\chi_{m_i} \oplus \chi_{-m_i})\boxtimes R[1])   \oplus \left(\chi_0 \boxtimes R[n-2(r+s)] \right)$. 
Le résultat suivant découle des calculs des $\epsilon(\pi)\in \widehat{A(\psi)}$ dans \cite{MR7} pour $\pi \in \Pi(\psi)$.  
On écrit  $\epsilon(\pi)=(\epsilon_1(\pi), \ldots,  \epsilon_{2(r+s)}(\pi), \epsilon_0(\pi))$,
où le $\epsilon_0$ vient du facteur  $\chi_0 \boxtimes R[n-2(r+s)]$ du paramètre et est donc absent si $n=2(r+s)$.
\begin{prop}\label{Ap63}
Soit $\psi$ comme ci-dessus et $\pi\in \Pi(\U(p,q),\psi)$. Alors $\pi$   est dans le spectre discret de $\scrX_{r,s}$ si et seulement si
$\epsilon(\pi)=(\epsilon_1(\pi), \ldots,  \epsilon_{2(r+s)}(\pi), \epsilon_0(\pi) )$
vérifie 
\[  \#\left\{ i\in \llbracket 1,2(r+s)\rrbracket \, \vert \,   \epsilon_i(\pi)=(-1)^{i-1}   \right\} =2r,  \]
 \[  \#\left\{ i\in \llbracket 1,2(r+s)\rrbracket \, \vert \,   \epsilon_i(\pi)=(-1)^{i}   \right\} =2s , \]
(rappelons que ceci dépend du choix d'une donnée forme intérieure pure quasi-déployée, un autre choix pourrait inverser les deux conditions)
et 
\[   \epsilon_{2i-1}(\pi)\epsilon_{2i} (\pi)=(-1)^{p+q-1}, \; (i=1,\ldots,r+s),\]
condition qui ne dépend pas du choix de la forme intérieure pure quasi-déployée.
\end{prop}

\dem On suppose d'abord que le caractère infinitésimal du paramètre $\psi$ est très régulier, c'est-à-dire
que les $m_i$ sont des entiers strictement positifs  assez distincts, que l'on a mis dans l'ordre décroissant.
Soit $\pi\in \Pi(\U(p,q),\psi)$. Alors $\pi$ est induit cohomologiquement à partir d'une paire $(\frqqq,L)$, où 
$L$ est le $c$-Levi fixé par l'espace symétrique, si et seulement si les conditions 
\[  \#\left\{ i\in \llbracket 1,2(r+s)\rrbracket \, \vert \,   \epsilon_i(\pi)=(-1)^{i-1}   \right\} =2r,  \quad
   \#\left\{ i\in \llbracket 1,2(r+s)\rrbracket \, \vert \,   \epsilon_i(\pi)=(-1)^{i}   \right\} =2s , \]
 sont vérifiées.
 
 Posons $t_i=\frac{m_i}{2}-\frac{n-2i+1}{2}$.  Les $m_i$ sont entiers, de la parité de $n-1$, et donc les $t_i$
 sont entiers.  Comme on a supposé les $m_i$ grands et suffisamment distincts,
les $t_i$ forment encore une suite strictement  décroissante.
 Considérons le caractère $\pi_L=\pi_L(t_1,-t_1, \ldots, t_{r+s},-t_{r+s}, 0)$ de $L=\U(1)^{2(r+s)}\times \U(r-r',s-s')$ trivial sur le facteur $\U(r-r',s-s')$
 et donné par un entier $t$ sur les facteurs $\U(1)$.
 Il existe une unique sous-algèbre parabolique $\theta$-stable $\frqqq_1=\frl\oplus \fru_1$
de $\frg$ telle que  $\pi_L(t_1,-t_1, \ldots, t_{r+s},-t_{r+s}, 0)$ soit dans le good range pour $\frqqq_1$.
 Il est clair que $\pi_L\in \caP(\frqqq_1)$.
La représentation $\pi=A_{\frqqq_1}(\pi_L)$ est alors dans le spectre discret de $\scrX$, et dans $\Pi(\U(p,q),\psi)$.
De plus $\epsilon(\pi)$ vérifie  les conditions voulues d'après \cite{MR7}.

Considérons maintenant un caractère $\pi_L=\pi_L(t_1,-t_1, \ldots t_{r+s},-t_{r+s}, 0)$ pour certains entiers $t_1,\ldots, t_{r+s}$
que l'on suppose seulement suffisamment distincts les uns des autres et de leurs opposés. 
A conjugaison près dans $K=\U(p,0)\times \U(0,q)$ on peut supposer les $t_i$ positifs
avec  $t_1>t_2\ldots> t_r$ et  $t_{r+1}>t_{r+2}\ldots> t_{r+s}$.
Ceci détermine une sous
algèbre parabolique $\theta$-stable $\frqqq=\frl\oplus \fru$
de $\frg$ telle que  $\pi_L(t_1,-t_1, \ldots t_{r+s},-t_{r+s}, 0)$ soit dans le good range pour $\frqqq_1$.
Si on pose  $t'_i=\frac{m_i}{2}-\frac{n-2i+1}{2}$ où la suite $(t'_i)_{i=1,\ldots r+s}$
est obtenue à partir des $t_i$ en les mettant dans l'ordre décroissant, alors $\pi= A_{\frqqq}(\pi_L)$
est dans le paquet $\Pi(\U(p,q),\psi)$ et elle contribue au spectre discret de $\scrX$ 
(on fait ici agir un élément de $N_K(\frt_0)\backslash N_{\mathbf{G}}(\frt_0)\simeq \frS_r\times \frS_s\backslash \frS_n$ pour ce ramener au calcul précédent)
et de même que ci-dessus, $\epsilon(\pi)$ vérifie la proposition voulue.
On peut conclure par un argument de comptage pour montrer que les conditions de la proposition sont suffisantes : le nombre d'éléments
du paquet $\Pi(\U(p,q),\psi)$ vérifiant ces conditions est le cardinal de 
$N_K(\frt_0)\backslash N_{\mathbf{G}}(\frt_0)\simeq \frS_r\times \frS_s\backslash \frS_n$.\qed

\begin{rmq}\label{mult163}
Les résultats de multiplicité un \cite{MR7} 
impliquent donc ici que le spectre discret de $\scrX_{r,s}$ à la propriété de multiplicité un.   
\end{rmq}

\begin{rmq}\label{ocrs}
Un élément du paquet d'Arthur $\Pi(\U(p,q),\psi)$ ne contribue au spectre discret  que d'au plus l'un  des $\scrX_{r,s}$ lorsque l'on fait varier
$r$ et $s$ .
\end{rmq}

\subsection*{Cas 4  :  $\scrX= \U(n,n)/\GL(n,\bbC)$}
Les paramètres s'écrivent 
$ \psi=\bigoplus_{i=1}^{n}  ((\chi_{m_i}\oplus \chi_{-m_i}) \boxtimes R[1])$.   

\begin{prop}\label{Ap64}Soit $\psi$ comme ci-dessus et $\pi\in \Pi(\U(n,n),\psi)$. Alors  $\pi$ est dans le spectre discret de $\scrX$ si et seulement si
$\epsilon(\pi)=(\epsilon_1(\pi), \ldots,  \epsilon_{2n}(\pi))$
vérifie 
\[ \#\left\{ i\in \llbracket 1,2n\rrbracket \, \vert \,   \epsilon_i(\pi)=(-1)^{i-1}   \right\} =n  
\text{ et } 
   \epsilon_{2i-1}(\pi)\epsilon_{2i} (\pi)=1, \, (i=1,\ldots, n).\]
\end{prop}

La démonstration est similaire à celle de la proposition \ref{Ap63}. Il en est de même des propositions analogues des cas qui suivent.

\begin{rmq}
Pour  $\scrX'=\U(2r,2s)/\U(r,s)\times \U(r,s)$, les paramètres 
d'Arthur sont les mêmes que pour $\scrX=\U(n,n)/\GL(n,\bbC)$ avec  $n=r+s$, 
 mais une  représentation $\pi$ d'un tel paquet  ne contribue au plus qu'à un seul des deux spectres : 
celui de $\scrX$ si $\epsilon(\pi)$ est comme dans la proposition \ref{Ap63},  
et celui de  $\scrX'$  si $\epsilon(\pi)$ est comme dans la proposition \ref{Ap64}, les deux cas étant mutuellement exclusifs.
 
 On a aussi, comme dans la remarque \ref{mult163}, la propriété de multiplicité un du  spectre discret.
 \end{rmq}

\subsection*{Cas 5  :  $\scrX= \U(2p,2q)/\Sp(p,q)$}   \label{etapi5}
Les paramètres sont de la forme
$  \psi=\bigoplus_{i=1}^{p+q}  (\chi_{m_i} \boxtimes R[2])$.   
\begin{prop}\label{Ap65}Soit $\psi$ comme ci-dessus et $\pi\in \Pi(\U(2p,2q),\psi)$. Alors 
$\pi$ est dans le spectre discret de $\scrX$ si et seulement si $    \epsilon_i(\pi)=-1, \, (i=1,\ldots, p+q)$.
\end{prop}

\subsection*{Cas 6  :  $\scrX= \U(n,n)/\Sp(2n,\bbR)$}
Les paramètres sont de la forme
$  \psi=\bigoplus_{i=1}^{n}  (\chi_{m_i} \boxtimes R[2])$.

\begin{prop}\label{Ap66} Soit $\psi$ comme ci-dessus et $\pi\in \Pi(\U(n,n),\psi)$.
 Alors 
$\pi$ est dans le spectre discret de $\scrX$ si et seulement si
$   \epsilon(\pi)$ est le caractère trivial.
\end{prop}

\begin{rmq}
Pour $\scrX'=\U(2n,2n)/\Sp(n,n)$, les paramètres 
d'Arthur sont les mêmes que ci-dessus, mais une  représentation d'un tel paquet  ne contribue au plus qu'à un seul des deux spectres, 
celui de $\scrX$ si  $  \epsilon_i(\pi)=-1$, \, (i=1,\ldots, n) , et celui de $\scrX'$  si   $\epsilon(\pi)$ est trivial.
\end{rmq}

 \subsection*{Cas 7 et 8  : $\SO(p,q)/\SO(r,s)\times \SO(r',s') $}
Les paramètres sont de la forme
\[\psi=\left( \bigoplus_{i=1}^{r+s}  \delta\left(\frac{m_i}{2}, -\frac{m_i}{2}\right) \boxtimes R[1]\right) \oplus (\Triv_{W_\bbR}\boxtimes R[2(n-r-s)]),  \]
si $p+q=2n+1$, et 
 \[\psi=\left( \bigoplus_{i=1}^{r+s}  \delta\left(\frac{m_i}{2}, -\frac{m_i}{2}\right) \boxtimes R[1]\right) \oplus (\Triv_{W_\bbR}\boxtimes R[2(n-r-s)-1]) 
\oplus (\sgn_{W_\bbR}^{p-n} \boxtimes R[1]) ,\]
si $p+q=2n$.

 \begin{prop}\label{Ap67}
Soit $\psi$ comme ci-dessus et $\pi\in \Pi(\SO(p,q),\psi)$. Alors $\pi$   est dans le spectre discret de $\scrX$ si et seulement si
$\epsilon(\pi)=(\epsilon_1(\pi), \ldots,  \epsilon_{(r+s)}(\pi) )$
vérifie 
\[  \#\left\{ i\in \llbracket 1,(r+s)\rrbracket \, \vert \,   \epsilon_i(\pi)=(-1)^{i-1}   \right\} =r,  \quad 
   \#\left\{ i\in \llbracket 1,(r+s)\rrbracket \, \vert \,   \epsilon_i(\pi)=(-1)^{i}   \right\} =s  \]
(rappelons que ceci dépend du choix d'une donnée de Whittaker ou d'une forme intérieure pure quasi-déployée, un autre choix pourrait inverser les deux conditions).

\end{prop}

\begin{rmq} En notant $\scrX_{r,s}$ l'espace symétrique ci-dessus, on peut faire varier $r$ et $s$, et on a alors le même résultat
que dans la remarque \ref{ocrs}.
\end{rmq}

 \subsection*{Cas 9 : $\SO(2p,2q)/\U(p,q)$} \label{epi9}
Les paramètres sont de la forme
$\psi=\bigoplus_{i=1}^{n/2  }   \delta\left(\frac{m_i}{2}, - \frac{m_i}{2}\right)\boxtimes R[2]$, 
si $n=p+q$ est pair et 
$ \psi=\left(\bigoplus_{i=1}^{\frac{n-1}{2}  }   \delta\left(\frac{m_i}{2}, - \frac{m_i}{2}\right)\boxtimes R[2]  \right) 
 \oplus \left( \Triv_{W_\bbR} \boxtimes R[1]  \right)\oplus \left( \sgn_{W_\bbR} \boxtimes R[1]  \right)$,  
si $n=p+q$ est impair. 

\begin{prop}\label{Ap69}
Soit $\psi$ comme ci-dessus et $\pi\in \Pi(\SO(2p,2q),\psi)$. Alors $\pi$   est dans le spectre discret de $\scrX$ si et seulement si
$\epsilon_i(\pi)=-1$, pour tout $i=1,\ldots n$.
\end{prop}

 \subsection*{Cas 10 : $\SO(n,n)/\GL(n, \bbR)$} \label{epi10}
Les paramètres sont de la forme
$ \psi= \bigoplus_{i=1}^{\frac{n}{2}  }   \delta\left(\frac{m_i}{2}, - \frac{m_i}{2} \right) \boxtimes R[2]$,   
si $n$ est pair et 
 $ \psi=\left(\bigoplus_{i=1}^{\frac{n-1}{2}  }   \delta\left(\frac{m_i}{2}, - \frac{m_i}{2}\right)\boxtimes R[2]  \right) 
 \oplus \left( \Triv_{W_\bbR} \boxtimes R[1] \right)\oplus \left( \Triv_{W_\bbR} \boxtimes R[1]
 \right)$,   
si $n$ est impair.
\begin{prop}\label{Ap610}
Soit $\psi$ comme ci-dessus et $\pi\in \Pi(\SO(n,n),\psi)$. Alors $\pi$   est dans le spectre discret de $\scrX$ si et seulement si
$\epsilon(\pi)$ est trivial.
\end{prop}

\begin{rmq} 
 Pour  $\scrX'=\SO(2n,2n)/\U(n,n)$,  les paramètres 
d'Arthur sont les mêmes que ci-dessus, mais une  représentation d'un tel paquet  ne contribue au plus qu'à un seul des deux spectres, 
celui de $\frX$ si  $\epsilon_i(\pi)=-1$, pour tout $i=1,\ldots n$, et celui de $\scrX'$  
si $\epsilon(\pi)$ est trivial.
\end{rmq}

 \subsection*{Cas 11 : $\Sp(2n,\bbR)/\Sp(2p,\bbR) \times \Sp(2(n-p),\bbR) $}   \label{epi11}
Les paramètres sont de la forme
\[\psi=\left( \bigoplus_{i=1}^p   \delta\left(\frac{m_i}{2}, -\frac{m_i}{2}\right) \boxtimes R[2]\right) \oplus (\Triv_{W_\bbR}\boxtimes R[2(n-2p))+1]) \]
avec les $m_i$ impairs.
\begin{prop}\label{Ap611}
Soit $\psi$ comme ci-dessus et $\pi\in \Pi(\Sp(2n,\bbR),\psi)$. Alors $\pi$   est dans le spectre discret de $\scrX$ si et seulement si
$\epsilon(\pi)$ est trivial.
\end{prop}

\subsection*{Cas 12 : $\Sp(4n,\bbR)/\Sp(2n,\bbC) $}   \label{epi12}
Les paramètres sont de la forme
\[\psi=\left( \bigoplus_{i=1}^n   \delta\left(\frac{m_i}{2}, -\frac{m_i}{2}\right) \boxtimes R[2]\right) \oplus (\Triv_{W_\bbR}\boxtimes R[1]) \]
avec les $m_i$ impairs.
\begin{prop}\label{Ap612}
Soit $\psi$ comme ci-dessus et $\pi\in \Pi(\Sp(2n,\bbR),\psi)$. Alors $\pi$   est dans le spectre discret de $\scrX$ si et seulement si
$\epsilon_i(\pi)=-1$, pour tout $i=1,\ldots, n$.
\end{prop}

\begin{rmqs} 
 Remarquons que le fair range est ici dans le weakly good range. On en déduit que le spectre discret de $\scrX$ à la propriété de multiplicité un.

--- Pour $\scrX'=\Sp(4n,\bbR)/\Sp(2n,\bbR)\times \Sp(2n,\bbR)$,  les paramètres 
d'Arthur sont les mêmes que ci-dessus, mais une  représentation  $\pi$ d'un tel paquet  ne contribue au plus qu'à un seul des deux spectres, 
celui de $\scrX$ si  $\epsilon_1(\pi)=-1$, pour tout $i=1,\ldots n$., et celui de $\scrX'$  
si $\epsilon(\pi)$ est trivial.
\end{rmqs}

\subsection*{Cas 13 : $\Sp(2n,\bbR)/\GL(n,\bbR) $}   \label{epi32}
Comme dans le cas 2, les séries discrètes de $\scrX$ sont les séries discrètes de $G$ distinguées par le fait que leur $K$-type minimal
a des invariants non triviaux sous $K\cap H=\Or(n)$.
Un paquet de Langlands de séries discrètes de $\Sp(2n,\bbR)$ est déterminé par un caractère infinitésimal entier régulier, 
c'est-à-dire dans un système de coordonnées usuel, un $n$-uplet $\lambda=(\lambda_1,\ldots, \lambda_{n})$
où les $\lambda_i$ forment une suite strictement décroissante constitués d'entiers. Notons $\Pi(\lambda)$ ce paquet de Langlands.

\begin{prop}
Soit $\Pi(\lambda)$ comme ci-dessus un paquet de Langlands de séries discrètes de $\Sp(2n,\bbR)$.
Alors soit tous les éléments de ce paquet sont des séries discrètes de $\scrX$, soit aucun ne l'est, et la première éventualité
se produit exactement quand les entiers $\lambda_i$ ont une parité qui alterne en commençant par un nombre impair, 
c'est-à-dire pour tout $i\in [1,n]$, $\lambda_i\equiv i[2]$.
\end{prop}

\begin{rmq}
Si les $\lambda_i$ ont bien une parité qui alterne mais en commençant par un nombre pair,
 on trouverait des séries discrètes pour $L^2_d(\Sp(2n,\mathbb{R})/\GL(n,\mathbb{R}))_\chi$ où $\chi$ 
 est le caractère ${\sgn\circ \det}$ de $\GL(n,\mathbb{R})$.
\end{rmq}

\dem
Remarquons que  $\rho=\left( n, n-1, \ldots , 1 \right)$ de sorte 
que $\lambda-\rho$ est toujours un n-uplet d'entiers.
Notons  $\pi_L(t_1,\ldots ,t_n)$ le caractère de $L=\U(1)^n$ donné par le $n$-uplet
d'entiers $(t_1,\ldots,t_n)$. La série discrète $A_\frqqq(\pi_L(t_1,\ldots,t_n))$, où $\frqqq=\frl\oplus \fru$ est la sous-algèbre de Borel
 telle que $\pi_L(t_1,\ldots,t_n)$ soit  dans le good range pour $\frqqq$ admet un 
$K$-type minimal ayant   des invariants non triviaux sous $K\cap H=\Or(n)$ si et seulement si 
$\pi_L\otimes \bigwedge^{\mathrm{top} } (\fru\cap \frp)$ a une restriction triviale à $L\cap K\cap H=\{\pm 1\}^n$.
Or on calcule facilement le caractère  $\bigwedge^{\mathrm{top}}(\fru\cap \frp)$ de $L\cap H\cap K$ : dans le système de coordonnées usuel
où les racines  (toutes imaginaires) sont les $\pm e_i \pm e_j$,  $1\leq i<j\leq n$, $\pm 2e_i$, $1\leq i\leq n$, les racines compactes
étant les $\pm (e_i- e_j)$. Les racines non compactes longues positives ($\pm 2e_i$, selon le signe de $t_i$)
ont une restriction triviale à à $L\cap K\cap H$ et la contribution des racines 
 non compactes courtes positives ($\pm (e_i+e_j)$, selon le signe de $t_i+t_j$) donne le caractère $\sgn^{n-1}$ sur chaque facteur $\{\pm 1\}$ de 
 $L\cap K\cap H$.  Ainsi la condition cherchée  est que 
  les $t_i$ sont tous pairs si $n$ est impair, et tous impairs si $n$ est pair.
 On en déduit le résultat. \qed

 
\section{Le sous-groupe $L$}\label{secL}

Dans cette section, nous déterminons  le sous-groupe $L=\mathrm{Cent}_G(\frt_0)$ ({\sl cf.}  section \ref{SDX}) par un argument de 
réduction au rang un, et nous commençons par là.
On utilise aussi le fait que l'on connait la forme générale des $c$-Levi des groupes classiques :

$G=\GL(n,\bbR)$, $L\simeq \left(\prod_{i=1}^m \GL(a_i,\bbC)\right) \times \left(\prod_{j=1}^s \GL(a'_j,\bbR)\right)$, avec 
$2\sum_{i=1}^m  a_i+\sum_{j=1}^s  a'_j=n$,

 $G=\U(p,q)$,  $L\simeq \prod_{i=1}^s \U(p_i,q_i)$, avec $\sum_{i=1}^s (p_i,q_i)=(p,q)$,

 $G=\Sp(2n,\bbR)$, $L\simeq \left( \prod_{i=1}^s \U(p_i,q_i)\right) \times \Sp(2a,\bbR) $, 
avec $2\left(\sum_{i=1}^s p_i+q_i\right)+a=n$,
 
 $G=\SO(p,q)$, $L\simeq \left( \prod_{i=1}^s \U(p_i,q_i)\right) \times \SO(r,s) $, 
avec $\left( \sum_{i=1}^s (2p_i, 2q_i )\right)+(r,s)=(p,q)$.

\subsection{Le rang un}\label{rangun} Nous calculons $L$ au cas par cas.

 \subsection*{Cas 1 :   $\scrX=\GL(2,\bbR)/\GL(1,\bbR)\times \GL(1,\bbR)$ }\label{Rg1GL}
 Le sous-groupe $H$ est  le sous-groupe des matrices diagonales de $\GL(2,\bbR)$.
 Un sous-espace de Cartan compact pour $\scrX$ est l'espace $\frt_0$ des matrices antisymétriques  
 et son centralisateur $L$ dans $G$ est 
 $ L=\left\{     \begin{pmatrix}  a&b\\ -b&a \end{pmatrix}, \; a,b\in \bbR, \, a^2+b^2\neq 0  \right\} $. 
 Ce groupe $L$ est isomorphe à $\bbC^\times$, et $L\cap H= \left\{     \begin{pmatrix}  a&0\\ 0&a \end{pmatrix}, \; a\in \bbR^\times
   \right\} \simeq \bbR^\times$.

\subsection*{Cas 2  :  $\scrX=\U(1)/\Or(1)$}  
Ce cas est trivial. On a $L=\U(1)$ et $L\cap H=H=\Or(1)$.

\subsection*{Cas 3 : $\scrX=G/H=\U(2)/ \U(1)\times \U(1)$}\label{Upq1}
On réalise  $G$ comme l'ensemble des matrices carrées de taille 2 à coefficients complexes $M$ telles que ${}^t\bar MM=I_2$,  et $H$
est le sous-groupe diagonal. L'involution $\sigma$ est la conjugaison par $\begin{pmatrix}1&0\\0&-1\end{pmatrix}$, et le sous-espace propre
$\frs$ de $\sigma$ dans $\frg$ pour la valeur propre $-1$ est l'ensemble des matrices antihermitiennes. Un sous-espace de Cartan est donc donné
par $\frt_0=\left\{ \begin{pmatrix}0&\lambda\\-\lambda&0\end{pmatrix}, \, \lambda\in \bbR\right\}$.
On a respectivement  
\[L=\left\{ \begin{pmatrix}z\cos \theta & z\sin \theta \\ -z\sin \theta & z \cos\theta \end{pmatrix}, \,\begin{matrix} \theta \in \bbR \\ z\in \U(1) \end{matrix}
\right\} \simeq \U(1)\times \U(1),  \, 
L\cap H=\left\{ \begin{pmatrix}z & 0 \\ 0 &  z\end{pmatrix}, \, z\in \U(1) \right\}\simeq \U(1).\]

 \subsection*{Cas 4 :    $\scrX=\U(1,1)/\GL(1,\bbC)$}
    On réalise $\U(1,1)$ comme le sous-groupe de $\GL(2,\bbC)$ préservant la forme hermitienne
    $q:\, (x,y)\mapsto  \bar x y+x\bar y$. Le sous groupe $H=\GL(1,\bbC)$ est alors réalisé comme l'ensemble des matrices 
    $\begin{pmatrix}  \lambda&0\\0&\bar \lambda^{-1} \end{pmatrix}$, $\lambda\in \bbC^\times$.
    L'involution $\sigma$ est la conjugaison par la matrice  $\begin{pmatrix}  1&0\\0&-1 \end{pmatrix}$.
    Le sous-espace propre de $\sigma$ pour la valeur propre $-1$ dans $\frg_0$ et un sous  sous-espace de Cartan sont donnés
    respectivement par  
    \[ \frs_0= \left\{   \begin{pmatrix}  0&b \\ c&0 \end{pmatrix}  , \; b,c \in \bbR \right\}, \quad \frt_0= \left\{   \begin{pmatrix}  0&b \\ b&0 \end{pmatrix}  , \; b \in \bbR  \right\}. \]
    Le centralisateur $L$ dans $G$ de $\frt_0$ est 
    \[  L= \left\{   \begin{pmatrix}  a&b \\ b& a \end{pmatrix}  \, \; a \bar a+ b \bar b=1,  \, a\bar b+\bar a b=0  \right\}=
    \,  \left\{ \begin{pmatrix}  \frac{z_1+z_2}{2}&   \frac{z_1-z_2}{2}\\  \frac{z_1-z_2}{2}&  \frac{z_1+z_2}{2} \end{pmatrix} , 
    z_1,z_2\in \U(1)  \right\}  \simeq \U(1)\times \U(1). \]
    et 
    $ L\cap H=\left\{ \begin{pmatrix}  z&0\\ 0&z \end{pmatrix}, \; z\in \U(1) \right \}\simeq \U(1)$.

\subsection*{Cas 5 :  $\scrX=G/H=\U(2,0)/\Sp(1,0)=\U(2)/\SU(2)$}
Ici le sous-espace de Cartan $\frt_0$ est central et donc $L=G=\U(2)$.

 \subsection*{Cas 6 : $\scrX=G/H=\U(1,1)/\Sp(2,\bbR)=\U(1,1)/\SU(1,1)$}
Ici le sous-espace de Cartan $\frt_0$ est central, et donc $L=G=\U(1,1)$.

 \subsection*{Cas 7 et 8 :  $\scrX=G/H=\SO(2)/\SO(1)\times \SO(1)$ }
Ce cas est trivial

\subsection*{Cas 9: $\SO(4)/\U(2)$} 
On réalise $\SO(4)$ de la façon suivante : soit 
$  \caM= \left\{ \begin{pmatrix}  a&b\\- \bar b &\bar a\end{pmatrix},\,  a,b\in \bbC\right\}$, 
que l'on peut voir soit comme un espace vectoriel réel de dimension $4$, muni de la forme quadratique
euclidienne donnée par le déterminant, soit comme un espace complexe de dimension 2, muni de la forme hermitienne, 
elle aussi donnée par le déterminant. On note $\caM ^\times$ les éléments inversibles de $\caM$, et l'on pose
$  G=   \left\{   (h,h')\in \caM^\times \times \caM^\times, \;  \det (h)=\det(h')  \right\}/\sim $,
la relation d'équivalence  $\sim$ étant définie par 
$(h,h')\sim (\lambda h, \lambda h')$, $h,h'\in \caM^\times$, $\lambda \in \bbR^\times$.
L'action de $(h,h')\in \caM^\times \times \caM^\times$ sur $m\in \caM$ est donnée par :
$   (h,h')\cdot m= hm(h')^{-1}$. 
Elle passe au quotient et définit une action de $G$ sur $\caM$ qui identifie $G$ et $\SO(\caM)=\SO(4)$.

Soit $H$ le sous-groupe de $G$ des éléments   qui commutent à l'action de $\bbC^\times$
sur $\caM$ donnée par la multiplication à gauche par la matrice 
$  \begin{pmatrix}  \lambda &0\\ 0 &\bar \lambda \end{pmatrix}$, 
pour tout $\lambda\in \bbC$. 
Ainsi
\[  H= \left\{   \left(  \begin{pmatrix}  a &0\\ 0 &\bar a \end{pmatrix} , h'\right) , \, a\in \bbC^\times, \, h'\in \caM^\times , \, \det (h')=a\bar a\right \}/\sim .  \]
Son action sur $\caM$ l'identifie à $\U(\caM)=\U(2)$. Notons 
$\mathrm{Det}   \left(  \begin{pmatrix}  a &0\\ 0 &\bar a \end{pmatrix} , h'\right) =\frac{a}{\bar a}$.
C'est le déterminant usuel via cette identification.
L'involution $\sigma$ de $G$  dont $H$ est le sous-groupe des points fixes est donnée par la conjugaison sur le premier facteur
 de $\caM^\times \times \caM^\times$ par la matrice $\begin{pmatrix}  1 &0\\ 0 &-1 \end{pmatrix}$.
Les algèbres de Lie $\frg$ et $\frh$ de $G$ et $H$ sont  respectivement :
\[  \frg=   \left\{   (X,X')\in \caM \times \caM, \;  \tr (X)=\tr(X')  \right\}/\sim, \] 
\[ \frh= \left\{   \left(  \begin{pmatrix}  \alpha &0\\ 0 &\bar \alpha \end{pmatrix} , X'\right) , \, \alpha\in \bbC, \, X'\in \caM, \, \tr(X)=\alpha+\bar \alpha \right \}/\sim.\]
L'espace propre $\frs$ pour la valeur propre $-1$ de $\sigma$ dans $\frg$ et le sous-espace de Cartan $\frt$ de $\frs$ sont donnés par 
\[\frs= \left\{   \left(  \begin{pmatrix}  0 &\beta \\ -\bar \beta  &0 \end{pmatrix} , 0 \right) , \, \beta \in \bbC  \right \}/\sim, \quad 
\frt= \left\{   \left(  \begin{pmatrix}  0 &\beta \\ -\beta  & 0 \end{pmatrix} , 0 \right) , \, \beta \in \bbR  \right \} /\sim. \]
On trouve alors 
\[ L  =  \left\{   \left(  \begin{pmatrix}  a &b\\ -b & a \end{pmatrix} , h'\right) , \, a,b \in \bbR, \, h'\in \caM^\times , \, \det (h')=a^2+b^2 \right \}/\sim \simeq \U(2),\] 
et 
\[  L\cap H= \left\{   \left(  \begin{pmatrix}  a &0\\ 0 & a \end{pmatrix} , h'\right) , \, a\in \bbR^\times, \, h'\in \caM^\times , \, \det (h')=a^2 \right \}/\sim \simeq \SU(2).\]

\subsection*{Cas 10 : $\SO(2,2)/\GL(2,\bbR)$}
On pose $\caM=\caM_2(\bbR)$ que l'on munit de la forme quadratique donnée par le déterminant, dont la signature est ici $(2,2)$. 
On définit alors $G$ et son action sur $\caM$ comme dans le cas précédent. Notons que $\caM^\times=\GL(2,\bbR)$.
Ceci identifie $G$ avec $\SO(2,2)$.
Soit $H$   le sous-groupe de $G$ défini par 
\[  H= \left\{   \left(  \begin{pmatrix}  a &0\\ 0 &b \end{pmatrix} , h'\right) , \, a,b\in \bbR^\times, \, h'\in \caM^\times , \, \det (h')=ab \right \}/\sim .  \]
Ce groupe s'identifie de manière évidente à $\GL(2,\bbR)$ et 
$ \mathrm{Det} : H \longrightarrow \bbR^\times ,   \left(  \begin{pmatrix}  a &0\\ 0 &b \end{pmatrix} , h'\right) \mapsto   \frac{b}{a}$
s'identifie au déterminant usuel.
L'involution $\sigma$ de $G$  dont $H$ est le sous-groupe des points fixes est donnée par la conjugaison sur le premier facteur
 de $\caM^\times \times \caM^\times$ par la matrice $\begin{pmatrix}  1 &0\\ 0 &-1 \end{pmatrix}$.
Les algèbres de Lie $\frg$ et $\frh$ de $G$ et $H$ sont  respectivement :
\[  \frg=   \left\{   (X,X')\in \caM \times \caM, \;  \tr (X)=\tr(X')  \right\}/\sim, \] 
\[ \frh= \left\{   \left(  \begin{pmatrix}  \alpha &0\\ 0 & \beta \end{pmatrix} , X'\right) , \, \alpha,\beta \in \bbC, \, X'\in \caM, \, \tr(X)=\alpha+\beta \right \}/\sim\]
L'espace propre $\frs$ pour la valeur propre $-1$ de $\sigma$ dans $\frg$ et le sous-espace de Cartan $\frt$ de $\frs$ sont donnés par 
\[\frs= \left\{   \left(  \begin{pmatrix}  0 &\alpha \\ - \beta  &0 \end{pmatrix} , 0 \right) , \, \alpha, \beta \in \bbR  \right \}/\sim, \quad \frt= \left\{   \left(  \begin{pmatrix}  0 &\beta \\ -\beta  &0 \end{pmatrix} , 0 \right) , \, \beta \in \bbR  \right \} /\sim. \]
 Son centralisateur dans $G$ est le sous-groupe 
\[ L  =  \left\{   \left(  \begin{pmatrix}  a &b\\ -b & a \end{pmatrix} , h'\right) , \, a,b \in \bbR, \, h'\in \caM^\times , \, \det (h')=a^2+b^2 \right \}/\sim \simeq \U(1,1), \]
et 
\[  L\cap H= \left\{   \left(  \begin{pmatrix}  a &0\\ 0 & a \end{pmatrix} , h'\right) , \, a\in \bbR^\times, \, h'\in \caM^\times , \, \det (h')=a^2 \right \}/\sim  \simeq \SU(1,1).\]

 \subsection*{Cas 11 :  $\scrX=G/H=\Sp(4,\bbR)/\Sp(2,\bbR)\times \Sp(2,\bbR)$}
On réalise $\Sp(4,\bbR)$ comme l'ensemble des matrices $M$ dans $\GL(4,\bbR)$
telle que ${}^tMJM=J$, où $J=\begin{pmatrix} 0&1&0&0\\ -1&0&0&0\\0&0&0&-1\\0&0&1&0\end{pmatrix}$.

 L'involution $\sigma$ est la conjugaison par la matrice $\begin{pmatrix} I_2&0_2\\0_2&-I_2\end{pmatrix}$.
Un sous-espace de Cartan est 
\[ \frt_0=\left\{ \tau(a)=  \begin{pmatrix} 0&0&a&0\\ 0&0&0&-a\\-a&0&0&0\\0&a&0&0\end{pmatrix}  , \, a\in \bbR \right\}  \]
et son centralisateur $L$ dans $G$ est l'ensemble des éléments commutant avec la matrice $\tau=\tau(1)$.
On remarque que ${}^t\tau=\tau^{-1}=-\tau$. Les éléments $M$ de $\Sp(4,\bbR)$ commutant avec $\tau$ vérifient
alors ${}^tM(-\tau) J M=(-\tau) J$, et préservent donc la forme bilinéaire symétrique  définie par la matrice
symétrique $-\tau J$. 
La forme symplectique définie par $J$ et la forme  bilinéaire symétrique définie par $-\tau J$ sont respectivement   
\[  ((x,y,z,t) ,(x',y',z',t'))\mapsto  xy'-yx'-zt'+tz', \quad    ((x,y,z,t) ,(x',y',z',t'))\mapsto xt'+yz'+zy'+tx' \] 
Identifions $\bbR^4$ et $\bbC^2$ via
$ (x,y,z,t)\mapsto (x+iz,t+iy)  $.
Les éléments de $L$ préservent alors la forme hermitienne
$     \bil{(x+iz,t+iy)}{(x'+iz',t'+iy')}= (\overline{x+iz}) (t'+iy')  + (\overline{t+iy})(x'+iz')    $
de signature $(1,1)$. On peut remonter les calculs et voir que les éléments de $\Sp(4,\bbR)$ qui préservent cette forme hermitienne
sont exactement ceux de $L$. Ainsi $L\simeq\U(1,1)$ et l'on vérifie que $L\cap H=\SU(1,1)$ dans cette identification.

 \subsection*{Cas 12 : $\scrX=G/H=\Sp(4,\bbR)/ \Sp(2,\bbC)$}
On identifie $\bbC$ et $\bbR^2$ à 
\[  \caE=\left\{    X=\begin{pmatrix}     x_1&x_2\\-x_2&x_1  \end{pmatrix}, x_1,x_2, \in \bbR \right\} \]
et $H=\SL(2,\bbC)$ à l'ensemble des matrices  
$ \begin{pmatrix}  \alpha &\beta  \\ \gamma& \delta  \end{pmatrix}, \quad  \alpha,\beta , \gamma,  \delta \in \caE, \; 
\alpha\delta-\beta\gamma= \begin{pmatrix}     1&0\\0&1  \end{pmatrix}$.  
C'est le groupe de symétrie de la forme symplectique sur $\caE^2$ définie par 
\[ \omega   \left(  \begin{pmatrix}  X\\Y \end{pmatrix},   \begin{pmatrix}  X'\\Y' \end{pmatrix}\right)
=XY'-YX' = \begin{pmatrix}  x_1y'_1-x_2y'_2    -y_1x'_1+y_2x_2'    &    x_2y'_1+x_1y'_2-y_2x'_1-y_1x'_2  \\ 
-  x_2y'_1-x_1y'_2+y_2x'_1+y_1x'_2 &  x_1y'_1-x_2y'_2    -y_1x'_1+y_2x_2'.  \end{pmatrix}\]
Introduisons les deux formes symplectiques réelles sur  $\bbR^4$ données par les parties réelles et imaginaires de $\omega$: 
\[  \omega_1 \left( {}^t(x_1,x_2,y_1,y_2), {}^t (x'_1,x'_2,y'_1,y'_2) \right)=  x_1y'_1-x_2y'_2    -y_1x'_1+y_2x_2',    \]
\[  \omega_2 \left( {}^t(x_1,x_2,y_1,y_2),  {}^t(x'_1,x'_2,y'_1,y'_2) \right) =x_2y'_1+x_1y'_2-y_2x'_1-y_1x'_2 , \]
et les matrices correspondantes 
\[ J_1=\begin{pmatrix}0&0&1&0\\    0&0&0&-1\\ -1&0&0&0\\0&1&0&0  \end{pmatrix} , \qquad J_2=
\begin{pmatrix}0&0&0&1\\    0&0&1&0\\ 0&-1&0&0\\-1&0&0&0  \end{pmatrix} \]
Remarquons que 
\begin{equation}\label{J1J2}
 {}^tJ_1=J_1^{-1}=-J_1, \quad   {}^tJ_2=J_2^{-1}=-J_2, \quad J_1J_2=-J_2J_1. \end{equation}
Soient $\sigma_1$ et $\sigma_2$ les involutions  de $\GL(4,\bbR)$ définies par :
\[ \sigma_1(g)=J_1({}^tg^{-1})J_1^{-1}, \quad  \sigma_2(g)=J_2({}^tg^{-1})J_2^{-1}. \]
Soit $G$ la réalisation de $\Sp(4,\bbR)$ donnée par la première de ces formes, c'est-à-dire
\[ G=\left\{   g\in\GL(4,\bbR) \, \vert\, {}^tgJ_1g=J_1  \right\}  .   \]
C'est le groupe des points fixes de $\sigma_1$. Les relations (\ref{J1J2})
impliquent que les involutions $\sigma_1$, $\sigma_2$, et $\theta$ (l'involution de Cartan $g\mapsto {}^tg^{-1}$) commutent, et $G$ est donc stable par 
$\theta$ et $\sigma_2$.
Soit $H$ le sous-groupe des points fixes de $\sigma_2$ dans $G$. Il s'identifie comme ci-dessus avec $\SL(2,\bbC)$.
 Sur les    algèbres de Lie, on obtient les  décompositions
\[\frg_0=\frk_0\stackrel{\theta}{\oplus} \frp_0, \quad  \frg_0=\frh\stackrel{\sigma_2}{\oplus} \frs_0, \quad
\frg_0=\frk_0\cap \frh_0 {\oplus} \frp_0\cap \frh_0\oplus \frk_0\cap\frs_0\oplus \frp_0\cap \frs_0.  \]
Ici 
\[ \frk_0\cap \frs_0=\frt_0 =\left\{ \begin{pmatrix}0&0&b&0\\    0&0&0&-b\\ -b&0&0&0\\0&b&0&0  \end{pmatrix} , \, b\in \bbR , \right\} \]
est un sous-espace de Cartan.
Son centralisateur dans $G$ est le $c$-Levi 
$  L=\left\{  g\in G\, \vert \, gJ_1=J_1g \right\}$, 
 car  $J_1$ est un élément générique de $\frt_0$, et finalement, on trouve
 \[ L=\{ g\in G\, \vert \, {}^tg^{-1}=g\} =K.\]
 Ainsi $L$ est le compact maximal de $\Sp(4,\bbR)$, donc isomorphe à $\U(2)$, 
  et par cet isomorphisme
 $L\cap H$ devient  $\SU(2)$.

 \subsection*{Cas 13 : $\scrX=G/H=\Sp(2,\bbR)/\GL(1,\bbR)$}  
Un sous-espace de Cartan est 
$ \frt_0=\left\{  \begin{pmatrix} 0&a\\-a&0\end{pmatrix}  , \, a\in \bbR \right\}  $
et  
$  L=\left\{ \begin{pmatrix} \cos \theta &\sin \theta\\-\sin \theta&\cos \theta \end{pmatrix}  , \, \theta\in \bbR \right\}\simeq \U(1), \quad L\cap H=\{  \pm I_2 \}$.

\subsection{Le cas général}\label{rangN}

Le calcul de $L$ en rang quelconque se ramène au rang un par un argument général que nous détaillons seulement dans le cas 1.
 L'involution $\sigma$ est la conjugaison par la matrice diagonale dont les coefficients sont 
 \[  (\underbrace{ 1,-1, 1,-1 \ldots ,1,-1, }_{2p}\underbrace{1,\ldots ,1}_{n-2p}). \]
 On introduit le sous-groupe $G'=\GL(2,\bbR)^{p}\times \GL(n-2p,\bbR)$, stable par l'involution $\sigma$, 
 et l'involution de Cartan $\theta$, de sorte que 
 $H'=G'\cap H=\left( \GL(1,\bbR)\times \GL(1,\bbR)\right)^{p}\times  \GL(n-2p,\bbR)$.
 L'espace symétrique $\scrX'=G'/H'$ est alors isomorphe à $\left( \GL(2,\bbR)/\GL(1,\bbR)\times \GL(1,\bbR)\right)^p$, c'est-à-dire un produit 
 d'espaces de rang un étudiés dans la section précédente. Soit $\frt_0$ un sous-espace de Cartan compact de $\frg'_0$ pour l'espace symétrique
$\scrX'$, et soit $L'$ son centralisateur dans $G'$. 
D'après  les calculs de la section  \ref{Rg1GL}, $L'$ est isomorphe à un 
produit $\left(\bbC^\times \right)^{p} \times \GL(n-2p,\bbR)$.
Par égalité des rangs de $\scrX$ et $\scrX'$ (tous deux égaux à $p$), $\frt_0$ est encore un sous-espace de Cartan compact de $\frg$ 
pour l'espace symétrique $\scrX$, et son centralisateur $L$ dans $G$ contient bien sur  $L'$.
Or les considérations de la remarque \ref{detL} nous donnent  l'isomorphie des complexifiés de $L'$ et $L$  et le fait que $L$ est produit
de groupes généraux linéaires complexes et réels. Ceci   implique alors   l'égalité de $L$ et $L'$, et de plus $L\cap H=L'\cap H=L'\cap H'$.
Pour résumer, on a donc 
\[ L= \left(\bbC^\times \right)^{p} \times \GL(n-2p,\bbR), \quad L\cap H = (\bbR^\times )^p\times \GL(n-2p,\bbR) \]

\medskip

Dans les autres cas, le résultat des calculs est indiqué dans la Table 2.

\begin{rmq}
Dans le cas 3 et 4, pour  l'inclusion des  facteurs $\U(1)$ de $L\cap H$  dans le facteur $\U(1)\times \U(1)$, voir le cas de rang 1.
\end{rmq}



\bibliographystyle{smfalpha}

\bibliography{MR8ver23}

\end{document}